\documentclass[12pt, a4paper]{article}
\usepackage[utf8]{inputenc}
\usepackage{titlesec}
\usepackage[english]{babel}
\usepackage{graphicx}
\usepackage{float}
\usepackage{fancyhdr}
\usepackage[matrix,arrow,curve]{xy}
\usepackage{pstricks} 
\usepackage{amsmath,amsfonts,verbatim,afterpage,theorem,euscript,mathrsfs,amssymb}
\usepackage{amsfonts}
\usepackage{amssymb}
\usepackage{array}
\usepackage{dsfont}
\usepackage{hyperref}
\usepackage{authblk}
\usepackage[titletoc]{appendix}

\newtheorem{Definition}{Definition}[section]
\newtheorem{Proposition}{Proposition}[section]

\newtheorem{Lemma}{Lemma}[section]
\newtheorem{Theorem}{Theorem}
\newtheorem{Corollary}{Corollary}[section]
\newtheorem{Remark}{Remark}[section]

\def \vu{\textbf{u}}
\def \vb{\textbf{b}}
\def \vv{\textbf{v}}
\def \vc{\textbf{c}}
\def \vd{\textbf{d}}
\def \vw{\textbf{w}}

\def \Ri{\mathcal{R}}

\def \X{\mathbb{X}}
\def \vx{\textbf{x}}
\def \vy{\textbf{y}}
\def \vz{\textbf{z}}
\def \vA{\textbf{A}}

\newcommand \te{\theta_\varepsilon}

\newcommand\vN{{\bf  \nabla}}

\def \F{\mathbb{F}}
\def \G{\mathbb{G}}

\def \Rt{\mathbb{R}^{3}}
\def \Rd{\mathbb{R}^{d}}

\def \finpv{\hfill $\blacksquare$  \newline }
\def \pv{{\bf{Proof.}}~}

\def \ds{\displaystyle}
\newcommand \Endproof{\hfill $\diamond$}

\begin{document}

\title{Weak-strong uniqueness in weighted $L^2$ spaces and weak suitable solutions in local Morrey spaces for the  MHD equations}
\author{Pedro Gabriel Fern\'andez-Dalgo\footnote{LaMME, Univ Evry, CNRS, Universit\'e Paris-Saclay, 91025, Evry, France } \footnote{e-mail : pedro.fernandez@univ-evry.fr} , Oscar Jarr\'in \footnote{Direcci\'on de Investigaci\'on y Desarrollo (DIDE), Universidad T\'ecnica de Ambato, Ambato, Ecuador} \footnote{e-mail : or.jarrin@uta.edu.ec}
}
\date{}\maketitle

\maketitle

\begin{abstract} 
We consider here the magneto-hydrodynamics (MHD) equations on the whole space. For the 3D case,  in the setting of the weighted $L^2$ spaces we obtain a weak-strong uniqueness criterion provided that the velocity field and the magnetic field belong to a fairly general multipliers space.  On the other hand, we study the local and global existence of weak suitable solutions for intermittent initial data, which is characterized through a local Morrey space. This large initial data space was also exhibit in a contemporary work \cite{BTK} in the context of 3D Navier-Stokes equations. Finally, we make a discussion on the local and global existence problem in the 2D case. 
\end{abstract}
 
\noindent{\bf Keywords : } MHD equations; Weak-strong uniqueness; Multipliers spaces, Local Morrey spaces;  Global weak solutions; Suitable solutions.\\ 
\noindent{\bf AMS classification : }  35Q30, 76D05.

\section{Introduction}
In a  recent work \cite{PF_PG}, P. G. Fernandez-Dalgo \& P.G. Lemarié-Rieusset obtained \emph{new energy controls} for the    homogeneous and incompressible Navier-Stokes (NS) equations, 
 which allowed them  to develop a theory to construct  weak solutions  for initial data $\vu_0$  belonging to the weighted space $L^{2}_{w_\gamma}(\Rt)=L^2(w_\gamma dx)$, where $0 < \gamma \leq 2$ and $w_\gamma(x):=(1+| x |)^{-\gamma}$. Moreover, this method also gives a \emph{new proof} of the existence of discretely self-similar solutions. \\ 
 
Due to the structural similarity between the (NS) equations and the magneto-hydrodynamics equations (see equations (MHD) below) it is quite natural to extend those recent results obtained for the (NS) equations to the more general setting  of the coupled magneto-hydrodynamics system which writes down as follows: 
\begin{equation*}
(\text{MHD}) \left\{ \begin{array}{ll}\vspace{2mm} 
\partial_t \vu = \Delta \vu - (\vu \cdot \nabla) \vu + (\vb \cdot \nabla) \vb - \nabla p + \vN \cdot \F,  \\ \vspace{2mm}
\partial_{t} \vb = \Delta \vb - (\vu \cdot \nabla) \vb + (\vb \cdot \nabla) \vu -\vN q + 
\vN  \cdot \G 
, \\ \vspace{2mm}
\vN \cdot  \vu=0, \,  \vN \cdot \vb =0, \\ \vspace{2mm}
\vu(0,\cdot)=\vu_0, \, \vb(0,\cdot)=\vb_0.  
\end{array}
\right.     
\end{equation*}
Here the fluid velocity $\vu:[0,+\infty)\times \Rt \to \Rt$, the fluid magnetic field  $\vb:[0,+\infty)\times \Rt \to \Rt$, the fluid pressure $p:[0,+\infty)\times \Rt \to \mathbb{R}$ and the term  $q: [0,+\infty)\times \Rt \to \mathbb{R}$ (which appears in physical models considering Maxwell's displacement currents \cite{Bat}, \cite{Shercliff})  are the unknowns. On the other hand, the data of the problem are given by  the fluid velocity at $t=0$:  $\vu_0:\Rt\to \Rt$; the magnetic field at $t=0$, $\vb_0:\Rt\to \Rt$; and  the tensors  $\F=(F_{i,j})_{1\leq i,j\leq 3}, \G=(G_{i,j})_{1\leq i,j \leq 3}$ (where $F_{i,j}, G_{i,j}: [0,+\infty)\times \Rt \to \mathbb{R}$) whose divergences: $\vN \cdot \F, \vN \cdot \G$, represent  volume forces applied to the fluids.\\

In the setting of this coupled system, in a previous work \cite{PF_OJ}, we adapted the \emph{energy controls} given in \cite{PF_PG} for the (NS) equations to the (MHD) equations and this approach allowed us to establish the existence of discretely self-similar solutions for discretely self-similar initial data belonging to $L^2_{loc}$; and moreover, the existence of global suitable weak solutions when the initial data $\vu_0, \vb_0$ belong to the weighted spaces $L^{2}_{w_\gamma}(\mathbb{R}^3)$, for $0< \gamma \leq 2$, and the tensor forces $\mathbb{F}, \mathbb{G}$ belong to the space $L^2((0, +\infty), L^{2}_{w_\gamma}(\mathbb{R}^3))$.  For all the details see Theorem $1$ and Theorem $2$ in \cite{PF_OJ}. In this paper, we continue with the research program started in \cite{PF_OJ} for the (MHD) equations in two directions. 
\subsection{Weak-strong uniqueness in the weighted $L^2$ spaces}
It is well-known that uniqueness of weak suitable solutions remains an outstanding open problem and in this sense the research community has attired the attention to look for supplementary assumptions in order to ensure the uniqueness of weak solutions. This kind of results are well-known as \emph{weak-strong uniqueness theorems}. In our first result,   we complement the study of the (MHD) equations in the framework of the weighted $L^2$ spaces with a weak-strong uniqueness theorem.  This result, is obtained in the setting of a multiplier space $\mathbb{X}_{T}$ we shall introduce as follows: for a time $0<T<+\infty$ fix, let us denote $E_T$  the \emph{energy space} of the time-dependent vector fields $\vv$ such that $\vv$ belongs to $L^\infty((0,T), L^2)$ and  moreover $\vN\vv$ belongs to $L^2((0,T),L^2)$. This space is doted by its usual norm 
$$\| \vv \|_{E_T}^2 = \sup_{0 \leq t \leq T} \|  \vv(t, .) \|^2_{L^2} + \int_0^T \|  \vN \vv(s, .) \|^2_{L^2} ds.$$
Then, we define  $\X_T$ the space of pointwise multipliers  on $(0,T) \times \mathbb{R}^3$ from $E_T$ to $L^2((0,T),L^2)$, which is a Banach space with  the norm:
$$ \| \vu \|_{\mathbb{X}_T} = \sup_{\| \vv \|_{E_T} \leq 1} \| \vu \vv \|_{L^2((0,T), L^2)}.$$
Moreover, we define  $ \X_T^{(0)}$ the space of multipliers $\vu \in \X_T$ such that for every $t_0 \in [0,T)$ we have 

$$\lim_{t_1 \to t_0^+} \| \mathds{1}_{(t_0, t_1)}(t) \vu(t,\cdot)  \|_{\mathbb{X}_T}= 0.$$

The multiplier space  $\X^{(0)}_{T}$ gives us a natural and general framework to prove a weak-strong uniqueness criterion. More precisely, based on the classical Prodi-Serrin's type condition \cite{GP,S} (for the (NS) equations) and considering the new energy controls on the (MHD) equations in the weighted spaces we have our first main  result. 
\begin{Theorem}[Weak-strong uniqueness]\label{ws} 
	Let $0\leq \gamma\leq 2$. Let $0<T<+\infty$. Let  $\vu_{0}, \vb_{0} \in L^2_{w_\gamma}(\mathbb{R}^3)$ be divergence-free vector fields, and moreover, consider forcing tensors   $\mathbb{F}(t,x)=\left(F_{i,j}(t,x)\right)_{1\leq i,j\leq 3} \in L^2((0,T),L^2_{w_\gamma})$ and $\mathbb{G}(t,x)=\left(G_{i,j}(t,x)\right)_{1\leq i,j\leq 3} \in L^2((0,T),L^2_{w_\gamma})$.\\  
	
	Let $(\vu, \vb, p, q)$ and $(\tilde \vu, \tilde \vb, \tilde p, \tilde q)$ two solutions of the problem (MHD)
	such that :
	\begin{itemize} 
	\item[$\bullet$]  $\vu,\vb, \tilde \vu, \tilde \vb$ belong to the space $L^\infty((0,T), L^2_{w_\gamma})$ and $ \vN \vu, \vN \tilde \vu, \vN\vb, \vN \tilde \vb  $ belong to $L^2((0,T),L^2_{w_\gamma})$
	\item[$\bullet$] the maps $t\in [0,T)\mapsto (\vu, \vb)(t,.)$ and $t\in [0,T)\mapsto (\tilde \vu, \tilde \vb)(t,.)$ are weakly continuous from $[0,T)$ to $L^2_{w_\gamma}(\Rt)$, and are strongly continuous at $t=0$ :
		\begin{equation*}
		\lim_{t\rightarrow 0} \|  (\vu(t,.)-\vu_0, \vb(t,.)-\vb_0) \|_{L^2_{w_\gamma}}=0,
		\end{equation*}
		and
		\begin{equation*}
		\lim_{t\rightarrow 0} \|  (\tilde \vu(t,.)-\vu_0, \tilde \vb(t,.)-\vb_0) \|_{L^2_{w_\gamma}}=0. 
		\end{equation*} 
		\item[$\bullet$]  there exist non-negative locally finite measures  $\mu$ and $\nu$ on $(0,T)\times\mathbb{R}^3$ such that
		\begin{equation}\label{energloc3}
		\begin{split}
		\partial_t(\frac { |\vu  |^2 
			+ |\vb|^2 }2)=&\Delta(\frac { |\vu  |^2 + |\vb|^2 }2)- |\vN \vu  |^2 - |\vN \vb|^2 - \vN\cdot\left( [\frac{ |\vu  |^2}2 + \frac{ |\vb  |^2}2 +p]\vu   \right) \\
		&  + \vN \cdot ([(\vu \cdot \vb)+ q] \vb ) + \vu \cdot(\vN\cdot\mathbb{F}) +\vb \cdot(\vN\cdot\mathbb{G})- \mu,
		\end{split}     
		\end{equation}
		and
		\begin{equation}\label{energloc4}
		\begin{split}
		\partial_t(\frac { | \tilde \vu  |^2 
			+ | \tilde \vb |^2 }2)=&\Delta(\frac { | \tilde \vu  |^2 + | \tilde \vb |^2 }2)- |\vN \tilde \vu  |^2 - |\vN \tilde \vb|^2 - \vN\cdot\left( [\frac{ | \tilde \vu  |^2}2 + \frac{ | \tilde \vb  |^2}2 + \tilde p] \tilde \vu   \right) \\
		&  + \vN \cdot ([( \tilde \vu \cdot \tilde \vb)+ \tilde q] \tilde \vb ) + \tilde \vu \cdot(\vN\cdot\mathbb{F}) + \tilde \vb \cdot(\vN\cdot\mathbb{G})- \nu.
		\end{split}     
		\end{equation}
	\end{itemize}
	If $\vu, \vb \in {\mathbb{X}} _T^{(0)}$ and if the products $\vu \cdot \partial_t \tilde \vu$  and $\vb \cdot \partial_t \tilde \vb$ 
	are well defined as distributions 
	then we have  $( \vu, \vb, p, q ) = ( \tilde \vu, \tilde \vb, \tilde p, \tilde q )$.\\
\end{Theorem} 

Some comments are in order. 
P.G. Lemarié-Rieusset used the space  $\X^{(0)}_T$  to prove a weak-strong uniqueness criterion of weak Leray solutions (see the Theorem $12.4$, page 359 of the book \cite{LR16}). Thus,  we improve  here this  result to the more general setting of the weighted spaces. \\

Let us observe some examples, where we have the embedding $E \subset \X^{(0)}_T$ for the following scale-invariant spaces $E$. For a proof of these embeddings see the Proposition $12.2$, page $361$ of \cite{LR16}. First,  for $2\leq p <+\infty$ such that $2/p+3/q=1$, the classical Prodi-Serrin criterion \cite{GP,S} considers  the space $E=L^p([0,T], L^q(\Rt))$. This result  was improved by P.G. Lemarié-Rieusset to next general spaces:   in \cite{LR02} it is considered $E=L^p([0,T], \mathcal{M}(\dot{H}^{3/q}\to L^2))$, where $\mathcal{M}(\dot{H}^{3/q}\to L^2)$ denotes the pointwise multiplier space from $\dot{H}^{3/q}(\Rt)$ to $L^2(\Rt)$. Thereafter, in \cite{PGLR3}  the multiplier space $\mathcal{M}(\dot{H}^{3/q}\to L^2)$ is substituted by an homogeneous Morrey space and we have $E=L^p([0,T], \dot{M}^{2,q}(\Rt))$. On the other hand, W. Von Wahl considers in \cite{VonWahl} the space $E=\mathcal{C}([0,T], L^3(\Rt))$ and thereafter, this result was generalized by P.G Lemarié-Rieusset in \cite{LR02} to the space $E=\mathcal{C}([0,T], \mathcal{V}^{1}_{0})$, where the space $\mathcal{V}^{1}_{0}$ is defined as the closure of the space $L^3(\Rt)$ in the multiplier spaces from $\dot{H}^{1}(\Rt)$ to $L^2(\Rt)$ denoted as $\mathcal{M}(\dot{H}^{1}\to L^2)$. \\

Finally, we may observe that we also need the assumption that $\vu \cdot \partial_t \tilde \vu$  and $\vb \cdot \partial_t \tilde \vb$ are well-defined in the distributional sense, which essentially is a technical requirement due to the general setting of the space $\X^{(0)}_{T}$ in which we state this result. However, in some particular cases this assumption is not longer required. For example, this is the case  of the space  $L^p((0,T),L^q) \subset \mathbb{X}
^{(0)}_{T}$ (with $2/p+3/q=1$ and $2<p<+\infty$) where the products above are well defined. To illustrate quickly  this fact, let us focus on the term $\vu \cdot \nabla \tilde p$. We remark first  that $ \vN \tilde p \in L^{p'}((0,T),L^{q'}_{\rm loc})$ where $ \frac{1}{p'}=1- \frac{1}{p}$, $ \frac{1}{q'}=1- \frac{1}{q}$. Indeed, let $ \frac{1}{\tilde p} = \frac{1}{p'} - \frac{1}{2} $ and $ \frac{1}{\tilde q} = \frac{1}{q'} - \frac{1}{2} $, by interpolation we find that $\sqrt{w_\gamma} \tilde\vu, \sqrt{w_\gamma} \tilde\vb \in L^{\tilde p} ((0,T), L^{\tilde q})$. Then, for a test function $\varphi \in \mathcal{D}'(\mathbb{R}^3)$, using the continuity of the Riesz transforms, and moreover, assuming that $\F=0$ (only for the sake of simplicity) we have  
\begin{equation*}
\begin{split}
 & \| \varphi  \vN \tilde p \|_{L^{p'} L^{q'}} \leq C_\varphi \| w_\gamma  \vN \tilde p \|_{L^{p'} L^{q'}} \leq C_\varphi \sum_{i,j,k} \| w_\gamma \partial_k (\tilde u _i \tilde u _j)+ \partial_k (\tilde b_i \tilde b_j)\|_{L^{ p'} L^{q '}}\\
  &\leq C (\| \sqrt{w_\gamma} \tilde \vu \|_{L^{\tilde p} L^{\tilde q}} \| \sqrt{w_\gamma} \vN \tilde \vu  \|_{L^2 L^2}+\| \sqrt{w_\gamma} \tilde \vb \|_{L^{\tilde p} L^{\tilde q}} \| \sqrt{w_\gamma} \vN \tilde \vb  \|_{L^2 L^2}).
\end{split}
\end{equation*}
Thus, if $\vu \in L^p((0,T),L^q)$ then we have $\vu \cdot \vN \tilde p \in L^{1}_{loc}([0,T]\times \Rt)$. 

\subsection{Weak suitable solutions for intermittent initial data}
For the (NS) equations,  in a recently  paper \cite{BTK} written by Bradshaw, Tsai \& Kukavika,  the main theorem on global existence given in \cite{PF_PG} by P.G. Fern\'andez-Dalgo \& P.G. Lemarié-Rieusset  is improved with respect to the initial data. 
We \emph{relax} the method developed in \cite{PF_PG}  to  \emph{enlarge} the initial data space and thus we generalize the previous works to the framework of the (MHD) equations. More precisely, following some ideas of \cite{Ba06} (for the (NS) equations) we define local Morrey-type space  $B_{2}(\mathbb{R}^3) \subset L^2_{ \rm loc}(\mathbb{R}^3)$ as the Banach space of all functions $ u \in L^2_{ \rm loc} $ such that :
\begin{equation*}
    \|u\|_{B_2}^2 =  \sup_{R \geq 1}   R^{-2} \int_{|x|\leq R} |u|^{2} \, dx   < +\infty.
\end{equation*}
Moreover, we define the time-space version of the local Morrey-type space:  $B_2 L^2(0,T)$, as   the  Banach space defined as the space of all functions $u  \in L^2_{ \rm loc} ((0,T)\times \mathbb{R}^3) $ such that $$ \| u \|_{B_2 L^2(0,T)}^2 = \sup_{R \geq 1 }  R^{-2}  \int_{|x|\leq R} \int_0^T |u|^{2} \, dt \, dx    < +\infty.$$
In this framework, our second main result reads as follows:    
\begin{Theorem}[Local and global weak suitable solutions]\label{B_2} Let $0<T <+\infty$. Let  $\vu_{0}, \vb_{0} \in B_{2}(\Rt)$ be  divergence-free vector fields. Let    $\mathbb{F}$ and $\G$ be    tensors  belonging to $ B_2 L^2 (0,T)$.  Then, there exists a time $0<T_0 <  T$ such that  the system (MHD) has a solution $(\vu, \vb,p,q)$  which satisfies :
 \begin{itemize} 
 \item[$\bullet$]  $\vu, \vb$ belong to $L^\infty((0,T_0), B_2   )$ and $\vN\vu, \vN\vb $ belong to $B_2   L^2 (0,T_0)$.
\item[$\bullet$] The pressure $p$ and the term $q$ are related to $\vu, \vb$, $\mathbb{F}$ and $\G$ by: 
 $$ p =\sum_{1 \leq i,j\leq 3} \Ri_i \Ri_j(u_iu_j-b_i b_j -F_{i,j})\,\,\,\, \text{and}\,\,\,\,  q =-\sum_{1 \leq i,j\leq 3} \Ri_i \Ri_j(G_{i,j}),$$
 where $\Ri_i= \frac{\partial_i}{\sqrt{-\Delta}}$ denotes the Riesz transform. 
 \item[$\bullet$] The map $t\in [0,T)\mapsto (\vu(t,\cdot ), \vu(t,\cdot ) )$ is $*$-weakly continuous from $[0,T)$ to $B_2(\Rt) $, and for all compact set $K \subset \mathbb{R}^3$ we have:
 $$ \lim_{t\rightarrow 0}    \Vert  ( \vu(t,\cdot) - \vu_{0} , \vb(t,\cdot) - \vb_{0} ) \Vert_{L^{2}(K)} = 0.$$
 \item[$\bullet$] The solution $(\vu, \vb, p, q)$ is suitable : there exists a non-negative locally finite measure $\mu$ on $(0,T)\times\mathbb{R}^3$ such that:
\begin{equation*} 
\begin{split}
   \partial_t(\frac {\vert\vu \vert^2 
    + |\vb|^2 }2)=&\Delta(\frac {\vert\vu \vert^2 + |\vb|^2 }2)-\vert\vN \vu \vert^2 - |\vN \vb|^2 - \vN\cdot\left( [\frac{\vert\vu \vert^2}2 + \frac{\vert\vb \vert^2}2 +p]\vu   \right) \\
  &  + \vN \cdot ([(\vu \cdot \vb)+ q] \vb ) + \vu \cdot(\vN\cdot\mathbb{F}) +\vb \cdot(\vN\cdot\mathbb{G})- \mu.
\end{split}     
\end{equation*}
 \end{itemize}
In particular we have the global control on the solution: for all $0\leq t \leq T_0$, 
 \begin{equation} \label{global-control}
 \begin{split}  
 &\max \{ \|(\vu, \vb )(t) \|_{B_2}^2 , \|\vN(\vu, \vb) \|_{B_2 L^2 (0,T_0)}^2  \} \\
 &\leq  C\| (\vu _0, \vb _0 ) \|^2_{B_{2}} + C\|  (\mathbb{F},\G) \|^2_{B_2 L^2 (0,t)}+ C \int_0^t \| (\vu, \vb) (s) \|^2_{B_{2 }} + \| (\vu, \vb) (s) \|^6_{B_{2 }} ds. \end{split}
\end{equation}
\item[$\bullet$] Finally,  if the data verify: $$
    \lim_{R \to +\infty}  R^{-2} \int_{|x| \leq R} |\vu_0 (x)|^2 + |\vb_0 (x)|^2 \, dx =0,
$$
and 
$$
    \lim_{R \to +\infty}  R^{-2} \int_0^{+\infty} \int_{|x| \leq R} |\mathbb{F}(t,x)|^2 + |\mathbb{G}(t,x)|^2 \, dx \, ds=0,
$$ then $(\vu, \vb, p,q)$ is a global weak solution.
 \end{Theorem}
 
 \begin{Remark}
 A vector field $\vu$ denotes the vector $(u_1, u_2, u_3)$ and for a tensor $\F = (F_{i,j})$ we use $\vN \cdot \F$ to denote the vector $(\sum_i \partial_i F_{i,1}, \sum_i \partial_i F_{i,2}, \sum_i \partial_i F_{i,3})$. Moreover, as $\vN \cdot \vb = 0$ then we can write $ ( \vb \cdot  \vN ) \vu = \vN \cdot (\vb \otimes \vu)$.
 \end{Remark}
 
It is worth to make the following comments on this result. 
Remark first that we prove a global control on the solutions (\ref{global-control}) which is not exhibited in  \cite{BTK}. This new control is also valid for the (NS) equations (taking  $\vb=0, \vb_0=0$ and $\mathbb{G}=0$ in the (MHD) system). On the other hand, it is interesting to note that the main difference between this result and our previous work  \cite{PF_OJ} 
is that, in the more general setting of the space $B_2(\Rt)$, the control on the pressure $p$ and the term $q$ is a little more technical, 
and so the method seems not to be adaptable to study the existence of self-similar solutions of equations (MHD) as done in Theorem $2$ in  \cite{PF_OJ}. \\  

Getting back to the (NS) equations,  the global existence and uniqueness of solutions for the 2D case with initial data $\vu_0 \in B_2  (\mathbb{R}^2)$ was an open problem proposed by A. Basson in \cite{Ba06}. In Appendix \ref{Sec:2D} we make a discussion on problematic arising on  the local and global existence for the 2D case, and moreover, we give a sketch of the proof of a result analogous to Theorem \ref{B_2} in dimension $2$. \\

\textbf{Organization of the paper}.  In Section  \ref{local-Morrey} we state some useful tools on the  weighted spaces and the local Morrey spaces.  Then, in Section \ref{Sec:WS} we give a proof of the weak-strong uniqueness criterion stated in Theorem \ref{ws}. Section \ref{Priori-estim-stability} is devoted to some \emph{a priori} estimates and stability results on the  (MHD) equations. In Section \ref{Proof-Th}, we use these \emph{a priori} results to prove the local and global existence of weak suitable solutions stated in Theorem \ref{B_2}. Finally, Appendix \ref{Sec:2D} is devoted to treat the 2D case. 


\section{The weighted spaces and local Morrey spaces}\label{local-Morrey} 
In order to understand how Theorem \ref{B_2} generalizes the results obtained by \cite{PF_PG}, in this section we state some previous results on the weighted spaces and the local Morrey spaces.  We state these results in the general space  $\mathbb{R}^d$.  
\begin{Definition}
Let $\gamma \geq 0$ and $1 <p < +\infty $. We denote $B^p_{\gamma}(\Rd)$ the Banach space  of all functions $ u \in L^p_{ \rm loc}(\Rd) $ such that :
\begin{equation*}
    \| u  \|_{B_{\gamma}^p} = \sup_{R \geq 1}  \left( \frac{1}{R^{\gamma }} \int_{B(0,R)} |u(x)|^{p} \, dx  \right)^{1/p} < +\infty.
\end{equation*}
Moreover, for $0<T\leq +\infty$,  $B^p_\gamma L^p (0,T)$ is the Banach space of all functions $u \subset (L^p_t L^p_x)_{ \rm loc}([0,T]\times \Rd)$ such that
\begin{equation*}
    \| u \|_{B^p_\gamma L^p (0,T)} = \sup_{R \geq 1}  \left( \frac{1}{R^\gamma} \int_0^T \int_{B(0,R)} |u(t,x)|^p \right)^{\frac{1}{p}} dx \, dt <+\infty.
\end{equation*}
\end{Definition}
In what follows,   we will denote $B^{p}_{\gamma}(\Rd)= B^{p}_{\gamma}$ and $B^{2}_2= B_2$.\\

Also, the space $B^p_{\gamma,0 }$ is defined as the subspace of all functions $u \in B^p_{\gamma }$ such that 
$$ \ds{\lim_{ R \to +\infty }  \frac{1}{ R^{\gamma }} \int_{B(0,R)} | u (x) |^p \, dx =0},$$ and similar,  $ B^p_{\gamma,0} L^p (0,T)$ is the subspace of all functions $u \in B^p_{\gamma} L^p (0,T)$ such that $$\ds{  \lim_{ R \to +\infty }  \frac{1}{ R^{\gamma }}  \int_0^T \int_{B(0,R)} | u (t,x) |^p \, dx\, dt =0. }$$ \\%

The following result shows how $B^p_{\gamma }$ is strongly lied with the weighted spaces $ L^p_{w_\gamma}=L^p(w_\gamma \, dx)$ (where  $w_\gamma=(1+\vert x\vert)^{-\gamma}$) considered in  \cite{PF_OJ} and \cite{PF_PG}.
\begin{Lemma}
\label{li}
Consider $\gamma \geq 0$ and let $  \gamma < \delta < + \infty $. We have the continuous embedding  $$L^p_{w_{\gamma}} \subset B^p_{\gamma,0 } \subset B^p_{\gamma } \subset L^p_{w_{\delta}}. $$ Moreover, for all $0<T\leq +\infty$ we have:   
 $$ L^p ((0, T), L^p_{w_\gamma}) \subset B^p_{\gamma,0} L^p (0,T) \subset B^p_\gamma L^p (0,T) \subset L^p ((0, T), L^p_{w_\delta}).$$  
\end{Lemma}
\pv
 Only the embedding $L^p ((0, T), L^p_{w_\gamma}) \subset B^p_{\gamma,0} L^p (0,T)$ is not proved in \cite{PF_PG2} and we prove it. Let $\lambda > 1$ and  $n\in \mathbb{N}$,  let $u_n(t,x) = u(t, \lambda^n x) $. We have:
 \begin{align*}
     \sup_{R \geq 1}  &\frac{1}{ (\lambda^n R)^\gamma} \int_0^{T} \int_{|x|\leq \lambda^n R} |u (t,  x)|^p \, dx \, dt = \sup_{R \geq 1} \frac{\lambda^{(d-\gamma)n}}{ R^\gamma} \int_0^{T} \int_{|x|\leq R} | u ( t, \lambda^n x)|^p \, dx \, dt \\
      = & \lambda^{(d-\gamma)n}\|u_{n}\|_{B^p_{\gamma} L^p (0,T)}^p 
       \leq C \lambda^{(d-\gamma)n} \|u_{n}\|_{L^p L^p_{w_\gamma} }^p \leq C \int_0^{T}\int \vert u(s,x)\vert^p  \frac{1}{(\lambda^n+\vert x\vert)^\gamma}  \, dx\, dt,
 \end{align*}
 and  we conclude by dominated convergence.\Endproof{}\\

Thereafter, we have the following result involving the interpolation theory of Banach spaces: 
\begin{Theorem}[\cite{PF_PG2}]\label{ti}
The space $B^p_\gamma$ can be obtained by interpolation: for all $0 < \gamma < \delta < \infty$ we have $\ds{  B^p_\gamma = [L^p, L^p_{w_{\delta}} ]_{\frac{\gamma}{\delta}, \infty}}$;  and the norms
$\| \cdot \|_{B^p_\gamma} $ and $ \| \cdot \|_{[L^p, L^p_{w_{\delta}} ]_{\frac{\gamma}{\delta}, \infty}}$ are equivalents.
\end{Theorem}
This theorem has a useful  corollary and in order to state it we recall first  the following result on the Muckenhoupt weights (see  \cite{Gr09} for a definition).

\begin{Lemma}[Muckenhoupt weights, \cite{PF_PG}]\label{Lemmuck} If $0<\delta<d$ and $1<p<+\infty$. Then, $w_\delta(x)=(1+\vert x \vert)^{-\delta}$ belongs to the Muckenhoupt class $\mathcal{A}_p(\mathbb{R}^d)$. Moreover we have: 
 \begin{itemize} 
    \item The Riesz transforms $R_j$ are bounded on $L^p_{w_\gamma} :
    \|R_jf\|_{L^p_{w_\gamma}}\leq C_{p,\delta} \|f\|_{L^p_{w_\gamma}} $
    
    \item The Hardy--Littlewood maximal function operator is bounded on $L^p_{w_\gamma}$ :
    \begin{equation*}
        \|\mathcal{M}_f\|_{L^p_{w_\gamma}}\leq C_{p,\delta} \|f\|_{L^p_{w_\gamma}}.
    \end{equation*}
 \end{itemize}
 \end{Lemma}

With this lemma at hand, the next important corollary of Theorem \ref{ti} follows:   

\begin{Corollary}
\label{operonb}
 If $0<\delta<d$ and $1<p<+\infty$, then we have:
 \begin{itemize} 
     \item The Riesz transforms $R_j$ are bounded on $B^p_{\delta}$ :
    $ \|R_jf\|_{B^p_{\delta}}\leq C_{p,\delta} \|f\|_{B^p_{\delta}} 
    $
    \item The Hardy--Littlewood maximal function operator is bounded on $B^p_\delta$ :
    \begin{equation*}
        \|\mathcal{M}_f\|_{B^p_{\delta}}\leq C_{p,\delta} \|f\|_{B^p_{\delta}}.
    \end{equation*}
 \end{itemize}
\end{Corollary}

\pv Remark that  Theorem  \ref{ti} implies $B^p_{\delta} = [L^p, L^p_{w_{\delta_0}} ]_{\frac{\delta}{\delta_0}, \infty}$, for some $\delta < \delta_0 < d$. So, we conclude directly by Lemma   \ref{Lemmuck}.
\Endproof{}\\

To close this section, let us recall the following useful Sobolev embedding for the weighted spaces. For a proof, see the Lemma $3$ of \cite{PF_PG}.
\begin{Lemma} Let $\gamma \geq 0$. Let $f\in L^2_{w_\gamma}$ such that $\vN f\in L^2_{w_\gamma}$ then $f\in L^6_{w_{3\gamma}}$ and
\begin{equation}\label{sobol}
 \|f\|_{L^6_{w_{3\gamma}}}\leq C_\delta (\|f\|_{L^2_{w_\gamma}}+ \|\vN f\|_{L^2_{w_\gamma}}).   
\end{equation}
\end{Lemma}
\section{Weak-strong uniqueness in weighted spaces}\label{Sec:WS}
We start by the following characterization of the pressure terms that we shall need later to prove our result on the weak-strong uniqueness. 

\subsection{Remark on the  pressure terms}

We may observe that the recent results and their proofs given in \cite{PF_PG2}, about the characterization of the pressure term in the (NS) equations can be immediately adapted to the framework of the coupled  (MHD) system. Here, we adapt the proof of the main result in $\cite{PF_PG2}$ to the (MHD) system.

\begin{Proposition}
  \label{prp}
 Let $d \in \{ 2 , 3\}$ be the spatial dimension. Let $0<T<+\infty$. Let   $\mathbb{F}, \G$ be  tensors $\F(t,x)=\left(F_{i,j}(t,x)\right)_{1\leq i,j\leq d}$, $\G(t,x)=\left(G_{i,j}(t,x)\right)_{1\leq i,j\leq d}$ such that $\F, \G$ belongs to $L^1((0,T),L^{1}_{w_{d+1}}(\Rd))$.
  Let $(\vu, \vb, {\bf S}, {\bf R} )$ be a solution of the following problem 
 \begin{equation*} \left\{ \begin{array}{ll}\vspace{2mm}
\partial_t \vu = \Delta \vu - \vN \cdot (\vu \otimes \vu) + \vN \cdot (\vb \otimes \vb) - {\bf S} + \vN \cdot \F,  \\ \vspace{2mm}
\partial_{t} \vb = \Delta \vb - \vN \cdot ( \vu \otimes \vb ) + \vN \cdot (\vb \otimes \vu) -{\bf R} + 
\vN  \cdot \G 
, \\ \vspace{2mm}
\vN \cdot  \vu=0, \quad  \vN \cdot \vb =0 \quad \vN\wedge {\bf R}=0 \quad \vN\wedge {\bf S}=0, 
 \\ \vspace{2mm}
  \vu(0,x)=\vu_0 (x), \quad  \vN \vb(0,x) =\vb_0 (x),
\end{array}
\right.     
\end{equation*}
such that :
  $\vu, \vb$ belong to $L^2 ((0,T),L^2_{w_{d+1}}(\Rd))$, ${\bf R, S}$ belong to $\mathcal{D}'( (0,T) \times \mathbb{R} ^d )$,  $ \lim_{t \to 0} \vu(t, .) = \vu_0 $ in $\mathcal{D}'(\mathbb{R}^d)$ and also $ \lim_{t \to 0} \vb(t, .) = \vb_0 $ . Moreover, let   $\varphi \in \mathcal{D}(\Rd)$ be an arbitrary test function  such that $\varphi(x )=1$ on a neighborhood of the origin. We define  the function  $$\Phi_{i,j,\varphi}=(1-\varphi) \partial_i \partial_j G_d,$$ where $G_d$ denotes the fundamental solution of the  operator $-\Delta$. \\ 
  
  Then,   there exist two functions $f(t), g(t)\in L^1((0,T))$ such that
   $${\bf S}=\vN p_\varphi +\partial_t g \quad and \quad {\bf R}=\vN q_\varphi +\partial_t f $$ with \begin{align*} p_\varphi= &  \sum_{i,j}(\varphi \partial_i \partial_j G_d) * (u_iu_j - b_i b_j -F_{i,j}) \cr &+   \sum_{i,j}   \int (\Phi_{i,j,\varphi}(x-y)-\Phi_{i,j,\varphi}(-y))  (u_i u_j - b_i b_j -F_{i,j})(t,y)\, dy \end{align*}
   and
   \begin{align*} q_\varphi= &  \sum_{i,j}(\varphi \partial_i \partial_j G_d) * (-G_{i,j}) \cr &+   \sum_{i,j}   \int (\Phi_{i,j,\varphi}(x-y)-\Phi_{i,j,\varphi}(-y))  (-G_{i,j}(t,y))\, dy. \end{align*}
   
   Moreover, the following facts hold:
   \begin{itemize}
   \item[$\bullet$] $\vN p_\varphi, \vN q_\varphi$ do  not depend on the choice of $\varphi$ : if we change $\varphi$ in $\psi$, then we have
   $$ p_\varphi(t,x)-p_\psi(t,x)=  \sum_{i,j}   \int (\Phi_{i,j,\psi}(-y)-\Phi_{i,j,\varphi}(-y))  (u_iu_j - b_i b_j-F_{i,j})(t,y)\, dy $$ and
   $$ q_\varphi(t,x)-q_\psi(t,x)=  \sum_{i,j}   \int (\Phi_{i,j,\psi}(-y)-\Phi_{i,j,\varphi}(-y))  (-G_{i,j}(t,y))\, dy $$
   \item[$\bullet$] $\vN p_\varphi$ is the unique solution of the Poisson problem
   $$ \Delta \vw= - \vN(\vN\cdot (\vN\cdot (\vu\otimes\vu - \vb\otimes\vb -\mathbb{F})) $$
   with 
   $$ \lim_{\tau\rightarrow +\infty} e^{\tau\Delta}\vw =0 \text{ in }\mathcal{D}',$$
   and
   $\vN q_\varphi$ is the unique solution of the Poisson problem
   $$ \Delta \vw= - \vN(\vN\cdot (\vN\cdot (-\G)) $$
   with 
   $$ \lim_{\tau\rightarrow +\infty} e^{\tau\Delta}\vw =0 \text{ in }\mathcal{D}'.$$
   \item If $\mathbb{F}$ belongs  to  $L^1 ((0,T),L^1_{w_{d}}(\Rd))$
and $\vu, \vb$ belong to $L^2 ((0,T),L^2_{w_{d}}(\Rd))$, then $g=0$ and $\vN p_\varphi=\vN p_0$, where we have\footnote{Remark that actually $p_0$ does not  depend on $\varphi$ and could have been defined as 
	$\ds{p_0=   \sum_{i,j}(  \partial_i \partial_j G_d) * (u_iu_j + b_i b_j -F_{i,j})}$.}
  \begin{equation*}
    p_0=   \sum_{i,j}(\varphi \partial_i \partial_j G_d) * (u_iu_j - b_ib_j -F_{i,j}) +   \sum_{i,j}((1-\varphi) \partial_i \partial_j G_d) * (u_iu_j - b_ib_j-F_{i,j}).
  \end{equation*}
  \item  If $\mathbb{G}$ belongs to  $L^1 ((0,T),L^1_{w_{d}}(\Rd))$ 
  and if $\vu, \vb$ belong to $L^2 ((0,T),L^2_{w_{d}}(\Rd))$, then $f=0$ and $\vN q_\varphi=\vN q_0$, where
  \begin{equation*}
    q_0=   \sum_{i,j}(\varphi \partial_i \partial_j G_d) * (-G_{i,j}) +   \sum_{i,j}((1-\varphi) \partial_i \partial_j G_d) * (-G_{i,j}).
  \end{equation*}
   \end{itemize}
      \end{Proposition}
      
\pv{}  First, we take the  divergence operator in the equation
  \begin{equation*}
      \partial_t  \vu = \Delta \vu - \vN \cdot (\vu\otimes\vu)+\vN\cdot(\vb \otimes \vb) -{\bf S} + \vN \cdot \F,
  \end{equation*}
and since we have  $div(\vu)=div(\vb)=0$, then we obtain the identity
  \begin{equation*}
      -\sum_{i,j} (\partial_i \partial_j (u_i u_j + b_i b_j)+  \partial_i \partial_j F_{i,j}) - \vN \cdot {\bf S}  =0,
  \end{equation*}
hence we have
  \begin{equation*}
      - \Delta {\bf S}  = \vN \left(\sum_{i,j} \partial_i \partial_j (u_i u_j - b_i b_j- F_{i,j} )  \right). 
  \end{equation*}
Following the same computations for the second equation in the system above we also get the identity  
  \begin{equation*}
      - \Delta {\bf R}  = \vN \left(\sum_{i,j} \partial_i \partial_j (- G_{i,j} ) \right). 
  \end{equation*}
On the other hand, let us set the functions $h_{i,j} = u_i u_j - b_i b_j - F_{i,j}$, and    
  $\Phi_{i,j,\varphi} = (1-\varphi) \partial_i \partial_j G_d $. Then, by Proposition 3 in \cite{PF_PG2} we can define the functions $p_\varphi$ and $q_\varphi$ as follows:
\begin{align*}
p_\varphi=  \sum_{i,j}(\varphi \partial_i \partial_j G_d) * h_{i,j} + \sum_{i,j} \int (\Phi_{i,j,\varphi}(x-y) - \Phi_{i,j\varphi}(-y) )  h_{i,j}(y) dy,
\end{align*}
and
\begin{align*}
    q_\varphi=  \sum_{i,j}(\varphi \partial_i \partial_j G_d) * (-G_{i,j}) + \sum_{i,j} \int ( \Phi_{i,j,\varphi}(x-y) - \Phi_{i,j\varphi}(-y) )  (-G_{i,j})(y) dy.
\end{align*}
Taking the gradient operator in these expressions we write  
\begin{equation*}
  \vN p_\varphi=  \vN \left(\sum_{i,j} (\varphi \partial_i \partial_j G_d) * h_{i,j}\right) + \vN\left( \sum_{i,j}((1-\varphi) \partial_i \partial_j G_d) * h_{i,j}\right)=U_1+U_2,
\end{equation*}
and 
\begin{equation*}
\begin{split}
\vN q_\varphi &=  \vN\left( \sum_{i,j}(\varphi \partial_i \partial_j G_d) * (-G_{i,j})\right) + \vN\left( \sum_{i,j}((1-\varphi) \partial_i \partial_j G_d) * (-G_{i,j})\right)\\
&=V_1+V_2.
\end{split}
\end{equation*}
Now, let us set  $\tilde S= {\bf S} -\vN p_\varphi$ and $\tilde R= {\bf R} -\vN q_\varphi$. We remark first  that we have the identities $\Delta (\vN p_\varphi) = \Delta {\bf S} $ and $\Delta (\vN q_\varphi) = \Delta {\bf R}$,  hence we obtain   $\Delta \tilde S = 0$ and $\Delta \tilde R = 0$ and  conclude that   $\tilde S , \tilde R$ are  harmonic  in $\mathcal{D}'(\Rd)$. Now, we will prove that $\tilde S$ and $\tilde R$ belong to the space $\mathcal{S}'(\Rd)$. For this, let $\alpha \in \mathcal{D}(\mathbb{R})$ be such that $\alpha (t)= 0$  for all $|t| \geq \varepsilon$, and moreover, let  $\beta\in\mathcal{D}(\mathbb{R}^d)$. Then, denoting by $*=*_{t,x}$ the convolution in the temporal and spatial variables, for $t \in (\varepsilon, T-\varepsilon)$  we write
 \begin{equation*} \begin{split} \tilde S(t) *_{t,x} &(\alpha\otimes\beta)  = ( \vu  *(-\partial_t\alpha\otimes\beta+\alpha\otimes\Delta\beta)\\  &+(-\vu\otimes\vu + \vb\otimes\vb + \F) \cdot *(\alpha\otimes\vN\beta))(t,\cdot )\\&  - \sum_{i,j} ( (h_{ij}) * ( \vN (\varphi \partial_i \partial_j G_d )*( \alpha\otimes\beta))) (t,\cdot) -  ( U_2 * (\alpha\otimes\beta)) (t,\cdot),\end{split} \end{equation*} 
 and
 \begin{equation*} \begin{split} \tilde R(t) *_{t,x} & (\alpha\otimes\beta)  = ( \vb *(-\partial_t\alpha\otimes\beta+\alpha\otimes\Delta\beta) \\
 &+ (-\vu \otimes \vb + \vb \otimes \vu + \G) \cdot *(\alpha\otimes\vN\beta))(t,\cdot )\\&  - \sum_{i,j} ( (-G_{ij}) * ( \vN (\varphi \partial_i \partial_j G_d )*( \alpha\otimes\beta))) (t,\cdot) -  ( V_2 * (\alpha\otimes\beta)) (t,\cdot).\end{split} \end{equation*} 
At this point, we may apply the  Proposition $1$ in \cite{PF_PG2} to  conclude that $\tilde S *( \alpha\otimes\beta)(t,.)$ and $\tilde R *( \alpha\otimes\beta)(t,.)$ belong to the space  $L^{1}_{w_{d+1}}(\Rd)$. As a consequence of this fact we have $\tilde S *( \alpha\otimes\beta)(t,.), \tilde R *( \alpha\otimes\beta)(t,.) \in \mathcal{S}'(\Rd)$ and  recalling that  they are harmonic we get that $\tilde S *( \alpha\otimes\beta)(t,.)$ and $\tilde R *( \alpha\otimes\beta)(t,.)$ are polynomials. But, as these functions belong to $L^{1}_{w_{d+1}}(\Rd)$ they are constants functions.\\

Now, if we assume that $\mathbb{F}$ belongs  to  $L^1 ((0,T),L^1_{w_{d}}(\Rd))$
and moreover  $\vu, \vb$ belong to $L^2 ((0,T),L^2_{w_{d}}(\Rd))$, then we find that the polynomial $\tilde S *( \alpha\otimes\beta)(t,.)$ belongs to  $L^{1}_{w_d}(\Rd)$, and thus it  is equals to $0$. Similarly, by the assumption that $\mathbb{G}$ belongs to  $L^1 ((0,T),L^1_{w_{d}}(\Rd))$
and $\vu, \vb$ belong to $L^2 ((0,T),L^2_{w_{d}}(\Rd))$, we get  that the corresponding polynomial $\tilde R *( \alpha\otimes\beta)(t,.)$ is equals to $0$. So we can use the identity approximation $\phi_\varepsilon = \frac 1{\varepsilon^4} \alpha	(\frac t\varepsilon)\beta(\frac x\varepsilon)$ and taking the limit when $\epsilon\to 0$,  we obtain that  $\tilde S$ and $\tilde R$ are constants in the spatial variable.  Thus, we have  ${\bf S}=\vN p_\varphi +a(t)$, where $a(t) =0$ under the assumptions $\mathbb{F} \in L^1 ((0,T),L^1_{w_{d}}(\Rd))$
and $\vu, \vb \in L^2 ((0,T),L^2_{w_{d}}(\Rd))$, and  ${\bf R} =\vN q_\varphi +b(t)$ where $b(t) =0$ under the assumptions $\mathbb{G} \in L^1 ((0,T),L^1_{w_{d}}(\Rd))$
and $\vu, \vb \in L^2 ((0,T),L^2_{w_{d}}(\Rd))$.\\ 

We remark now that as the functions $a$ and $b$ above do not depend the spatial variable then we may take a function $\beta\in\mathcal{D}(\Rd)$ such that $\int\beta\, dx=1$ and write $a=a*_x\beta$. Thus, for $t_0$  a Lebesgue point of the functions $t  \mapsto \| \vu (t, .) \|_{L^2(w_{d+1} dx)}$ and $t  \mapsto \| \vb (t, .) \|_{L^2(w_{d+1} dx)}$, we get the identities 
$$ a(t)=\partial_t (\vu_{t_0}*\beta-\vu*\beta+\int_{t_0}^t \vu*\Delta\beta -(\vu\otimes\vu - \vb\otimes\vb -\mathbb{F})\cdot*\vN\beta-p_\varphi* \vN\beta \, ds)=\partial_t g,$$
and
$$ b(t)=\partial_t (\vb_{t_0}*\beta-\vb*\beta+\int_{t_0}^t \vb*\Delta\beta -(\vu \otimes \vb - \vb \otimes \vu -\G)\cdot*\vN\beta-q_\varphi* \vN\beta \, ds)=\partial_t f.$$
Here, we observe that as we have $\partial_t\partial_j g=\partial_j a=0$ and $\partial_j g(0,.)=0$, then  we find that $g$ depends only on temporal variable. Moreover, the expression of  $g$ given above  proves that $g\in L^1((0,T))$. Using the same argument we also get  $f\in L^1((0,T))$.\\

\Endproof

\subsection{Proof of Theorem \ref{ws}}
By the characterization of the pressure terms given above (with the dimension $d=3$), we know that $p$ and $q$ can be taken as follows, let us choose   $\varphi \in \mathcal{D}(\Rt)$   such that $\varphi(x )=1$ on a neighborhood of the origin  and let us denote $h_{i,j} = u_i u_j - b_i b_j - F_{i,j}$, and 
  $\Phi_{i,j,\varphi} = (1-\varphi) \partial_i \partial_j G_3 $, then we can take
\begin{align*}
    p=  \sum_{i,j}(\varphi \partial_i \partial_j G_3) * h_{i,j} + \sum_{i,j} \int ( \Phi_{i,j,\varphi}(x-y) - \Phi_{i,j\varphi}(-y) )  h_{i,j}(y) dy,
\end{align*}
and
\begin{align*}
    q=  \sum_{i,j}(\varphi \partial_i \partial_j G_3) * (-G_{i,j}) + \sum_{i,j} \int (\Phi_{i,j,\varphi}(x-y) - \Phi_{i,j\varphi}(-y) )  (-G_{i,j})(y) dy.
\end{align*}
In fact, let us emphasize that  we can characterize the pressure terms in a almost general context only using the Riesz transforms, we refer to \cite{BA_PA_PF} for more details. However, the general characterization given in the Proposition \ref{prp} is very useful to study the (MHD) equations with initial data in $B_2$ (see Appendix).

For instance, as $\sqrt{w_\gamma}\vu,\sqrt{w_\gamma}\vb \in L^2((0,T), L^2)$ and $\sqrt{w_\gamma} \vN \vu,\sqrt{w_\gamma} \vN \vb \in L^2((0,T), L^2)$, we obtain by interpolation that $w_\gamma u_i u_j $ and $w_\gamma b_i b_j$ belongs to $L^{\hat a} L^{\hat b}$ with $\frac{2}{\hat a}+ \frac{3}{\hat b}= \frac{3}{2}$.  Taking $r \in \left(1, \min \{ \frac{3}{2}, \frac{3}{\gamma} \}\right)$, 
and $a \in \mathbb{R}$ satisfying $\frac{2}{\hat a} + \frac{3}{r} = 3$, we get that 
\begin{equation*}
   \mathcal{R}_i \mathcal{R}_j (u_i u_j),  \mathcal{R}_i \mathcal{R}_j (b_i b_j)  
  \in L^{a}((0,T),L^{r}_{w_{ r \gamma}}(\mathbb{R}^d)), 
\end{equation*}
and  by the continuity of the Riesz transforms on $L^2_{w_{ \gamma}}(\mathbb{R}^3)$ we have
\begin{equation*}
   \mathcal{R}_i \mathcal{R}_j F_{i,j},  \mathcal{R}_i \mathcal{R}_j G_{i,j} 
  \in L^{2}((0,T),L^{2}_{w_\gamma}(\mathbb{R}^d)).
\end{equation*}
Indeed, the following estimate holds: taking $\hat b$ given by $\frac{2}{\hat a}+ \frac{3}{\hat b}= \frac{3}{2}$, we can write
\begin{align}  \label{Eq:2}
  & \Big{\|} \,  \mathcal{R}_i \mathcal{R}_j (u_i u_j)  
    \, \Big{\|}_{L^{\hat a}((0,T), L^{r}_{w_{r \gamma}}(\mathbb{R}^3))}
    \nonumber \\
  & \leq C_{ \gamma, r} 
     \| u_i u_j \|_{L^{\hat a}((0,T), L^{r}_{w_{r \gamma}}(\mathbb{R}^3))} 
    \nonumber \\
  & \leq C_{ \gamma, r} 
    \|\sqrt{w_\gamma} \, u_i \|_{L^{\infty} \left((0,T), L^{2}(\mathbb{R}^3) \right)} 
    \, \cdot \, 
    \| \sqrt{w_\gamma} \, u_j \|_{L^{\hat a} \left((0,T), L^{\hat b}(\mathbb{R}^3) \right)}
    \nonumber \\
  & \leq \gamma^{\frac{1}{\hat a}} \; {\tilde{C}}_{\gamma, r}  \; 
    \| {\textbf{u}} \|_{L^{\infty}((0,T), L^{2}_{w_\gamma}(\mathbb{R}^3))}^{1+\frac{\hat a-2}{ \hat a}} 
    \nonumber \\
  & \;\;\;\; \times     
    \left( \int_0^T ( \| {\textbf{u}}(s) \|_{L^2_{w_\gamma}(\mathbb{R}^3)} 
    + \| \nabla {\textbf{u}}(s) \|_{L^2_{w_\gamma}(\mathbb{R}^3)})^2 ds 
    \right)^{\frac{1}{\hat a}}. 
\end{align}
Thus, as the Riesz transforms are well-defined for all the terms  composing the pressure terms, we have necessarily  the identities
\begin{equation}\label{Charact-p}
     p = \sum_{1 \leq i,j\leq 3} \Ri_i \Ri_j(u_iu_j-b_i b_j -F_{i,j}),
\end{equation}
and
\begin{equation*}
    q = \nabla \left(\sum_{1 \leq i,j\leq 3} \Ri_i \Ri_j( -G_{i,j}) \right).
\end{equation*}
The same identities hold true for the pressure terms $\tilde p$ and $\tilde q$.\\

Now, let $\vv = \vu - \tilde \vu$, $\vw= \vb - \tilde \vb$, $a=p-\tilde p$ and $b= q - \tilde q$. We remark first that by the characterization of the terms $q$ and $\tilde{q}$ 
we have $q=\tilde{q}$ and then $b=0$. So, we will prove the identities $\vv=0$, $\vw=0$ and $a=0$. \\

Using the  identity 
\begin{equation}\label{IdBase}
\frac{\vert \vv \vert^2+ \vert \vw \vert^2}{2}=\frac{\vert \vu \vert^2 +\vert \vb \vert^2}{2}+ \frac{\vert \tilde{\vu} \vert^2+ \vert \tilde{\vb} \vert^2}{2}-\vu \cdot \tilde{\vu}-\vb \cdot \tilde{\vb},
\end{equation} 
we   write   
\begin{equation}\label{eq-aux}
\begin{split}
\partial_t (\frac { | \vv  |^2 + | \vw  |^2}2) =& \partial_t (\frac { | \vu  |^2 + | \vb  |^2}2) + \partial_t (\frac { | \tilde \vu  |^2 + | \tilde \vb  |^2}2)  \\
& - \vu \cdot \partial_t  \tilde \vu - \tilde \vu \cdot \partial_t  \vu - \vb \cdot \partial_t \tilde \vb - \tilde \vb \cdot \partial_t  \vb.
\end{split}
\end{equation}
Recall that by assumption the terms $\vu \cdot \partial_t  \tilde \vu$ and $\vb \cdot \partial_t \tilde \vb$ are well-defined as distributions so it remains to must verify that the terms $\tilde \vu \cdot \partial_t  \vu$ and $\tilde \vb \cdot \partial_t  \vb$ are also well-defined   in the distributional sense.  For this we have the following technical lemma.

\begin{Lemma} Within the framework of Theorem \ref{ws}, as we have $ \vu,  \vb  \in \X_T$ then we get $\tilde \vu \cdot \partial_t \vu\in L^{1}_{loc}([0,T]\times \Rt)$ and $\tilde \vb \cdot \partial_t \vb\in L^{1}_{loc}([0,T]\times \Rt)$. 
\end{Lemma}	
\pv We shall verify that we have $\tilde \vu \cdot \partial_t \vu\in L^{1}_{loc}([0,T]\times \Rt)$. The treatment for the other term $\tilde \vb \cdot \partial_t$ follows the same lines. \\

As we have $\ds{\partial_t  \vu= \Delta  \vu - ( \vu  \cdot \nabla)  \vu + (\vb \cdot \nabla)  \vb -\nabla  p + \nabla \cdot \F}$, then we formally write 
\begin{equation}\label{EQ}
 \tilde \vu \cdot \partial_t  \vu= \tilde \vu \cdot \Delta  \vu - \tilde\vu \cdot (( \vu \cdot \nabla)  \vu) + \tilde\vu \cdot (( \vb \cdot \nabla)  \vb) -\tilde \vu \cdot \nabla  p + \tilde\vu \cdot (\nabla \cdot \F),
\end{equation}
where we must prove that each term in the right side belong to $L^{1}_{loc}([0,T] \times \Rt)$.  Since all the terms are  treated in a similar way it is enough to detail the computations for the terms $ \tilde\vu \cdot (( \vu \cdot \nabla)  \vu) $ and $\tilde \vu \cdot \nabla  p $.  Let $0<t\leq T$ and let $\varphi \in \mathcal{D}(\Rt)$ be an arbitrary  test function. For the function $\varphi$ given, we set  $\psi \in \mathcal{D}(\Rt)$ such that $0 \leq \psi \leq 1$ and $\psi=1$ on $supp(\varphi)$. \\

For the term $\tilde \vu \cdot ((\vu \cdot \nabla) \vu)$,  since $div(\vu)=0$ then $(\vu \cdot \nabla) \vu = \vN \cdot (\vu \otimes \vu)$ and we have 
\begin{equation*}
\begin{split}
\int_{0}^{t}\int \varphi \tilde \vu \cdot (( \vu \cdot \nabla)  \vu)\, dx \, ds &= \sum_{i,j}\int_{0}^{t}\int \varphi \tilde u_i \partial_j ( u_j  u_i)\, dx \, ds\\
& = \sum_{i,j}\int_{0}^{t}\int \varphi \tilde  u_i \partial_j (\psi \,  u_j  u_i)\, dx \, ds,
\end{split}
\end{equation*}hence we can write 
\begin{equation}\label{Estim-aux}
\begin{split}
&\left\vert  \int_{0}^{t}\int \varphi \tilde  \vu \cdot (( \vu  \cdot \nabla)  \vu)\, dx \, ds\right\vert \leq C \sum_{i,j} \int_{0}^{t} \Vert \varphi \tilde  u_i \Vert_{\dot{H}^{1}} \Vert \partial_j(  \psi  u_j   u_i) \Vert_{\dot{H}^{-1}}\, ds   \\
\leq & C\int_{0}^{t} \Vert \nabla (\varphi \tilde  \vu) \Vert_{L^2} \Vert \psi ( \vu \otimes  \vu ) \Vert_{L^2}\, ds \leq \Vert \nabla (\varphi \tilde  \vu) \Vert_{L^{2}_{t}L^{2}_{x}} \Vert  \vu \otimes  (\sqrt{w_\gamma} \vu ) \Vert_{L^{2}_{t}L^{2}_{x}}\\
\leq & C_{\gamma,T} ( \Vert \sqrt{w_\gamma} \tilde  \vu \Vert_{L^{\infty}L^{2}_{x}} + \Vert \sqrt{w_\gamma} \nabla \tilde  \vu \Vert_{L^{2}L^{2}_{x}}) \Vert  \vu \Vert_{\X_T}\Vert \sqrt{w_\gamma} \vu \Vert_{E_T}\, ds \\
\leq & C_{\gamma,T} (\Vert \sqrt{w_\gamma} \tilde \vu \Vert_{L^{\infty}_{t}L^{2}_{x}} +\Vert \sqrt{w_\gamma} \nabla \tilde \vu\Vert_{L^{2}_{t}L^{2}_{x}})\\
&\times \Vert \vu \Vert_{\X_T}  (\Vert \sqrt{w_\gamma}  \vu \Vert_{L^{\infty}_{t}L^{2}_{x}} +\Vert \sqrt{w_\gamma} \nabla \vu\Vert_{L^{2}_{t}L^{2}_{x}})<+\infty.
\end{split}
\end{equation}
We study now the  term $\tilde \vu \cdot \nabla p$. Recall first that by Proposition $2.1$ of \cite{PF_OJ} we have 
\begin{equation*}
\begin{split}
\nabla p =& \nabla \left(\sum_{1 \leq i,j\leq 3} \Ri_i \Ri_j(u_iu_j-b_i b_j -F_{i,j}) \right)\\
=&  \nabla \left(\sum_{1 \leq i,j\leq 3} \Ri_i \Ri_j(u_iu_j-b_i b_j) \right)- \nabla \left( \sum_{1 \leq i,j\leq 3} \Ri_i \Ri_j(F_{i,j}) \right)\\
 = &\nabla p_1 + \nabla p_2.
\end{split}
\end{equation*}
The term $\tilde \vu \cdot \nabla p_2$ is easily estimated by the hypothesis on the tensor $\F$ and the computations above. Thereafter, for the term $\nabla p_1$, following the same estimates performed in (\ref{Estim-aux}), and using the fact that $\Ri_i \Ri_j$ is a bounded operator in the space $L^{2}$, we find
\begin{equation*}
\begin{split}
\left\vert \int_{0}^{t} \int \varphi \tilde \vu \cdot \nabla p_1 \, dx \,ds\right\vert \leq & C_{\gamma,T}  (\Vert \sqrt{w_\gamma} \tilde \vu \Vert_{L^{\infty}_{t}L^{2}_{x}} +\Vert \sqrt{w_\gamma} \nabla \tilde \vu\Vert_{L^{2}_{t}L^{2}_{x}}) \\
&\times \left(\Vert \vu \Vert_{\X_T} (\Vert \sqrt{w_\gamma}  \vu \Vert_{L^{\infty}_{t}L^{2}_{x}} +\Vert \sqrt{w_\gamma} \nabla  \vu\Vert_{L^{2}_{t}L^{2}_{x}})\right. \\
&\left.+ \Vert \vb \Vert_{\X_T} (\Vert \sqrt{w_\gamma}  \vb \Vert_{L^{\infty}_{t}L^{2}_{x}} +\Vert \sqrt{w_\gamma} \nabla  \vb\Vert_{L^{2}_{t}L^{2}_{x}})\right)<+\infty.
\end{split}
\end{equation*}\finpv

Once all the terms in (\ref{eq-aux}) are well-defined as distributions, by  the locally energy balances  (\ref{energloc3}) and (\ref{energloc4}) we have 
\begin{equation*}
\begin{split}
& \partial_t (\frac { | \vv  |^2 + | \vw  |^2}2) +\mu + \nu = \Delta(\frac { |\vu  |^2 + |\vb|^2 }2)- |\vN \vu  |^2 - |\vN \vb|^2  \\
&- \vN\cdot\left( [\frac{ |\vu  |^2}2 + \frac{ |\vb  |^2}2 +p]\vu   \right)  + \vN \cdot ([(\vu \cdot \vb)+ q] \vb ) + \vu \cdot(\vN\cdot\mathbb{F}) +\vb \cdot(\vN\cdot\mathbb{G}) \\
&+\Delta(\frac { | \tilde \vu  |^2 + | \tilde \vb |^2 }2)- |\vN \tilde \vu  |^2 - |\vN \tilde \vb|^2 -  \vN\cdot\left( [\frac{ | \tilde \vu  |^2}2 + \frac{ | \tilde \vb  |^2}2 + \tilde p] \tilde \vu   \right) \\
&  + \vN \cdot ([( \tilde \vu \cdot \tilde \vb)+ \tilde q] \tilde \vb ) + \tilde \vu \cdot(\vN\cdot\mathbb{F}) + \tilde \vb \cdot(\vN\cdot\mathbb{G})\\
&- \vu \cdot \partial_t  \tilde \vu - \tilde \vu \cdot \partial_t  \vu - \vb \cdot \partial_t \tilde \vb - \tilde \vb \cdot \partial_t  \vb,
\end{split}
\end{equation*}
which can be rewritten as:
\begin{equation*}
\begin{split}
& \partial_t (\frac { | \vv  |^2 + | \vw  |^2}2) +\mu + \nu = \underbrace{ \Delta(\frac { |\vu  |^2 + |\vb|^2 }2) +\Delta(\frac { | \tilde \vu  |^2 + | \tilde \vb |^2 }2)}_{(1)} \\
& \underbrace{- |\vN \vu  |^2 - |\vN \vb|^2  - |\vN \tilde \vu  |^2 - |\vN \tilde \vb|^2 - \vu \cdot \partial_t  \tilde \vu - \tilde \vu \cdot \partial_t  \vu - \vb \cdot \partial_t \tilde \vb - \tilde \vb \cdot \partial_t  \vb}_{(2)} \\
&- \vN\cdot\left( [\frac{ |\vu  |^2}2 + \frac{ |\vb  |^2}2 +p]\vu   \right)  - \vN\cdot\left( [\frac{ | \tilde \vu  |^2}2 + \frac{ | \tilde \vb  |^2}2 + \tilde p] \tilde \vu   \right)  \\
& + \vN \cdot ([(\vu \cdot \vb)+ q] \vb ) + \vN \cdot ([( \tilde \vu \cdot \tilde \vb)+ \tilde q] \tilde \vb ) \\
&+ \vu \cdot(\vN\cdot\mathbb{F}) +\vb \cdot(\vN\cdot\mathbb{G})  + \tilde \vu \cdot(\vN\cdot\mathbb{F}) + \tilde \vb \cdot(\vN\cdot\mathbb{G}).
\end{split}
\end{equation*}

Now, we use  (\ref{IdBase}) to treat the terms $(1)$ and $(2)$, and thus we obtain 
\begin{equation*}
\begin{split}
& \partial_t (\frac { | \vv  |^2 + | \vw  |^2}2) +\mu + \nu = \Delta (\frac{\vert \vv \vert^2 +\vert \vw \vert^2}{2}) - \vert \nabla \vv \vert^2-\vert \nabla \vw \vert^2\\ 
&\underbrace{- (\partial_t \vu - \Delta \vu )\cdot \tilde{\vu} - (\partial_t \tilde{\vu} - \Delta \tilde{\vu} )\cdot \vu - (\partial_t \vb - \Delta \vb )\cdot \tilde{\vb} -  (\partial_t \tilde{\vb} - \Delta \tilde{\vb} )\cdot \vb}_{(3)} \\
&- \vN\cdot\left( [\frac{ |\vu  |^2}2 + \frac{ |\vb  |^2}2 +p]\vu   \right)  - \vN\cdot\left( [\frac{ | \tilde \vu  |^2}2 + \frac{ | \tilde \vb  |^2}2 + \tilde p] \tilde \vu   \right)  \\
& + \vN \cdot ([(\vu \cdot \vb)+ q] \vb ) + \vN \cdot ([( \tilde \vu \cdot \tilde \vb)+ \tilde q] \tilde \vb ) \\
&+ \vu \cdot(\vN\cdot\mathbb{F}) +\vb \cdot(\vN\cdot\mathbb{G})  + \tilde \vu \cdot(\vN\cdot\mathbb{F}) + \tilde \vb \cdot(\vN\cdot\mathbb{G}).
\end{split}
\end{equation*}
Thereafter, to study the term $(3)$ we use the fact that  $(\vu, \vb, p, q)$ and $(\tilde \vu, \tilde \vb, \tilde p, \tilde q)$ are  two solutions of the equations (MHD). 
Then we find 
\begin{equation*}
 \begin{split}
   & \partial_t (\frac { | \vv  |^2 + | \vw  |^2}2) +\mu+ \nu  = \Delta(\frac { | \vv  |^2 + | \vw  |^2}2) - | \vN \vv  |^2 - | \vN \vw  |^2  \\
    &\underbrace{- \vN\cdot\left( \frac{ |\vu |^2}2\vu+ \frac{ |\tilde \vu |^2} 2 \tilde \vu \right) + ((\vu \cdot \vN) \vu) \cdot \tilde \vu + (( \tilde \vu \cdot \vN) \tilde \vu ) \cdot \vu}_{(4)} \\
    &\underbrace{- \vN\cdot\left( \frac{ |\vb |^2}2\vu+ \frac{ | \tilde \vb |^2 } 2 \tilde \vu \right) + ((\tilde \vu \cdot \vN) \tilde \vb) \cdot \vb + (( \vu \cdot \vN)  \vb) \cdot \tilde \vb}_{(5)} \\
    & + \vN \cdot((\vu \cdot \vv) \vb) - ( ( \tilde \vb \cdot  \vN ) \tilde \vb) \cdot \vu - ((\vb \cdot \vN ) \vb) \cdot \tilde \vu\\
    & + \vN \cdot((\tilde \vu \cdot \tilde \vv) \tilde \vb) - ((\tilde \vb \cdot \vN ) \tilde \vu) \cdot\vb - ((\vb \cdot \vN ) \vu ) \cdot\tilde \vb \\
    &-\vN\cdot(a \vv). 
 \end{split} 
 \end{equation*}
We will verify that each term in  the right side  can be written as a sum of terms of the form
\begin{equation*}
\vN \cdot ((\vx \cdot \vy) \vz)
\end{equation*}
where at least two elements of $\{ \vx, \vy, \vz \}$ belong to $ \{ \vv, \vw  \} $, or  terms of the form
\begin{equation*}
( ( \vx \cdot \vN ) \vy ) \cdot \vz
\end{equation*}
where $\vy \in \{ \vv , \vw \}$ and at least one element of $\{ \vx, \vz \}$ belongs to $ \{ \vu, \vb  \} $. As we will see this fact permit to use the hypothesis $\vu, \tilde \vb \in \mathbb{X}_T $  to obtain a good control and use the Gr\"onwall inequality.\\

We start by studying the terms $(4)$ and $(5)$. For the term $(4)$, remark that we can write
\begin{equation*}
\begin{split}
    ((\vu& \cdot \vN) \vu) \cdot \tilde \vu + ((\tilde \vu \cdot \vN) \tilde \vu) \cdot \vu \\
    =& ((\vu \cdot \vN)  \vv) \cdot \tilde \vu +  ((\vu \cdot \vN)  \tilde \vu) \cdot \tilde \vu  - ((\tilde \vu \cdot \vN)  \vv ) \cdot \vu + ((\tilde \vu \cdot \vN)  \vu ) \cdot \vu \\
    =& \vN \cdot \left(\frac{| \tilde \vu |^2}{2} \vu + \frac{| \vu |^2}{2} \tilde \vu \right) - ((\tilde \vu \cdot \vN)  \vv ) \cdot \vu - ((\vu \cdot \vN)  \vv)  \cdot  \vv + ((\vu \cdot \vN)  \vv) \cdot  \vu  \\
    =& \vN \cdot \left(\frac{| \tilde \vu |^2}{2} \vu + \frac{| \vu |^2}{2} \tilde \vu  - \frac{|\vv|^2}{2} \vu \right) + (( \vv \cdot \vN)  \vv ) \cdot \vu ,
\end{split}
\end{equation*}
hence we obtain 
 \begin{equation*}
    (4) = - \vN\cdot\left( \frac{ \vv \cdot (\vu + \tilde \vu ) }2 \vv + \frac{|\vv|^2}{2} \vu \right) + ((\vv \cdot \vN) \vv) \cdot \vu 
 \end{equation*}
In a similar way, for the term $(5)$, observe that we have  
\begin{equation*} 
\begin{split}
    ((  \tilde  \vu & \cdot \vN) \tilde \vb) \cdot \vb + (( \vu \cdot \vN)  \vb) \cdot \tilde \vb \\
    =& -((\tilde \vu \cdot \vN)  \vw) \cdot \vb + ((\tilde \vu \cdot \vN)  \vb) \cdot \vb  + (( \vu \cdot \vN)  \vw ) \cdot \tilde \vb + (( \vu \cdot \vN)  \tilde \vb ) \cdot \tilde \vb \\
    =& \vN \cdot \left(\frac{| \vb |^2}{2} \tilde \vu + \frac{| \tilde \vb |^2}{2}  \vu \right) -((\tilde \vu \cdot \vN)  \vw) \cdot \vb + (( \vu \cdot \vN)  \vw ) \cdot \tilde \vb
\end{split}
\end{equation*}
and thus we obtain
 \begin{equation*}
(5) = - \vN \cdot \left( \frac{ \vw \cdot (\vb +  \tilde \vb ) }2 \vv  + \frac{|\vw|^2}{2} \vu\right) + (( \vv \cdot \vN)  \vw ) \cdot  \vb 
\end{equation*}
With these identities on the terms $(4)$ and $(5)$ we have 
\begin{equation*}
 \begin{split}
   & \partial_t (\frac { | \vv  |^2 + | \vw  |^2}2) +\mu+ \nu= \Delta(\frac { | \vv  |^2 + | \vw  |^2}2) - | \vN \vv  |^2 - | \vN \vw  |^2  \\
    &- \vN\cdot\left( \frac{ \vv \cdot (\vu + \tilde \vu ) }2 \vv + \frac{|\vv|^2}{2} \vu \right) + ((\vv \cdot \vN) \vv) \cdot \vu  \\
    &- \vN \cdot \left( \frac{ \vw \cdot (\vb +  \tilde \vb ) }2 \vv + \frac{|\vw|^2}{2} \vu \right)  + (( \vv \cdot \vN)  \vw ) \cdot  \vb \\
    & \underbrace{+ \vN \cdot (( \vu \cdot \vb) \vb)   - ((\tilde \vb \cdot  \vN ) \tilde \vb) \cdot \vu - ((\vb \cdot \vN ) \vb ) \cdot \tilde \vu}_{(6)}\\
    &\underbrace{ + \vN \cdot (( \tilde \vu \cdot \tilde \vb )  \tilde \vb)  - ((\tilde \vb \cdot \vN ) \tilde \vu ) \cdot\vb - ((\vb \cdot \vN ) \vu) \cdot\tilde \vb }_{(7)},  \\
    &-\vN\cdot(a \vv) ,
 \end{split}
\end{equation*}
where, it remains to treat the terms $(6)$ and $(7)$.  As
\begin{equation*}
\begin{split}
    \vN \cdot & (( \vu \cdot \vb) \vb) + \vN \cdot (( \tilde \vu \cdot \tilde \vb )  \tilde \vb) \\
    =& \vN \cdot (( \vu \cdot \vw) \vw) + \vN \cdot (( \vu \cdot \vw) \tilde \vb) + \vN \cdot (( \vu \cdot \tilde \vb) \vb) \\
    & + \vN \cdot ((  \vv \cdot \tilde \vb )   \vw) - \vN \cdot ((  \vv \cdot  \tilde \vb )  \vb) + \vN \cdot ((  \vu \cdot \tilde \vb )  \tilde \vb) \\
    =&  \vN \cdot ((\vu \cdot \vw)  \vw )+ \vN \cdot ((\vv \cdot \tilde   \vb) \vw) + \vN \cdot  ((\vu \cdot \vb ) \tilde \vb)  + \vN \cdot ((\tilde \vu \cdot \tilde \vb) \vb)
\end{split}
\end{equation*}
we  have  
\begin{equation*}
\begin{split}
(6)+ (7)=& \vN \cdot ((\vu \cdot \vw)  \vw )+ \vN \cdot ((\vv \cdot \tilde   \vb) \vw) + \vN \cdot  ((\vu \cdot \vb ) \tilde \vb)  + \vN \cdot ((\tilde \vu \cdot \tilde \vb) \vb) \\
& - ((\tilde \vb \cdot \vN ) \tilde \vb) \cdot \vu - ( (\vb \cdot \vN ) \vb ) \cdot \tilde \vu - ((\tilde \vb \cdot \vN ) \tilde \vu ) \cdot\vb - ( (\vb \cdot \vN ) \vu ) \cdot\tilde \vb \\
=&  \vN \cdot ((\vu \cdot \vw)  \vw ) + \vN \cdot ((\vv \cdot \tilde \vb) \vw) \\
&+ ((\tilde \vb  \cdot \vN \vu) \cdot \vb) + ((\tilde \vb \cdot  \vN \vb) \cdot \vu) + (( \vb \cdot \vN \tilde \vu) \cdot \tilde \vb) + (( \vb \cdot \vN \tilde \vb) \cdot \tilde \vu)
  \\
& - ( (\tilde \vb \cdot  \vN) \tilde \vb) \cdot \vu - ( (\vb \cdot \vN ) \vb) \cdot \tilde \vu - ((\tilde \vb \cdot \vN ) \tilde \vu) \cdot\vb - ( (\vb \cdot \vN ) \vu) \cdot \tilde \vb \\
= & \vN \cdot ((\vu \cdot  \vw) \vw) + \vN \cdot ((\vv \cdot \tilde \vb) \vw) +    \vN \cdot((\vv \cdot \vw) \vb )  \\
 & - ( ( \vw \cdot \vN ) \vw)\cdot \vu   - ( ( \vw \cdot \vN ) \vv )\cdot \vb \\
\end{split}
\end{equation*}

 Thus, by this  identity   we are able to write 
\begin{equation*}
 \begin{split}
    &\partial_t (\frac { | \vv  |^2 + | \vw  |^2}2) +\mu+ \nu = \Delta(\frac { | \vv  |^2 + | \vw  |^2}2) - | \vN \vv  |^2 - | \vN \vw  |^2  \\
    &- \vN\cdot\left( \frac{ \vv \cdot (\vu + \tilde \vu ) }2 \vv + \frac{|\vv|^2}{2} \vu \right) + ((\vv \cdot \vN) \vv) \cdot \vu  \\
    &- \vN \cdot \left( \frac{ \vw \cdot (\vb +  \tilde \vb ) }2 \vv + \frac{|\vw|^2}{2} \vu \right)  + (( \vv \cdot \vN)  \vw ) \cdot  \vb \\
    & + \vN \cdot ((\vv \cdot \tilde \vb) \vw) + \vN \cdot ((\vu \cdot  \vw) \vw) +   \vN \cdot((\vv \cdot \vw) \vb )  \\
    & - ( ( \vw \cdot \vN ) \vw )\cdot \vu   - ( ( \vw \cdot \vN ) \vv )\cdot \vb \\
    &-\vN\cdot(a \vv), 
 \end{split}
\end{equation*}

which can be rewritten as 
\begin{equation}\label{energloc}
 \begin{split}
    & \partial_t (\frac { | \vv  |^2 + | \vw  |^2}2)  + | \vN \vv  |^2 + | \vN \vw  |^2 + \mu+ \nu \\
    = &  \Delta(\frac { | \vv  |^2 + | \vw  |^2}2)  - \vN\cdot \underbrace{ \left( \frac{ \vv \cdot (\vu + \tilde \vu ) }2 \vv  +   \frac{ \vw \cdot (\vb +  \tilde \vb ) }2 \vv  + \frac{|\vv|^2 +|\vw|^2  }{2} \vu   \right)}_{\vA_1} \\
    & + \vN \cdot \underbrace{ \left(  (\vv \cdot \tilde \vb) \vw + (\vu \cdot  \vw) \vw  +  (\vv \cdot \vw) \vb \right)}_{\vA_2} \\
    & \underbrace{ - ( ( \vw \cdot \vN ) \vw)) \cdot \vu   - ( ( \vw \cdot \vN ) \vv)) \cdot  \vb + ((\vv \cdot \vN) \vv)  \cdot  \vu + (( \vv \cdot \vN)  \vw ) \cdot \vb }_{A_3} \\
    & -\vN\cdot(a \vv) . 
 \end{split}
 \end{equation}
We will apply \eqref{energloc} to a suitable test function which we shall define as follows. First, we consider a function  $\alpha_{\eta,t_0,t_1}$ which converges almost everywhere to $\mathds 1 _{[t_0, t_1]}$ when $\eta \to 0$  and such that $\partial_t \alpha_{\eta,t_0,t_1}$ is the difference between two identity approximations, the first one in $t_0$ and the second one in $t_1$. For this, let $
\alpha\in \mathcal{C}^\infty(\mathbb{R})$ be a non-decreasing function which is equal to $0$ on $(-\infty,\frac{1}{2})$ and equal to $1$ on $(1, +\infty)$.  For $0<\eta< \min(\frac {t_0}2,T-t_1) $ we set 
\begin{equation}
\label{alpha}
\alpha_{\eta,t_0,t_1}(t)=\alpha( \frac{t-t_0}\eta)-\alpha(\frac{t-t_1}\eta).
\end{equation}
Thereafter,  we consider a non-negative function $\phi\in\mathcal{D}(\mathbb{R}^3)$ which is equal to $1$ for $ | x |\leq 1$ and to $0$ for $ | x |\geq 2$ and we set 
\begin{equation}
\label{phi}
\phi_R(x)=\phi(\frac x R).
\end{equation}
Finally, for $\epsilon>0$ we define the function  $w_{\gamma,\epsilon}= \frac 1{(1+\sqrt{\epsilon^2+ | x |^2})^\delta}$. We may observe that  $\alpha_{\eta,a,s}(t)\phi_R(x) w_{\gamma,\epsilon}(x)$ belongs to $\mathcal{D}((0,T)\times\mathbb{R}^3)$ and $\alpha_{\eta,a,s}(t)\phi_R(x) w_{\gamma,\epsilon}(x) \geq 0$. Thus, using the local energy balance (\ref{energloc}) with this particular test function  we obtain 
\begin{equation*}
\begin{split}
-\iint &\frac{ \vert \vv \vert^2+  |\vw  |^2}2 \partial_t\alpha_{\eta,t_0,t_1} \phi_R w_{\gamma,\epsilon}\, dx\, ds + \iint (\vert \vN \vv \vert^2+  | \vN \vw  |^2)\, \,  \alpha_{\eta,t_0,t_1} \phi_R w_{\gamma,\epsilon} dx\, ds \\ \leq&-\sum_{i} \iint \partial_i( \vv \cdot \vv + \vw\cdot \vw )\,  \alpha_{\eta,t_0,t_1}  (w_{\gamma,\epsilon}\partial_i \phi_R+\phi_R \partial_iw_{\gamma,\epsilon})\, dx\, ds \\
&- \sum_{i} \iint ( \vA_1+\vA_2 )_i \, \alpha_{\eta,t_0,t_1}  (w_{\gamma,\epsilon}\partial_i \phi_R+\phi_R \partial_i w_{\gamma,\epsilon})\, dx\, ds \\
&+ \iint  A_3 \, \alpha_{\eta,t_0,t_1}  \phi_R w_{\gamma,\epsilon} \, dx\, ds \\
&  + \sum_{i} \iint (a \vv)_i \, \alpha_{\eta,t_0,t_1}  (w_{\gamma,\epsilon}\partial_i \phi_R+\phi_R \partial_i w_{\gamma,\epsilon})\, dx\, ds. 
\end{split}
\end{equation*}
In this inequality, we  take the limit when $\eta \to 0$.  By the dominated convergence theorem we obtain (when the limit in the left side is well-defined)
\begin{equation*}
\begin{split}
- \lim_{\eta \to 0}& \iint \frac{ \vert \vv \vert^2+  |\vw  |^2}2 \partial_t\alpha_{\eta,t_0,t_1} \phi_R w_{\gamma,\epsilon}\, dx\, ds + \int_{t_0}^{t_1} \int (\vert \vN \vv \vert^2+  | \vN \vw  |^2)\, \,   \phi_R w_{\gamma,\epsilon} dx\, ds \\ \leq&-\sum_{i} \int_{t_0}^{t_1} \int \partial_i( \vv \cdot \vv + \vw\cdot \vw )\,   (w_{\gamma,\epsilon}\partial_i \phi_R+\phi_R \partial_iw_{\gamma,\epsilon})\, dx\, ds \\
&- \sum_{i} \int_{t_0}^{t_1} \int ( \vA_1+\vA_2 )_i \,   (w_{\gamma,\epsilon}\partial_i \phi_R+\phi_R \partial_i w_{\gamma,\epsilon})\, dx\, ds +  \int_{t_0}^{t_1} \int A_3 \,  \phi_R w_{\gamma,\epsilon} \, dx\, ds \\
&  + \sum_{i} \int_{t_0}^{t_1} \int (a \vv)_i  \,   (w_{\gamma,\epsilon}\partial_i \phi_R+\phi_R \partial_i w_{\gamma,\epsilon})\, dx\, ds. 
\end{split}
\end{equation*}
Now, if $t_0$ and $t_1$ are Lebesgue points of the measurable function $$ A_{R,\epsilon}(t)=\int  (\vert \vv (t,x) \vert^2+  | \vw (t,x) |^2   )\phi_R(x)  w_{\gamma,\epsilon}(x)\, dx,$$
and moreover,  as we have 
$$-\int \int \frac{\vert \vv \vert^2+  | \vw  |^2}2 \partial_t\alpha_{\eta,a,s}   \phi_R w_{\gamma,\epsilon} \, dx\, ds=-\frac 1 2\int \partial_t\alpha_{\eta,a,s}
A_{R,\epsilon}(s) \, ds,$$
then  we obtain 
$$ \lim_{\eta\rightarrow 0}  -\int \int \frac{\vert \vv \vert^2 +  | \vw  |^2}2 \partial_t\alpha_{\eta,a,s}  \phi_R w_{\gamma,\epsilon} \, dx\, ds=\frac 1 2 (  A_{R,\epsilon}(t_1)- A_{R,\epsilon}(t_0)).$$
Thereafter, the continuity at $0$ of $\vv$ and $\vw$ permit to let $t_0$ go to $0$ and thus we  replace $t_0$ by $0$ in this inequality. Moreover, if we let $t_1$ go to $t$, where $t \in (0,T)$, then by the  weak continuity  we obtain
$\ds{ A_{R,\epsilon}(t)\leq \liminf_{t_1\rightarrow t }  A_{R,\epsilon}(t_1)}$, so  we may as well replace $t_1$ by $t$. Thus,  for every $t\in (0,T)$  we have 
\begin{equation*}
\begin{split}
\int & \frac{ \vert \vv(t,\cdot)\vert^2+ | \vw(t,\cdot)  |^2}2   \phi_R w_{\gamma,\epsilon} \, dx + \int_{t_0}^{t_1} \int (\vert \vN \vv \vert^2+  | \vN \vw  |^2)\, \,   \phi_R w_{\gamma,\epsilon} dx\, ds \\ \leq&-\sum_{i} \int_{0}^{t} \int \partial_i( \vv \cdot \vv + \vw\cdot \vw )\,   (w_{\gamma,\epsilon}\partial_i \phi_R+\phi_R \partial_iw_{\gamma,\epsilon})\, dx\, ds \\
&- \sum_{i} \int_{0}^{t} \int ( \vA_1+\vA_2 )_i \,   (w_{\gamma,\epsilon}\partial_i \phi_R+\phi_R \partial_i w_{\gamma,\epsilon})\, dx\, ds \\
&+  \int_{0}^{t} \int A_3 \,   w_{\gamma} \, dx\, ds  + \sum_{i} \int_{0}^{t} \int  (a \vv)_i \,   (w_{\gamma,\epsilon}\partial_i \phi_R+\phi_R \partial_i w_{\gamma,\epsilon})\, dx\, ds. 
\end{split}
\end{equation*}
In this inequality, we take now  limit when $R \to +\infty $, and moreover, the limit when  $\epsilon \to 0$ to obtain: 
\begin{equation*}
\begin{split}
\int & \frac{ \vert \vv(t,\cdot)\vert^2+ | \vw(t,\cdot)  |^2}2   w_\gamma \, dx + \int_{0}^{t} \int (\vert \vN \vv \vert^2+  | \vN \vw  |^2)\,  w_{\gamma} dx\, ds \\ 
\leq&-\int_{0}^{t} \int \vN ( \vv \cdot \vv + \vw\cdot \vw ) \cdot \vN w_\gamma \, dx\, ds - \underbrace{ \int_{0}^{t} \int ( \vA_1+\vA_2 ) \cdot \vN w_\gamma \, dx\, ds}_{I_1} \\
&+ \underbrace{ \int_{0}^{t} \int A_3  \,  w_{\gamma} \, dx\, ds}_{I_2}  + \underbrace{ \int_{0}^{t} \int  a \vv \cdot \vN w_\gamma\, dx\, ds}_{I_3}. 
\end{split}
\end{equation*}
At this point, in the following technical lemmas  we  estimate the terms $I_1$, $I_2$ and $I_3$. To make the notation more convenient we write 
$$ \int   ( | \vv |^2+ \vert \vw\vert^2 )w_\gamma dx=\Vert \sqrt{w_\gamma}(\vv, \vw)\Vert^{2}_{L^2}, \,\int (| \vN \vv |^2 +| \vN \vw |^2)w_\gamma dx= \Vert \sqrt{w_\gamma}\vN(\vv, \vw)\Vert^{2}_{L^2},$$
and then we have 
\begin{Lemma} $\ds{\vert I_1 \vert \leq   C_\gamma \int_0^t \Vert \sqrt{w_\gamma}(\vv, \vw)(s,\cdot)\Vert^{2}_{L^2}\, ds  + \frac{1}{4} \int_0^t \Vert \sqrt{w_\gamma}\vN(\vv, \vw)(s,\cdot)\Vert^{2}_{L^2} \, ds }$. 
\end{Lemma}	 
\pv  We observe first that we have $ |\vN w_{\gamma}| \leq C_\gamma w_{\frac{3}{2}\gamma} $ and then we get 
$$ \vert I_1 \vert \leq C_\gamma \int_{0}^{t}\int \vert \vA_1+ \vA_2 \vert w_{3/2 \gamma} dx\, ds.$$
We observe moreover that each term in the expression $\vA_1+\vA_2$ writes down as the product of three vectors: $(\vx \cdot \vy)\vz$ where at least two vectors belong to $\{ \vv, \vw\}$ and the third one belongs to $\{\vu, \vb , \tilde \vu, \tilde \vb \}$.  So, we will estimate the generic expression $(\vx \cdot \vy)\vz$, where, without loss of generality we may assume that $\vx, \vz \in \{\vv, \vw \}  $ and $\vy \in \{\vu, \vb , \tilde \vu, \tilde \vb \}$.   Remark that for $\delta > 0$ (which we will set later), by  the H\"older inequalities and the Young inequalities we have 
\begin{equation*}
\begin{split}
&\int \vert (\vx \cdot \vy)\vz \vert w_{3/2 \gamma} dx \leq \int (\sqrt{w_\gamma} \vert \vx \vert ) (\sqrt{w_\gamma} \vert \vz \vert ) (\sqrt{w_\gamma} \vert \vy \vert ) dx  \\
\leq & \Vert \sqrt{w_\gamma} \vx \Vert_{L^3}  \Vert \sqrt{w_\gamma} \vz \Vert_{L^6} \Vert \sqrt{w_\gamma} \vy \Vert_{L^2} \leq \underbrace{ \delta^{-1}  \Vert \sqrt{w_\gamma} \vx \Vert^{2}_{L^3}+ \delta  \Vert \sqrt{w_\gamma} \vz \Vert^{2}_{L^6} \Vert \sqrt{w_\gamma} \vy \Vert^{2}_{L^2}}_{(a)}.
\end{split}
\end{equation*}
Moreover, by the interpolation inequalities and the Young inequalities  we have 
$$ \Vert \sqrt{w_\gamma} \vx \Vert^{2}_{L^3} \leq \delta^{-2} \Vert\sqrt{w_\gamma}  \vx  \Vert^{2}_{L^2}+ \delta^2 \Vert \sqrt{w_\gamma} \vx \Vert^{2}_{L^6}, $$ hence we can write 
$$ (a) \leq \delta^{-3} \Vert \sqrt{w_\gamma} \vx \Vert^{2}_{L^2} + \delta \Vert \sqrt{w_\gamma} \vx \Vert^{2}_{L^6}+  \delta  \Vert \sqrt{w_\gamma} \vz \Vert^{2}_{L^6} \Vert \sqrt{w_\gamma} \vy \Vert^{2}_{L^2}=(b).$$
At this point, we use the Sobolev embedding (\ref{sobol}) to estimate the terms $\Vert \sqrt{w_\gamma} \vx \Vert^{2}_{L^6}$ and $ \Vert \sqrt{w_\gamma} \vz \Vert^{2}_{L^6} $ and we obtain 
\begin{equation*}
\begin{split}
&(b)\leq  \delta^{-3} \Vert \sqrt{w_\gamma} \vx \Vert^{2}_{L^2}  + \delta (\Vert \sqrt{w_\gamma} \vx \Vert^{2}_{L^2}+  \Vert \sqrt{w_\gamma} \vN \vx \Vert^{2}_{L^2})\\
&+\delta (\Vert \sqrt{w_\gamma} \vz \Vert^{2}_{L^2}+  \Vert \sqrt{w_\gamma} \vN \vz \Vert^{2}_{L^2}) \Vert \sqrt{w_\gamma} \vy \Vert^{2}_{L^2}.
\end{split}
\end{equation*}
Now,  by the previous estimate we get 
\begin{equation*}
\begin{split}
(b) \leq & \delta^{-3} \Vert \sqrt{w_\gamma}(\vv, \vw)\Vert^{2}_{L^2} + \delta (\Vert \sqrt{w_\gamma}(\vv, \vw)\Vert^{2}_{L^2}+ \Vert \sqrt{w_\gamma}\vN(\vv, \vw)\Vert^{2}_{L^2})\\
&+ \delta (\Vert \sqrt{w_\gamma}(\vv, \vw)\Vert^{2}_{L^2}+ \Vert \sqrt{w_\gamma}\vN(\vv, \vw)\Vert^{2}_{L^2}) \left( \sup_{0<s<T} \Vert \sqrt{w_\gamma} \vy \Vert^{2}_{L^2} \right).
\end{split}
\end{equation*}
We observe that we can set  the parameter $\delta$ small enough such that it verifies $\ds{\max\left[\delta, \delta \left( \sup_{0<s<T} \Vert \sqrt{w_\gamma} \vy \Vert^{2}_{L^2} \right)\right]\leq 1/64}$, and by the previous estimate we finally get 
$$  \int \vert (\vx \cdot \vy)\vz \vert w_{3/2 \gamma} dx \leq C_\gamma \Vert \sqrt{w_\gamma}(\vv, \vw)\Vert^{2}_{L^2} + \frac{1}{64}\Vert \sqrt{w_\gamma}\vN(\vv, \vw)\Vert^{2}_{L^2},$$ hence, integrating in the temporal variable we have the desired estimate.  \finpv

To estimate the term $I_2$ we use the information $\vu, \vb \in {\mathbb{X}} _T$ and we have the following result.
\begin{Lemma} Assume that $\vu, \vb \in {\mathbb{X}} _T$. Then, for all $0<t<T$ we have 
	\begin{equation*}
	\begin{split}
	\vert I_2 \vert  \leq & C\| \vu   \|_{{\mathbb{X}} _T}  \left( \sup_{ 0 \leq s \leq t}\|\sqrt{w_\gamma} (\vv, \vw) (s) \,  \|_{L^2}^2 + \int_0^{t}  \| \sqrt{w_\gamma}    \vN (\vv, \vw)(s,\cdot)  \|_{L^2}^2 \, ds \right)\\
	&+ C\| \vb   \|_{{\mathbb{X}} _T}  \left( \sup_{ 0 \leq s \leq t}\| \sqrt{w_\gamma}  (\vv, \vw) (s,\cdot) \,  \|_{L^2}^2  + \int_0^{t}  \|  \sqrt{w_\gamma}    \vN (\vv, \vw)(s,\cdot)  \, \|_{L^2}^2 \, ds \right).  
	\end{split}
	\end{equation*}
\end{Lemma}	
\pv  We observe that each term in the expression $A_3$ writes down as $\ds{\sum_{i,j=1}^{3} x_j (\partial_j z_i) y_i}$, where, without loss of generality we way assume that $\{\vx, \vz\}$ belong to $\{ \vv, \vw \}$ and $\vy \in \{ \vu, \vb \}$.   Then, to estimate each term, by the Cauchy-Schwarz inequalities (in the spatial and temporal variables) we write 
\begin{equation*}
\begin{split}
\int_{0}^{t}\Vert w_\gamma \vert \vx \vert \vert \vN \vz \vert \vert \vy \vert  (s,\cdot)\Vert_{L^1}\, ds  \leq & \left( \int_{0}^{t} \Vert \sqrt{w_\gamma} \vert \vx \vert \vert \vy \vert(s,\cdot) \Vert^{2}_{L^2}\,ds \right)^{1/2}\\
&\times  \left( \int_{0}^{t} \Vert \sqrt{w_\gamma} \vert \vN\vz \vert (s,\cdot)\Vert^{2}_{L^2}\,ds \right)^{1/2}.  
\end{split}
\end{equation*}
Thereafter, in order to estimate the first term in the right side, by definition of the multiplier space $\mathbb{X}_T$, and moreover, by definition of the energy space $E_T$   we can write 
\begin{equation*}
\begin{split}
&\int_{0}^{t}\Vert w_\gamma \vert \vx \vert \vert \vN \vz \vert \vert \vy \vert  (s,\cdot)\Vert_{L^1}\, ds  \leq \Vert \vy \Vert_{\mathbb{X}_T} \Vert \sqrt{w_\gamma} \vx \Vert_{E_T} \, \left( \int_{0}^{t} \Vert \sqrt{w_\gamma} \vert \vN\vz \vert (s,\cdot)\Vert^{2}_{L^2}\,ds \right)^{1/2}\\
\leq & C\,\Vert \vy \Vert_{\mathbb{X}_T} \left( \Vert \sqrt{w_\gamma} \vx \Vert^{2}_{E_T}+ \int_{0}^{t} \Vert \sqrt{w_\gamma} \vert \vN\vz \vert (s,\cdot)\Vert^{2}_{L^2}\,ds\right),
\end{split}
\end{equation*}
hence we get the desired estimate. \finpv  

Finally, for the term $I_3$ we have: 
\begin{Lemma}  $\ds{	\vert I_3 \vert  \leq  C \int_0^t   \|  \sqrt{w_\gamma} ( \vv, \vw)(s,\cdot) \|_{L^2}^2\, ds + \frac{1}{4} \int_0^t \Vert \sqrt{w_\gamma}  \vN (\vv,\vw)(s,\cdot) \Vert^{2}_{L^2}\, ds}$.
\end{Lemma}	
\pv Remark first that by the characterization of the pressure terms $p$ and $\tilde p$ given in (\ref{Charact-p}) we have $\ds{p - \tilde p  = \sum_{i,j} \Ri_i \Ri_j (u_iu_j - \tilde u_i \tilde u_j - b_i b_j+ \tilde b_i \tilde b_j)}$, where, as $\vv=\vu-\tilde \vu$ and $\vw=\vb-\tilde \vb$ then we can write 
\begin{equation*}
\begin{split}
 u_iu_j - \tilde u_i \tilde u_j - b_i b_j+ \tilde b_i \tilde b_j = v_i u_j +\tilde u_i v_j - w_j b_j -\tilde b_i w_j.   
\end{split}    
\end{equation*}
Then, as we have $\vert \vN w_\gamma \vert \leq C_\gamma w_{\frac{3}{2}\gamma}$,  by the H\"older inequalities  we can  write 
\begin{equation*}
\begin{split}
\vert I_3 \vert     \leq&  C_\gamma \left\Vert   w_\gamma \left( \sum_{i,j} \Ri_i \Ri_j (u_iu_j - \tilde u_i \tilde u_j - b_i b_j+ \tilde b_i \tilde b_j) \right) \right\Vert_{L^{\frac{6}{5}}} \Vert  \sqrt{w_\gamma}  \vv \Vert_{L^6} \\
\leq&  C_\gamma \left\Vert   w_\gamma \left( \sum_{i,j} \Ri_i \Ri_j (v_i u_j +\tilde u_i v_j - w_j b_j -\tilde b_i w_j) \right) \right\Vert_{L^{\frac{6}{5}}} \Vert  \sqrt{w_\gamma}  \vv \Vert_{L^6}\\
\leq& C_\gamma( \| w_\gamma (|\vu| + |\tilde \vu|) |\vv|  \,  \|_{L^{\frac{6}{5}}} + \| w_\gamma (|\vb| + |\tilde \vb|) |\vw|  \,  \|_{L^{\frac{6}{5}}} )  \| \sqrt{ w_\gamma}  \vv \,  \|_{L^6} \\
\leq& C_\gamma \| \sqrt{w_\gamma} (|\vu| + |\tilde \vu|)\|_{L^2} \| \sqrt{w_\gamma} \vv  \,  \|_{L^3} \|  \sqrt{w_\gamma}   \vv \,  \|_{L^6} \\
& +  C_\gamma \| \sqrt{w_\gamma} (|\vb| + |\tilde \vb|)\|_{L^2} \| \sqrt{w_\gamma} \vw  \,  \|_{L^3} \| \sqrt{w_\gamma}  \vv \,  \|_{L^6} \\
\leq& C_\gamma \| \sqrt{w_\gamma} \vv  \,  \|_{L^3} \|  \sqrt{w_\gamma}   \vv \,  \|_{L^6} \| \sqrt{w_\gamma} (|\vu| + |\tilde \vu|)\|_{L^2}  \\
&+ C_\gamma \| \sqrt{w_\gamma} \vv  \,  \|_{L^3} \|  \sqrt{w_\gamma}   \vw \,  \|_{L^6} \| \sqrt{w_\gamma} (|\vb| + |\tilde \vb|)\|_{L^2}\\
&=(a). 
\end{split}
\end{equation*}
Now, for $\delta>0$ we use first  the Young inequalities, and thereafter, we use  the Sobolev embedding (\ref{sobol}) to estimate the terms $\|  \sqrt{w_\gamma}   \vv\,  \|_{L^6} $ and $\|  \sqrt{w_\gamma}   \vw \,  \|_{L^6} $. Thus  we get 
\begin{equation*}
\begin{split} 
(a)   \leq&  C_\gamma \delta^{-1} \|  \sqrt{w_\gamma}  \vv \,  \|_{L^3}^2 + C_\gamma \delta  \| \sqrt{w_\gamma} \vv  \|_{L^6}^2 \| \sqrt{w_\gamma} (| \vu| + |\tilde \vu|)  \|_{L^2}^2  \\  
& +  C_\gamma \delta^{-1} \|  \sqrt{w_\gamma}  \vw \,  \|_{L^3}^2 + C \delta  \| \sqrt{w_\gamma}  \vv  \|_{L^6}^2 \| \sqrt{w_\gamma} (| \vb| + |\tilde \vb|)  \|_{L^2}^2  \\ 
\leq&  C\delta^{-3}  \|  \sqrt{w_\gamma}  \vv \,  \|_{L^2}^2 + C_\gamma\delta  (\| \sqrt{w_\gamma}  \vv  \|_{L^2}^2+\| \sqrt{w_\gamma} | \vN \vv |   \|_{L^2}^2 ) \\
&+ C_\gamma\delta  (\| \sqrt{w_\gamma}  \vv  \|_{L^2}^2+\| \sqrt{w_\gamma} | \vN \vv |   \|_{L^2}^2 )\left( \sup_{0 \leq s \leq T}  \| \sqrt{w_\gamma} (| \vu| + |\tilde \vu|)  \|_{L^2}\right)  \\ 
& +  C_\gamma\delta^{-3}  \|  \sqrt{w_\gamma} \vw \,  \|_{L^2}^2 + C\delta  (\| \sqrt{w_\gamma}  \vw  \|_{L^2}^2+\| \sqrt{w_\gamma} | \vN \vw |   \|_{L^2}^2 ) \\
&+ C_\gamma\delta  (\| \sqrt{w_\gamma}  \vw  \|_{L^2}^2+\| \sqrt{w_\gamma} | \vN \vw |   \|_{L^2}^2 )\left( \sup_{0 \leq s \leq T}  \| \sqrt{w_\gamma} (| \vb| + |\tilde \vb| \|_{L^2}\right).
\end{split}
\end{equation*}
Hence, setting the parameter $\delta$ small enough such it verifies: 
$$ \max\left[ \delta, \delta \left( \sup_{0 \leq s \leq T}  \| \sqrt{w_\gamma} (| \vu| + |\tilde \vu| ) \|_{L^2}\right), \delta \left( \sup_{0 \leq s \leq T}  \| \sqrt{w_\gamma} (| \vb| + |\tilde \vb| ) \|_{L^2}\right) \right] \leq 1/4,$$
we have the desired estimate. \finpv 

Thus, for $0 \leq t_0 < t_1 <T $ if we suppose that $\vv = \vw = 0$ on $[0,t_0]$, we have for all $t \in [t_0 , t_1] $
\begin{equation*}
\begin{split}
\| & (\vv,\vw)(t,\cdot) \,  \|_{L^2_{w_\gamma}}^2  + \frac{1}{2} \int_0^t \| \vN  (\vv, \vw)   \, \|_{L^2_{w_\gamma}}^2 \,  ds \\
\leq&  C \int_0^t   \|  (\vv, \vw)(s) \,  \|_{L^2_{w_\gamma}}^2 ds \\
&+ C\| \mathds{1}_{(t_0 , t_1)} \vu   \|_{X_T}  \left( \sup_{ 0 \leq s \leq t_1}\| (\vv, \vw) (s) \,  \|_{L^2_{w_\gamma}}^2 ds + \int_0^{t}  \|    \vN (\vv, \vw)(s)   \, \|_{L^2_{w_\gamma}}^2 ds \right)\\
&+ C\|\mathds{1}_{(t_0 , t_1)} \vb   \|_{X_T}  \left( \sup_{ 0 \leq s \leq t_1}\| (\vv, \vw) (s) \,  \|_{L^2_{w_\gamma}}^2  + \int_0^{t}  \|    \vN (\vv, \vw)(s)   \, \|_{L^2_{w_\gamma}}^2 ds \right).
\end{split}
\end{equation*}
Thus, we can take the supreme on $[0,t_1]$  to obtain
\begin{equation*}
\begin{split}
\sup_{0 \leq t \leq t_1} \| & (\vv,\vw)(t,\cdot) \,  \|_{L^2_{w_\gamma}}^2  + \frac{1}{2} \int_{0}^{t_1} \| \vN  (\vv, \vw)   \, \|_{L^2_{w_\gamma}}^2 \,  ds \\
\leq&  C \int_{0}^{t_1}   \|  (\vv, \vw)(s) \,  \|_{L^2_{w_\gamma}}^2 ds \\
&+ C\| \mathds 1_{(t_0 , t_1)} \vu   \|_{X_T}  \left( \sup_{ 0 \leq s \leq t_1}\| (\vv, \vw) (s) \,  \|_{L^2_{w_\gamma}}^2  + \int_{0}^{t_1}  \|    \vN (\vv, \vw)(s)   \, \|_{L^2_{w_\gamma}}^2 ds \right)\\
&+ C\| \mathds 1 _{(t_0 , t_1)}  \vb  \|_{X_T}  \left( \sup_{ 0 \leq s \leq t_1}\| (\vv, \vw) (s) \,  \|_{L^2_{w_\gamma}}^2 + \int_{0}^{t_1}  \|    \vN (\vv, \vw)(s)   \, \|_{L^2_{w_\gamma}}^2 ds\right).
\end{split}
\end{equation*}
We set the quantity
\begin{equation*}
f(t_1)= \sup_{0 \leq t \leq t_1} \|  (\vv,\vw)(t) \,  \|_{L^2_{w_\gamma}}^2  + \frac{1}{2} \int_{0}^{t_1} \|     \vN (\vv, \vw)(s)   \, \|_{L^2_{w_\gamma}}^2 \,  ds. \\
\end{equation*}
Then, let $T_0 \in (t_0, T)$ such that for   $t_1\in (t_0,T_0)$  the quantities  $\|  \mathds{1}_{(t_0 , t_1)}\, \vu  \|_{X_T}$ and  $\|  \mathds{1}_{(t_0 , t_1)}\, \vb \|_{X_T}$ are small enough. Thus, for all $t_0 < t_1  < T_0$ we find:
\begin{equation*}
f(t_1) \leq  C \int_0^{t_1}   f(s) ds.
\end{equation*}
We conclude by the Gr\"onwall's lemma that $( \vu, \vb, p, q ) = ( \tilde \vu, \tilde \vb, \tilde p, \tilde q )$ on $[t_0, T_0]$. Finally,  as $t_0 \in [0 , T)$ is arbitrary then we have the identities $( \vu, \vb, p, q ) = ( \tilde \vu, \tilde \vb, \tilde p, \tilde q )$ on $[0,T)$. \Endproof{}

\section{Some results for the ($MHD^*$) system}\label{Priori-estim-stability}  
Our main theorem bases on the two following  results  for the  equations: 
\begin{equation*}\label{MHDa}
(\text{$MHD^{*}$}) \left\{ \begin{array}{ll}\vspace{2mm} 
\partial_t \vu = \Delta \vu - (\vv \cdot \nabla) \vu + (\vc \cdot \nabla) \vb - \nabla p + \vN \cdot \F,  \\ \vspace{2mm}
\partial_{t} \vb = \Delta \vb - (\vv \cdot \nabla) \vb + (\vc \cdot \nabla) \vu -\vN q + 
\vN  \cdot \G 
, \\ \vspace{2mm}
\vN \cdot  \vu=0, \,  \vN \cdot \vb =0, \\ \vspace{2mm}
\vu(0,\cdot)=\vu_0, \, \vb(0,\cdot)=\vb_0.   
\end{array}
\right.     
\end{equation*}
In this system,  we shall consider two cases for the functions $(\vv, \vc)$:
\begin{enumerate}
\item[$\bullet$] when we will consider the  (MHD) equations we take $(\vv,\vc)= (\vu, \vb)$.
\item[$\bullet$] when we will consider the regularized (MHD) equations we take  $(\vv,\vc) = (\vu* \theta_\epsilon, \vb* \theta_\epsilon)$, where, for $0 < \varepsilon < 1$ and for a   fixed,   non-negative and  radially  non increasing test function $\theta \in \mathcal{D}(\mathbb{R}^3)$  which is equals to 0 for $|x| \geq 1 $ and $\int \theta \, dx = 1$; we define $\theta_\varepsilon(x)=\frac{1}{\varepsilon^3}\theta (x / \varepsilon)$. 
\end{enumerate}

\subsection{A priori estimates}\label{Priori-estim}
\begin{Theorem}\label{estimatesB2}
 Let $0<T<+\infty$. Let  $\vu_{0}, \vb_{0} \in B_2$ be a divergence-free vector fields and let  $\mathbb{F},\mathbb{G} $ be  tensors  such that $\mathbb{F}, \mathbb{G} \in B_2 L^2(0,T)$. Moreover, let $(\vu, \vb, p, q) $ be a solution of the problem ($MHD^*$). \\ 

We suppose that:
 \begin{itemize} 
 \item[$\bullet$]  $\vu, \vb$ belongs to $L^\infty((0,T), B_2   )$ and $\vN\vu, \vN\vb$ belongs to $B_2   L^2 (0,T)$.
 \item[$\bullet$] The pressure $p$ and the term $q$ are related to $\vu$, $\vb$, $\mathbb{F}$ and $\G$ by 
 \begin{equation*}
     p  = \sum_{1 \leq i,j\leq 3} \mathcal{R}_i \mathcal{R}_j (v_{i} u_{j} - c_{i} b _{j} - F_{i,j})\,\,\, \text{and} \,\,\,
q  = \sum_{1 \leq i,j\leq 3} \mathcal{R}_i \mathcal{R}_j ( v_{i} b_{j} - c_{j} u _{i}- G_{ij}).   
 \end{equation*}

 \item[$\bullet$] The map $t\in [0,T)\mapsto \vu(t,\cdot)$ is $*$-weakly continuous from $[0,T)$ to $B_2 $, and for all compact set $K \subset \mathbb{R}^3$ we have:
 $$ \lim_{t\rightarrow 0}    \Vert  ( \vu(t,\cdot) - \vu_{0} , \vb(t,\cdot) - \vb_{0} ) \Vert_{L^{2}(K)} = 0.$$
 \item[$\bullet$] The solution $(\vu, \vb, p, q)$ is suitable : there exists a non-negative locally finite measure $\mu$ on $(0,T)\times\mathbb{R}^3$ such that
\begin{equation}\label{energloc1}
\begin{split}
   \partial_t(\frac {\vert\vu \vert^2 
    + |\vb|^2 }2)=&\Delta(\frac {\vert\vu \vert^2 + |\vb|^2 }2)-\vert\vN \vu \vert^2 - |\vN \vb|^2 - \vN\cdot\left( (\frac{\vert\vu \vert^2}2 + \frac{\vert\vb \vert^2}2) \vv  +p\vu   \right)\\
    &+ \vN \cdot \left( (\vu \cdot \vb) \vc + q \vb \right) + \vu \cdot(\vN\cdot\mathbb{F}) +\vb \cdot(\vN\cdot\mathbb{G})- \mu.
\end{split}     
\end{equation}
 \end{itemize}
Then, exists a constant $C \geq 1$, which  does not depend on $T$, and not on $\vu_0$, $\vb_0$ $\vu, \vb$, $\mathbb{F}$, $\mathbb{G}$ nor $\epsilon$, such that: 
\begin{enumerate}
    \item[$\bullet$] We have the following control on $[0,T)$:
 \begin{equation} \label{control}
 \begin{split}  
 &\max \{ \|(\vu, \vb )(t) \|_{B_2}^2 , \|\vN(\vu, \vb) \|_{B_2 L^2 (0,t)}^2  \} \\
 & \leq  C\| (\vu _0, \vb _0 ) \|^2_{B_{2}} + C\|  (\mathbb{F},\G) \|^2_{B_2 L^2 (0,t)}+ C \int_0^t \| (\vu, \vb) (s) \|^2_{B_{2 }} + \| (\vu, \vb) (s) \|^6_{B_{2 }} ds. \end{split}
\end{equation}
\item[$\bullet$] Moreover,  if $T_0<T$ is small enough: 
$$ C \left( 1 +\|(\vu_0, \vb_0)\|_{B_2}^2+  \|( \mathbb{F}, \G ) \|_{B_2 L^2 (0,T_0)}^2 \right)^2 \, T_0\leq 1, $$ then the following control respect to the data  holds: 
\begin{equation}\label{energ2}
\begin{split}
 \sup_{0\leq t\leq T_0} & \max \{ \|(\vu, \vb )(t,.)\|_{B_2}^2,  {  \|\vN(\vu, \vb)\|_{B_2 L^2 (0,t)}^2 } \} \\
&\leq C \left(1 +\|( \vu_0, \vb_0 )\|_{B_2}^2+  \| (\mathbb{F}, \G) \|_{B_2 L^2 (0,T_0)}^2 \right). \end{split}    
\end{equation} 
\end{enumerate}
\end{Theorem}
 \pv
 In this proof, we will focus only in the  case $(\vv,  \vc )= (\vu * \te ,\vb * \te  )$   (the case $(\vv, \vc)=(\vu, \vb)$  can be treated in a similar way). 
 
We start by proving the global control (\ref{control}). The idea is to apply the energy balance (\ref{energloc1}) to a suitable test function. Let $0<t_0<t_1<T$. We consider a function $\alpha_{\eta,t_0,t_1}$ which converges almost everywhere to $\mathds 1 _{[t_0, t_1]}$ and such that $\partial_t \alpha_{\eta,t_0,t_1}$ is the difference between two identity approximations, the first one in $t_0$ and the second one in $t_1$. For this, we take a non-decreasing function $
\alpha\in \mathcal{C}^\infty(\mathbb{R})$ which is equals to $0$ on $(-\infty,\frac{1}{2})$ and is equals to $1$ on $(1, +\infty)$. Then,   for $0<\eta< \min(\frac {t_0}2,T-t_1) $ we set the function $$\ds{ \alpha_{\eta,t_0,t_1}(t)=\alpha( \frac{t-t_0}\eta)-\alpha(\frac{t-t_1}\eta)}.$$
On the other hand, we consider  a non-negative function $\phi\in\mathcal{D}(\mathbb{R}^3)$ which is equals to $1$ for $\vert x\vert\leq 1/2$ and is equals to $0$ for $\vert x\vert\geq 1$; and for $R\geq 1$ we set
\begin{equation}
\label{varphidef}
    \ds{\phi_R(x)=\phi(\frac x R)}.
\end{equation}

Thus, by the energy balance (\ref{energloc1}) applied to test function $\alpha_{\eta,t_0,t_1}   \phi_R$, we can write
\begin{equation*} \begin{split}
  -\iint & \frac{\vert \vu \vert^2}2+ \frac{\vert \vb \vert^2}2  \partial_t\alpha_{\eta,t_0,t_1} \phi_R  \, dx\, ds + \iint \vert\vN\vu \vert^2 + \vert\vN\vb \vert^2 \, \,  \alpha_{\eta,t_0,t_1} \phi_R   dx\, ds
\\ \leq & \int \int \frac{ | \vu |^2 + | \vb |^2}{ 2 }  \,  \alpha_{\eta,t_0,t_1}  \Delta \phi_R \,  dx\, ds
 \\
 &+ \sum_{i=1}^3 \iint  [(\frac{\vert\vu \vert^2}2+ \frac{|\vb |^2}{2} ) v_{ i} + p u_{i} ] \alpha_{\eta,t_0,t_1}  \partial_i \phi_R \, dx\, ds
\\
 &+ \sum_{i=1}^3 \iint  [(\vu \cdot \vb) c_{ i}  + q b_{i} ] \alpha_{\eta,t_0,t_1}  \partial_i \phi_R \, dx\, ds 
\\&- \sum_{1 \leq i,j \leq 3} ( \iint  F_{i,j} u_{j}  \alpha_{\eta,t_0,t_1}  \partial_i \phi_R \, dx\, ds + \iint   F_{i,j}\partial_iu_{j}\  \alpha_{\eta,t_0,t_1} \phi_R \, dx\, ds )
\\&- \sum_{1 \leq i,j \leq 3} ( \iint  G_{i,j} b_{j}  \alpha_{\eta,t_0,t_1}  \partial_i \phi_R \, dx\, ds + \iint   G_{i,j}\partial_ib_{j}\  \alpha_{\eta,t_0,t_1} \phi_R \, dx\, ds ),
\end{split}\end{equation*}
and taking the limit  when $\eta$ goes to $0$, by the dominated convergence theorem we obtain  (when the limit in the left side is well-defined):
\begin{equation*} \begin{split}
  - \lim_{\eta\rightarrow 0} \iint & \frac{\vert \vu \vert^2}2+ \frac{\vert \vb \vert^2}2  \partial_t\alpha_{\eta,t_0,t_1} \phi_R  \, dx\, ds + \int_{t_0}^{t_1} \int \vert\vN\vu \vert^2 + \vert\vN\vb \vert^2 \, \,  \phi_R   dx\, ds
\\ \leq & \int_{t_0}^{t_1} \int \frac{ | \vu |^2 + | \vb |^2}{ 2 }  \,   \Delta \phi_R \, dx\, ds
 \\
 &+ \sum_{i=1}^3 \int_{t_0}^{t_1}  \int [(\frac{\vert\vu \vert^2}2+ \frac{|\vb |^2}{2} ) v_{ i} + p u_{i} ]   \partial_i \phi_R \, dx\, ds
\\
 &+ \sum_{i=1}^3 \int_{t_0}^{t_1} \int [(\vu \cdot \vb) c_{ i}  + q b_{i} ]   \partial_i \phi_R \, dx\, ds 
\\&- \sum_{1 \leq i,j \leq 3} ( \int_{t_0}^{t_1} \int F_{i,j} u_{j}    \partial_i \phi_R \, dx\, ds + \int_{t_0}^{t_1} \int   F_{i,j}\partial_iu_{j}\   \phi_R \, dx\, ds )
\\&- \sum_{1 \leq i,j \leq 3} ( \int_{t_0}^{t_1} \int G_{i,j} b_{j}   \partial_i \phi_R \, dx\, ds + \int_{t_0}^{t_1}  \int  G_{i,j}\partial_ib_{j}\   \phi_R \, dx\, ds ).
\end{split}\end{equation*}
We define now the quantity
  $$\ds{A_{R}(t)=\int ( \vert \vu(t,x)\vert^2  + |\vb(t,x)  \vert^2 )   \phi_R(x) \, dx},$$ hence,  if $t_0$ and $t_1$ are Lebesgue points of $A_{R}(t)$ and moreover,  due to the fact that 
$$-\iint (\frac{|\vu|^2}2 + \frac{ |\vb|^2}2 ) \partial_t\alpha_{\eta,t_0,t_1}   \phi_R \, dx\, ds=-\frac 1 2\int \partial_t\alpha_{\eta,t_0,t_1}
A_{R}(s) \, ds,$$
 we have  $$ \lim_{\eta\rightarrow 0}  -\iint ( \frac{\vert\vu\vert^2}2 + \frac{\vert\vb\vert^2}2 )  \partial_t\alpha_{\eta,t_0,t_1}  \phi_R \, dx\, ds=\frac 1 2 (  A_{R}(t_1)- A_{R}(t_0)).$$ 

Then, since $\phi_R$ is a support compact function  we can let $t_0$ go to $0$ and thus we can replace $t_0$ by $0$ in this inequality. Moreover, if we let $t_1$ go to $t$, then by the $*$-weak continuity  we have $ A_{R}(t)\leq \lim_{t_1\rightarrow t }  A_{R}(t_1)$, and thus  we may replace $t_1$ by $t\in (0,T)$. In this way, for every $t\in (0,T)$ we can write: 
\begin{equation}\label{localineq} \begin{split}
\int & \frac{\vert \vu(t,x)\vert^2 +  \vert \vb(t,x)\vert^2}2 \phi_R \, dx    +  \int_0^t  \int ( \vert\vN\vu\vert^2+  \vert\vN\vb\vert^2)\, \,  \phi_R \, ds\, dx
 \\ \leq & \int \frac{\vert \vu_0(x)\vert^2 + \vert \vb_0(x)\vert^2 }2 \phi_R \, dx   + \int_0^t  \int \frac{ | \vu |^2 + | \vb |^2}{ 2 }  \,   \Delta \phi_R \, ds\, dx
\\
 &+ \sum_{i=1}^3 \int_{0}^{t}  \int [(\frac{\vert\vu \vert^2}2+ \frac{|\vb |^2}{2} ) v_{ i} + p u_{i} ]   \partial_i \phi_R \, dx\, ds
\\
 &+ \sum_{i=1}^3 \int_{0}^{t} \int [(\vu \cdot \vb) c_{ i}  + q b_{i} ]   \partial_i \phi_R \, dx\, ds 
\\&- \sum_{1 \leq i,j \leq 3} ( \int_{0}^{t} \int F_{i,j} u_{j}    \partial_i \phi_R \, dx\, ds + \int_{0}^{t} \int   F_{i,j}\partial_iu_{j}\   \phi_R \, dx\, ds )
\\&- \sum_{1 \leq i,j \leq 3} ( \int_{0}^{t} \int G_{i,j} b_{j}   \partial_i \phi_R \, dx\, ds + \int_{0}^{t}  \int  G_{i,j}\partial_ib_{j}\   \phi_R \, dx\, ds ) .
\end{split}\end{equation}
In this inequality, we still need to estimate the terms in the right-hand side. For the second term, as $R\geq 1$ we write 
\begin{equation*}
\label{clb2}
    \frac{1}{R^2} \int (|\vu|^2+ |\vb|^2)  \Delta \phi_R  \, dx  \leq  \frac{ C }{R^4} \int_{B(0,R)} ( |\vu|^2+|\vb|^2) \, dx \leq  C ( \| \vu \|^2_{B_{2 }}+\| \vb \|^2_{B_{2 }}).
\end{equation*}
The third and fourth terms are estimates as follows. We consider first the expressions where the pressure terms $p$ and $q$ do not appear. Using the H\"older inequalities and the Sobolev embeddings we have: 
\begin{align*}
    \sum_{i=1}^3  \int &  \frac{ (\vu\cdot \vb ) }2 (b_i * \theta_\epsilon)      \partial_i \phi_R \,  dx  
    \leq \| \vu \|_{L^{\frac{12}{5}} (B(0,R)) }  \| \vb \|_{L^{\frac{12}{5}} (B(0,R)) }  \| \vb * \theta_\epsilon \|_{L^6 (B(0,R))} \| \nabla \phi_R \|_{ L^\infty } \nonumber\\ 
    & \leq  \frac{C}{R} \| \vu \|^{3/4}_{L^2 (B(0,R)) }   \| \vu \|^{1/4}_{L^6 (B(0,R)) } \| \vb \|^{3/4}_{L^2 (B(0,R)) }   \| \vb \|^{5/4}_{L^6 (B(0,R+1)) } \\
    & \leq  \frac{C}{R} \| \vb \|^{3/4}_{L^2 (B(0,R)) }   \| \vu \|^{3/4}_{L^2 (B(0,R)) }   U^{1/4} B^{5/4},
\end{align*}
where we have denoted the quantities 
$$ U= \left( \int |\phi_{2R} \nabla \vu|^2 dx \right)^{1/2} + \left(\int_{|x| \leq 2R} |\vu|^2 dx \right)^{1/2} $$
and
$$ B=  \left(\int |\phi_{2(R+1)} \nabla \vb|^2 dx \right)^{1/2} + \left(\int_{|x| \leq 2(R+1)} |\vb|^2 dx \right)^{1/2}.$$
Thus, we can write (by the Young's inequalities for products with  $1 = \frac{1}{8} + \frac{1}{8} + \frac{1}{8} + \frac{5}{8}$):
\begin{align*}
\label{cub2}
    \frac{1}{R^2} & \sum_{i=1}^3  \int  \frac{ (\vu\cdot \vb ) }2 (b_i * \theta_\epsilon)      \partial_i \phi_R \,  dx  \nonumber \\
    & \leq C  (  \frac{\| \vu \|_{L^2 (B(0,R)) }}{R } )^{3/4}  (  \frac{\| \vb \|_{L^2 (B(0,R)) }}{R } )^{3/4}  ( \frac{U  }{R } )^{1/4}  ( \frac{ B }{R } )^{5/4} \nonumber\\
    & \leq C \| (\vu, \vb) \|_{B_{2 }}^6 
    +  C \| (\vu, \vb)  \|_{B_{2 }}^2  +  \frac{C_0}{R^2} \int |\phi_{2R} \nabla \vu|^2 + |\phi_{2(R+1)} \nabla \vb|^2  \, dx
\end{align*}
where $C_0>0$ is an arbitrarily small constant.\\

Now, in order to estimate  the expressions where the pressure terms $p$ and $q$ appear, we need the following technical lemma which will be proved at the end of this section.  
\begin{Lemma}
\label{lemmat} Within the hypothesis of Theorem  \ref{estimatesB2},  
the terms  $p$ and $q$ belong $  L^{3/2}_{\rm loc} $. Moreover,  there exist an arbitrarily small constant $C_0>0$ and a constant $C>0$,  which do not depend on $T$, $\vu$, $\vb$, $\vu_0$, $\vb_0$, $\mathbb{F}$, $\mathbb{G}$ nor $\epsilon$; such that for all $R \geq 1$ and for all  $0\leq t \leq T$ we have: 
\begin{align*}
     \frac{1}{R^2} & \sum_{i=1}^3 \int_0^t    \int     (p u_i +  q b_i) \,  \partial_i \phi_R \, ds\, dx   \\
     \leq &    C \|(\mathbb{F}, \G)\|_{B_2 L^2 (0,t)}^2   + C \int_0^t  \| (\vu, \vb)(s) \|_{B_{2 }}^2 + \| (\vu, \vb)(s) \|_{B_{2 }}^6   \\ 
     & 
    + \frac{C_0}{R^2} \int \int_0^t  |\varphi_{2(5R+1)} \nabla \vu|^2 + |\varphi_{2(5R+1)} \nabla \vb|^2 \, dx.
\end{align*} 
\end{Lemma}
Finally, the fifth and sixth terms (which involve the tensor forces $\mathbb{F}$ and $\mathbb{G}$) are easily estimate as follows. We will write down  only the estimates for  $\mathbb{F}$ since the estimates for   $\mathbb{G}$ are completely similar:  

  \begin{equation*}\begin{split} \left\vert \frac{1}{R^2} \sum_{1 \leq  i,j \leq 3} \int_0^t\int F_{i,j} (\partial_i u_j) \phi_R \, dx\, ds \right\vert\leq& \,  C  \|\mathbb{F}\|_{B_2 L^2 (0,t)}^2  +  \frac{C_0}{R^2} \int_0^t \int_{|x|<R}  |\vN \vu|^2 \, dx \, ds , 
 \end{split}\end{equation*}
and 
   \begin{equation*}\begin{split} \left\vert \frac{1}{R^2} \sum_{1 \leq  i,j \leq 3} \int_0^t\int F_{i,j} u_i \partial_j(\phi_R)\, dx\, ds\right\vert\leq &  \, C \|\mathbb{F}\|_{B_2 L^2 (0,t)}^2 + C \int_0^t \|\vu(s )\|_{B_2}^2\, ds  . 
 \end{split}\end{equation*}
 where $C_0>0$ always denote a small enough constant. \\
 
 Once we dispose of all these estimates, we are able to write 

 \begin{equation*} \begin{split}  \int & ( \frac{\vert \vu(t,x)\vert^2}2 + \frac{\vert \vb(t,x)\vert^2}2 ) \phi_R \, dx     +  \int_0^t  \int  (\vert\vN\vu\vert^2 + \vert\vN\vb\vert^2 )\, \,  \phi_R \, ds \, dx  \\ \leq& \int (\frac{\vert \vu(0,x)\vert^2}2 + \frac{\vert \vb(0,x)\vert^2}2 ) \phi_R \, dx   + C \|( \mathbb{F}, \G )\|_{B_2 L^2 (0,t)}^2 \, ds \\ & + C  \int_0^t \| (\vu, \vb) (s, \cdot) \|^2_{B_{2 }} + \| (\vu, \vb ) (s, \cdot) \|^6_{B_{2 }} \, ds  \\
 & + \frac{C_0}{R^2} \int \int_0^t  |\varphi_{2(5R+1)} \nabla \vu|^2 + |\varphi_{2(5R+1)} \nabla \vb|^2 \, dx, 
 \end{split}\end{equation*} 
where the desired energy control  \eqref{control} follows. To finish this proof,  the estimate (\ref{energ2})  follows directly from \eqref{control} and the Lemma $3.1$ in \cite{PF_OJ} (see the proof of Corollary $3.3$, page 17, for all the details).  
\Endproof
\\

\textbf{Proof of Lemma \ref{lemmat}.} As in the proof of the theorem above, we will consider only the case $(\vv, \vc)=(\vu* \te ,\vb* \te  )$. Moreover, we will focus only on the expression which involves the pressure $p$, since the computations for the other expression, where the term $q$ appears, are completely similar. \\  

We write $\ds{ \frac{1}{R^2} \sum_{k=1}^3 \int_0^t \int_{\vert x \vert \leq R} \vert     p u_k \vert \vert \  \partial_k \phi_R \vert \,  dx \, ds\leq   \frac{c}{R^3} \sum_{k=1}^3 \int_0^t \int_{\vert x \vert \leq R} \vert     p u_k \vert  \,  dx \, ds}$, 
and recalling that   $p  = \sum_{1 \leq i,j\leq 3} \mathcal{R}_i \mathcal{R}_j ( (u_{i}*\te) u_{j} - (b_{i}*\te)  b _{j} - F_{i,j} )$, the last expression allow us to write
\begin{equation*}
\begin{split}
 \frac{1}{R^2}& \sum_{k=1}^3 \int_0^t \int_{\vert x \vert \leq R} \vert     p u_k \vert \vert \  \partial_k \phi_R \vert  \,  dx \, ds  \\
 \leq& \frac{c}{R^3} \sum_{k=1}^3 \int_0^t \int_{\vert x \vert \leq R} \vert u_k   \sum_{i,j=1}^3  \mathcal{R}_i \mathcal{R}_j ( (u_{i}*\te) u_{j} )  \vert  \,  dx \, ds \\
 &+ \frac{c}{R^3} \sum_{k=1}^3 \int_0^t \int_{\vert x \vert \leq R} \vert u_k  \sum_{i,j=1}^3 \mathcal{R}_i \mathcal{R}_j ( (b_{i}*\te)  b _{j} - F_{i,j}  )  \vert  \,  dx \, ds,   
\end{split}    
\end{equation*}
 and since we have the same information on $\vu$ and $\vb$ it is enough to study the last term above. For $R\geq 1$ we define the following expressions:  
 $$  p_1  = \sum_{i,j} \mathcal{R }_i \mathcal{R}_j   ( \mathds{1}_{|y|< 5R}  (\theta_{\epsilon} * b_i) b_j ), \,\,\, p_2  = - \sum_{i,j} \mathcal{R }_i \mathcal{R}_j  (  \mathds{1}_{|y|\geq 5R}  (\theta_{\epsilon} * b_i) b_j ),$$ and 
 $$   p_3  = - \sum_{i,j} \mathcal{R }_i \mathcal{R}_j   ( \mathds{1}_{|y|< 5R}  F_{i,j} ), \,\,\,  p_4  = \sum_{i,j} \mathcal{R }_i \mathcal{R}_j  (  \mathds{1}_{|y|\geq 5R}  F_{i,j}  ),$$ and then, by the Young's inequalities  (for products), we have
\begin{equation*}
\begin{split}
      \frac{c}{R^3}& \sum_{k=1}^3 \int_0^t \int_{\vert x \vert \leq R} \vert u_k \sum_{i,j=1}^3   \mathcal{R}_i \mathcal{R}_j ( (b_{i}*\te)  b _{j} - F_{i,j}  )  \vert  \,  dx \, ds \\
  &\leq \frac{C}{R^3} \int_0^t  \int_{|x|\leq R}(|p_1|^{3/2} + |p_2|^{3/2} + | \vu|^3 +  |p_3|^{2} + |p_4|^{2}+|\vu|^2) dx\, ds ,
\end{split}
\end{equation*}
where we will study each term separately. \\

To study $p_1$, by the continuity of $\mathcal{R}_i$ on $L^{\frac{3}{2}}(\Rt)$, since the test function $\theta_\varepsilon$ verifies $\int \theta_\varepsilon(x) dx = 1$ and $supp (\theta_\epsilon)  \subset \overline{B(0,1)}$ and moreover, by the Fubini's theorem   we can write 

\begin{align*}
\label{cpb2}
\int_{|x|\leq R} |p_1|^{3/2}  dx
    &\leq C \int |p_1|^{3/2}  dx  \leq   C \int | ( \mathds{1}_{|x|< 5R} (\theta_{\epsilon} * \vb) \otimes \vb)|^{3/2}  dx  \nonumber \\
    &\leq C \left( \int |  \mathds{1}_{|x|< 5R} (\theta_{\epsilon} * \vb) |^{3}  dx\right)^{1/2} \left( \int |  \mathds{1}_{|y|< 5R}  \vb |^{3} \, dx  \right)^{1/2}  \nonumber  \\
    &\leq C \left( \int_{|x|\leq 5R}  \int_{|x-z|\leq 1 }  \theta_{\epsilon}(x-z) |\vb(z)|^{3}  dz  dx\right)^{1/2} \left( \int | ( \mathds{1}_{|y|< 5R}  \vb)|^{3} \, dx  \right)^{1/2}  \nonumber \\
    &\leq C \left( \int_{|x|\leq 5R}   \int_{|z|\leq 5R+1 }  \theta_{\epsilon}(x-z) |\vb(z)|^3 dz   dx\right)^{1/2} \left( \int | ( \mathds{1}_{|y|< 5R}  \vb)|^{3} \, dx  \right)^{1/2}  \nonumber\\
    &\leq C \int_{|z|\leq 5R+1} |\vb|^3 \, dz. 
\end{align*}
With this estimate at hand, we see that
\begin{equation*}
    \int_{|x|\leq R} |\vu|^3   + |p_1|^{3/2} \, dx 
   \leq C \int_{|x|\leq 5R+1} |\vu|^3+|\vb|^3  \, dx,
 \end{equation*}
and using the Sobolev embedding   we  write  
 \begin{align*}
    &\frac{C}{R^3} \int_{|x|\leq 5R+1} |\vu|^3  \, dx  \leq  \frac{C}{R^3} \| \vu \|^{3/2}_{L^2 (B(0,5R+1)) }   \| \vu \|^{3/2}_{L^6 (B(0,5R+1)) } \\
    & \leq \frac{C}{R^{3/2}} \| \vu \|^{3/2}_{L^2 (B(0,5R+1)) }    \left(  \left( \frac{1}{R^2} \int |\phi_{2(5R+1)} \nabla  \vu|^2 dx \right)^{1/2} + \left( \frac{1}{R^2} \int_{|x| \leq 2(5R+1)} |\vu|^2 dx \right)^{1/2}  \right)^{3/2} \nonumber\\
    & \leq C \| \vu \|_{B_{2 }}^6 
    +  C_0 \| \vu  \|_{B_{2 }}^2  +   \frac{C_0}{R^2} \int |\phi_{2(5R+1)} \nabla  \vu|^2 dx,
\end{align*}
where $C_0>0$ is a arbitrarily small constant. Similar bounds works for $\vb$. \\

We study now the term  $p_2$. Remark first  that there exist a constant $C>0$ (which does not depend on $R>1$) such that for all $|x|\leq R$ and all $|y|\geq 5R$, the kernel $\mathbb{K}_{i,j}$ of the  operator $\mathcal{R}_i \mathcal{R}_j$ verifies
$|\mathbb{K}_{i,j} (x-y)|  \leq \frac{C }{|y|^3}$ (see \cite{Gr09} for a proof)  and then we write:
\begin{align*}
   &\left(\int_{|x|\leq R} |p_2|^{3/2}  dx \right)^{2/3}\\
   &\leq C \sum_{i,j} \left(\int_{|x|\leq R}  \left( \int  | \mathbb{K}_{i,j} (x-y) |  \, | (\theta_{\epsilon} * b_i)(y)  b_j(y) | \,  \mathds{1}_{|y|\geq 5R} \, dy \right)^{3/2}  dx \right)^{2/3} \\
   &\leq C \left( \int_{|x|\leq R}   \left( \int_{|y|\geq 5R}  \frac{1}{|y|^3} | (\theta_{\epsilon} * \vb) \otimes \vb |  \, dy \right)^{3/2}  dx  \right)^{2/3}\\
   &\leq C R^2 \int_{|y|\geq 5R}     \frac{1}{|y|^3} | (\theta_{\epsilon} * \vb) \otimes \vb |   dy  \\
   & \leq C R^2 \left(\int_{|y|\geq 5R}     \frac{1}{|y|^3} | \theta_{\epsilon} * \vb |^2   dy \right)^{1/2} \left( \int_{|y|\geq 5R}     \frac{1}{|y|^3} | \vb |^2   dy \right)^{1/2}    \\
    & \leq C R^2 \left(\int_{|y|\geq 5R}     \frac{1}{|y|^3} \int_{|y-z|<1} \theta_{\epsilon}(y-z) | \vb (z) |^2 dz\,  dy \right)^{1/2} \left( \int_{|y|\geq 5R}     \frac{1}{|y|^3} | \vb |^2   dy \right)^{1/2}    \\
    & \leq C R^2 \left(\int_{|y|\geq 5R}      \int_{|z|\geq 5R-1} \frac{1}{|z|^3} \theta_{\epsilon}(y-z) | \vb (z) |^2 dz\,  dy \right)^{1/2} \left( \int_{|y|\geq 5R}     \frac{1}{|y|^3} | \vb |^2   dy \right)^{1/2}    \\
   & \leq C R^2  \int_{|z|\geq 5R-1}     \frac{1}{|z|^3} | \vb |^2   dz.    \\
\end{align*}
With this estimate, and the fact that $B_2(\Rt) \subset L^{2}_{w_3}(\Rt)$, we finally obtain
\begin{equation*}
    \frac{ C}{R^3} \int_{|y| \leq R} |p_2|^{3/2}  dx \leq C \left( \int \frac{1}{(1+|z|)^3} |\vb|^{2} \right)^{3/2}
    \leq C \| \vb \|^3_{B_2}.
\end{equation*}
It remains to estimate the terms $p_3$ and $p_4$ which involve the tensor $\mathbb{F}$. For $p_3$, using the continuity of the Riesz transform  $\mathcal{R}_i $ on $L^2$, we obtain directly: 
\begin{equation*}
    \frac{ c }{R^3}   \int_0^t \int_{|x|\leq R} |p_3|^{2}  dx\,ds
    \leq  \frac{ C }{R^3} \sum_{i,j} \int_0^t \int_{|x|< 5R} |  \mathbb{F}_{i,j}|^{2} dx\, ds 
    \leq C \|  \mathbb{F} \|_{B_2 L^2 (0,t)}^2.
\end{equation*}
For the term $p_4$,  remark first  that we have 
\begin{align*}
   \left(\int_{|x|\leq R} |p_4|^{2}  dx \right)^{1/2}
   &
   C \leq \sum_{i,j} \left(\int_{|x|\leq R}  \left( \int_{|y|\geq 5R}  | \mathbb{K}_{i,j}(x-y)    \mathbb{F}_{i,j} | \, dy \right)^{2}  dx \right)^{1/2} \\
   &\leq C \sum_{i,j} \left( \int_{|x|\leq R}   \left( \int_{|y|\geq 5R}  \frac{1}{|y|^3} | \mathbb{F}_{i,j}  | \, dy \right)^{2}  dx  \right)^{1/2} \\
   & \leq C \sum_{i,j}  R^{3/2} \int_{|y|\geq 5R}     \frac{1}{|y|^3} | \mathbb{F}_{i,j} | \,  dy ,
\end{align*}
and then, for $0<\delta < 1 $,  and by the H\"older inequalities  we can write:
\begin{align*}
    \frac{C}{R^3} \int_0^t \int_{|x| \leq R} |p_4|^{2}  \, dx \, ds 
    & \leq  C \sum_{i,j}  \int_0^t \left( \int \frac{1}{(1+|x|)^3} | \mathbb{F}_{i,j} | \, dx \right)^{2} \, ds \\
    & \leq  C  \sum_{i,j} \int_0^t \int \frac{1}{(1+|x|)^{2 + \delta}} | \mathbb{F}_{i,j} |^2 \, dx  \, ds \\
    &  \leq    C  \sum_{i,j} \int \frac{1}{(1+|x|)^{2 + \delta}} \int_0^t | \mathbb{F}_{i,j}  |^2  \, ds \, dx  \\
    & \leq   C \| \mathbb{F} \|^2_{B_2 L^2(0,t)} .
\end{align*}
 The lemma is proven.
\Endproof{}
\subsection{A stability result}
\begin{Theorem}\label{stability}
 Let $0<T<+\infty$. Let  $\vu_{0,n}, \vb_{0,n}$ be  divergence-free vector fields such that $(\vu_{0,n}, \vb_{0,n})\in B_2$. Let    $\mathbb{F}_n$ and $\G_n$ be    tensors  such that $(\mathbb{F}_n, \G_n ) \in B_2 L^2(0,T)$. Let $(\vu_n, \, \vb_n, \, p_n\, q_n ) $ be a solution  of the $(MHD^{*})$  problem: 
 \begin{equation}\label{ADn}
  \left\{ \begin{array}{ll}\vspace{2mm} 
 \partial_t \vu_n = \Delta \vu_n - (\vv_n \cdot \nabla) \vu_n + (\vc_n \cdot \nabla) \vb_n - \nabla p_n + \vN \cdot \F_n,  \\ \vspace{2mm}
 \partial_{t} \vb_n = \Delta \vb_n - (\vv_n \cdot \nabla) \vb_n + (\vc_n \cdot \nabla) \vu_n -\vN q_n + 
 \vN  \cdot \G_n 
 , \\ \vspace{2mm}
 \vN \cdot  \vu_n=0, \, \,   \vN \cdot \vb_n =0, \\ \vspace{2mm}
 \vu_n(0,\cdot)=\vu_{0,n}, \,\, \vb_n(0,\cdot)=\vb_{0,n}.  
 \end{array}
 \right.     
 \end{equation}
  which verifies the same hypothesis of Theorem \ref{estimatesB2}. \\  

If $(\vu_{0,n}, \, \vb_{0,n})$ is strongly convergent to $(\vu_{0,\infty}, \, \vb_{0,\infty})$ in  $B_2$, and  if the sequence $(\mathbb{F}_n, \, \G_n )$ is strongly convergent to $(\mathbb{F}_\infty, \, \G_\infty)$  in $B_2L^2(0,T)$; then there exists $(\vu_\infty, \vb_\infty, p_\infty, q_\infty)$   and an increasing  sequence $(n_k)_{k\in\mathbb{N}}$ with values in $\mathbb{N}$ such that:
  \begin{itemize} 
 \item[$\bullet$] $(\vu_{n_k}, \vb_{n_k})$ converges *-weakly to $(\vu_\infty, \vb_\infty)$ in $L^\infty((0,T), B_2)$, $(\vN\vu_{n_k}, \vN\vb_{n_k})$ converges weakly to $(\vN\vu_\infty, \vN\vb_\infty) $ in $B_2L^2(0,T)$.
 \item[$\bullet$] $(\vu_{n_k}, \vb_{n_k})$ converges strongly  to $(\vu_\infty, \vb_\infty)$ in  $L^2_{\rm loc}([0,T)\times\mathbb{R}^3)$. 
 \item[$\bullet$]For $2<\gamma < 5/2$, the sequence $(p_{n_k},q_{n_k}) $converges weakly to $(p_\infty, q_\infty) $ in  $L^3((0,T), L^{6/5}_{w_{\frac{6 \gamma}{5}}}) + L^2((0,T), L^{2}_{w_{\gamma}})$.
\end{itemize}
 
 Moreover, $(\vu_\infty, \vb_\infty, p_\infty, q_\infty)$ is  a solution of the  problem $(MHD^{*})$:   
 \begin{equation}\label{ADinf}
 \left\{ \begin{array}{ll}\vspace{2mm} 
\partial_t \vu_\infty = \Delta \vu_\infty - (\vu_\infty \cdot \nabla) \vu_\infty + (\vb_\infty \cdot \nabla) \vb_\infty - \nabla p_\infty + \vN \cdot \F_\infty,  \\ \vspace{2mm}
\partial_{t} \vb_\infty = \Delta \vb_\infty - (\vu_\infty \cdot \nabla) \vb_\infty + (\vb_\infty \cdot \nabla) \vu_\infty -\vN q_\infty + 
\vN  \cdot \G_\infty 
, \\ \vspace{2mm}
\vN \cdot  \vu_\infty=0, \, \,   \vN \cdot \vb_\infty =0, \\ \vspace{2mm}
\vu_\infty(0,\cdot)=\vu_{0,\infty}, \,\, \vb_\infty(0,\cdot)=\vb_{0,\infty},  
\end{array} 
\right.     
\end{equation}
and verifies all the hypothesis of Theorem \ref{estimatesB2}.
\end{Theorem}
\pv 
We will verify that the sequence $(\vu_n, \vb_n)$ satisfy the hypothesis of the Rellich lemma (see Lemma 6 in \cite{PF_PG}). Remark first that: since for $2<\gamma$ we have that  $\vu_n, \vb_n$ is bounded in $L^\infty((0,T), B_2) \subset L^\infty((0,T), L^2_{w_\gamma} ) $ and moreover, since we have that  $\vN\vu_n, \vN\vb_n$ is bounded in $B_2   L^2 (0,T) \subset L^2((0,T), L^2_{w_\gamma} ) $, then  for all $\varphi\in \mathcal{D}(\mathbb{R}^3)$ we have that $(\varphi \vu_n, \varphi \vb_n)$ are bounded in  $L^2((0,T), H^1)$. On the other hand, for the pressure $p_n$ and the term $q_n$ we write 
$p_n=p_{n,1}+ p_{n,2}$  with
$$ p_{n,1} =  \sum_{i=1}^3\sum_{j=1}^3 R_iR_j(v_{n,i}u_{n,j} - c_{n,i} b_{n,j}), \quad   p_{n,2}=-\sum_{i=1}^3\sum_{j=1}^3 R_iR_j(F_{n,i,j}),$$
and we write $q_n=q_{n,1}+q_{n,2}$ with
$$ q_{n,1} =  \sum_{i=1}^3\sum_{j=1}^3 R_iR_j(v_{n,i}b_{n,j} - c_{n,i} u_{n,j}), \quad   q_{n,2}=-\sum_{i=1}^3\sum_{j=1}^3 R_iR_j(G_{n,i,j}).$$  From now on we fix $\gamma \in(2,\frac{5}{2})$, and  using the interpolation inequalities and the continuity of the Riesz transforms in the Lebesgue weighted spaces we 
get that the sequence  $(p_{n,1}, q_{n,1})$ is bounded in   $L^{3}((0,T),L^{6/5}_{w_{\frac {6\gamma}5}})$. Indeed, for the term $p_{n,1}$ recall that by Lemma \ref{Lemmuck} we have that for $0<\gamma < 5/2$ the weight $w_{6\gamma/5 }$ belongs to the Muckenhoupt class $\mathcal{A}_{p}(\Rt)$ (with $1<p<+\infty$) and then we can write: 
 \begin{align*}
        \label{productbu}
        \|  \sum_{i,j} \mathcal{R }_i \mathcal{R}_j ( \vu_{n,i}  \vu_{n,j} ) w_\gamma   \|_{ L^{6/5}} & \leq
       \|  ( \vu_n \otimes \vu_n )  w_\gamma   \|_{ L^{6/5}}  \leq \| \sqrt{w_\gamma} \vu_n \|^{\frac{3}{2}}_{ L^{2} } \| \sqrt{w_\gamma} \vu_n \|^{\frac{1}{2}}_{ L^{6} } \\
       &\leq \| \sqrt{w_\gamma} \vu \|^{\frac{3}{2}}_{ L^{2} } ( \| \sqrt{w_\gamma} \vu \|_{ L^{2} } +  \| \sqrt{w_\gamma} \vN \vu \|_{ L^{2} })^{\frac{1}{2}}. 
  \end{align*} 
The term $q_{n,1}$ is estimated in a similar way. Moreover we have that the sequence and $(p_{n,2}, q_{n,2})$ is bounded in $L^{2}((0,T),L^{2}_{w_\gamma})$. With these information, by equation (\ref{ADn}) we obtain that $(\varphi \partial_t \vu_n, \varphi \partial_t \vb_n)$ are bounded in the space  $L^2 L^2 + L^2 W^{-1,6/5}+ L^2 H^{-1}\subset L^2((0,T), H^{-2})$. Thus, we can apply  the Rellich lemma and there exists an increasing  sequence $(n_k)_{k\in\mathbb{N}}$ in $\mathbb{N}$, and there exist a couple of functions $(\vu_\infty, \vb_\infty)$   such that $(\vu_{n_k},\vb_{n_k}) $ converges strongly  to $(\vu_\infty, \vb_\infty)$ in  $L^2_{\rm loc}([0,T)\times\mathbb{R}^3)$.
We also have that $(\vv_{n_k},\vc_{n_k}) = (\vv_{n_k} * \theta_{\epsilon_{n_k}}, \vc_{n_k} * \theta_{\epsilon_{n_k}}) $ converges strongly  to $(\vu_\infty, \vb_\infty)$ in  $L^2_{\rm loc}([0,T)\times\mathbb{R}^3).$\\ 
 
 As $(\vu_n, \vb_n)$ are bounded in $L^\infty((0,T), L^2_{w_\gamma})$ and $(\vN\vu_n, \vN\vb_n)$ are bounded  in $L^2((0,T), L^2_{w_\gamma})$, we have  that  $(\vu_{n_k},\vb_{n_k}) $ converges *-weakly to $(\vu_\infty, \vb_\infty) $ in $L^\infty((0,T), L^2_{w_\gamma})$, and  $(\vN\vu_{n_k},\vN\vb_{n_k})$ converges weakly to $(\vN\vu_\infty, \vN\vb_\infty)$ in $L^2((0,T),L^2_{w_\gamma})$. Moreover, by the Sobolev embeddings and the interpolation inequalities we have that  $(\vu_{n_k}, \vb_{n_k})$ converges weakly to $(\vu_\infty, \vb_\infty) $ in $L^3((0,T), L^3_{w_{3\gamma/2}})$.
Also $(\vv_{n_k}, \vc_{n_k}) = (\vv_{n_k} * \theta_{\epsilon_{n_k}}, \vc_{n_k} * \theta_{\epsilon_{n_k}})$ converges weakly to $(\vu_\infty, \vb_\infty) $ in $L^3((0,T), L^3_{w_{3\gamma/2}})$, since it  is bounded in $L^3((0,T), L^3_{w_{3\gamma/2}})$. In particular, we may observe that  the  terms   $v_{n_k, i} u_{n_k, j}$, $c_{n_k, i} b_{n_k, j}$, $v_{n_k, i} b_{n_k, j}$ and $c_{n_k, i} u_{n_k, j}$ are weakly convergent in $(L^{6/5}L^{6/5})_{\rm loc}$ and thus in $\mathcal{D}'((0,T)\times \mathbb{R}^3)$. \\

As those terms are bounded in  $L^{3}((0,T),L^{6/5}_{w_{\frac {6\gamma}5}})$, they are weakly convergent in  $L^{3}((0,T),L^{6/5}_{w_{\frac {6\gamma}5}})$; and defining $p_\infty= p_{\infty,1}+p_{\infty,2}$ with  $$ p_{\infty,1} =  \sum_{i=1}^3\sum_{j=1}^3 R_iR_j(v_{\infty,i}u_{\infty,j} - c_{\infty,i} b_{\infty,j}), \quad   p_2=-\sum_{i=1}^3\sum_{j=1}^3 R_iR_j(F_{\infty,i,j}),$$
and $q_\infty= q_{\infty,1}+ q_{\infty,2}$ with 
$$ q_{\infty,1} =  \sum_{i=1}^3\sum_{j=1}^3 R_iR_j(v_{\infty,i}b_{\infty,j} - c_{\infty,i} u_{\infty,j}), \quad   q_2=-\sum_{i=1}^3\sum_{j=1}^3 R_iR_j(G_{\infty,i,j}),$$  we obtain  that $(p_{n_k,1}, q_{n_k,1}) $ are weakly convergent in  $L^{3}((0,T),L^{6/5}_{w_{\frac {6\gamma}5}})$  to $(p_{\infty,1}, q_{\infty,1}) $, and moreover, we get that  $(p_{n_k,2}, q_{n_k,2})$  is strongly convergent in  $L^2((0,T),L^2_{w_\gamma})$  to $(p_{\infty,2}, q_{\infty,2})$. So, we have that $(\vu_\infty, \vb_\infty, p_\infty, q_\infty)$ verify  the three first equations in the system $(MHD^{*})$ in $\mathcal{D}'((0,T)\times\mathbb{R}^3)$.
\\

It remains to verify the conditions at the time $t=0$. Remark that  $(\partial_t\vu_\infty, \partial_t\vb_\infty)$ are locally in $L^2 H^{-2}$, and then $(\vu_\infty, \vb_\infty)$ have representatives such that  $t\mapsto (\vu_\infty(t,.), \vb_\infty(t,.))$ is continuous from $[0,T)$ to $\mathcal{D}'(\mathbb{R}^3)$ (hence *-weakly continuous from $[0,T)$ to $B_2$) and moreover, they coincide with $\vu_\infty(0,.)+\int_0^t \partial_t \vu_\infty\, ds$ and $\vb_\infty(0,.)+\int_0^t \partial_t \vb_\infty\, ds$. Thus, in $\mathcal{D}'((0,T)\times \mathbb{R}^3)$, we have that
\begin{align*}
&\vu_\infty(0,.)+\int_0^t \partial_t \vu_\infty\, ds  =\vu_\infty
     =\lim_{n_k\rightarrow +\infty} \vu_{n_k} =\lim_{n_k\rightarrow +\infty}
\vu_{n_k,0}+ \int_0^t \partial_t \vu_{n_k}\, ds\\
&=\vu_{\infty,0}+\int_0^t \partial_t\vu_\infty\, ds,
\end{align*}
which implies that $\vu_\infty(0,.)=\vu_{\infty,0}$. Similar we have the identity $\vb_\infty(0,.)=\vb_{\infty,0}$. We conclude that  $(\vu_\infty, \vb_\infty, p_\infty, q_\infty)$ is a solution of the $(MHD^{*})$ equations.\\ 

Our next task is to verify the local energy balance. We define the quantity 
\begin{equation*}
\begin{split}
A_{n_k} =& - \partial_t(\frac {\vert\vu_{n_k} \vert^2 
 	+ |\vb_{n_k}|^2 }2)+ \Delta(\frac {\vert\vu_{n_k}\vert^2 + |\vb_{n_k}|^2 }2) -  \vN\cdot\left( (\frac{\vert\vu_{n_k} \vert^2}2 + \frac{\vert\vb_{n_k} \vert^2}2)\vv_{n_k} \right)\\
 	& -\vN \cdot(p_{n_k} \vu_{n_k}) - \vN \cdot(q_{n_k} \vb_{n_k})+ \vN \cdot ((\vu_{n_k} \cdot \vb_{n_k}) \vc_{n_k}) \\
 &+ \vu_{n_k} \cdot(\vN\cdot\mathbb{F}_{{n_k}}) +\vb_{n_k} \cdot(\vN\cdot\mathbb{G}_{{n_k}}).
\end{split}    
\end{equation*}
Remark that by the information on $(\vu_n, \vb_n)$ and by interpolation we have $(\vu_n, \vb_n)$ are bounded in  $L^{10/3}((0,T), L^{10/3}_{w_{5\gamma/3}})$ and then $(\vu_{n_k}, \vb_{n_k})$ are  locally bounded in $L^{10/3}_{t}L^{10/3}_{x}$ and locally strongly convergent in $L^{2}_{t} L^{2}_{x}$. So, $(\vu_{n_k}, \vb_{n_k})$  converges strongly  in $(L^{3}_{t} L^{3}_{x})_{loc}$. Moreover, by Lemma \ref{lemmat} we have that $(p_{n_k}, q_{n_k})$ are locally bounded in $L^{3/2}_{t}L^{3/2}_{x}$. Thus the quantity $A_{n_k}$ converges in the distributional sense to 
\begin{equation*}
\begin{split}
A_{\infty} =& - \partial_t(\frac {\vert\vu_{\infty} \vert^2 
 	+ |\vb_{\infty}|^2 }2)+ \Delta(\frac {\vert\vu_{\infty}\vert^2 + |\vb_{\infty}|^2 }2) -  \vN\cdot\left( (\frac{\vert\vu_{\infty} \vert^2}2 + \frac{\vert\vb_{\infty} \vert^2}2)\vv_{\infty} \right)\\
 	&-\vN \cdot(p_{\infty} \vu_{\infty}) - \vN \cdot(q_{\infty} \vb_{\infty})+ \vN \cdot ((\vu_{\infty} \cdot \vb_{\infty}) \vc_{\infty}) \\
 &+ \vu_{\infty} \cdot(\vN\cdot\mathbb{F}_{{\infty}}) +\vb_{\infty} \cdot(\vN\cdot\mathbb{G}_{{\infty}}). 
\end{split}    
\end{equation*}
Moreover, recall that  by hypothesis of this theorem there exist $\mu_{n_k}$ a non-negative locally finite measure on $(0,T)\times \mathbb{R}^3$ such that 
\begin{equation*}
 \begin{split}
& \partial_t(\frac {\vert\vu_{n_k} \vert^2 
 	+ |\vb_{n_k}|^2 }2)= \Delta(\frac {\vert\vu_{n_k}\vert^2 + |\vb_{n_k}|^2 }2)-\vert\vN \vu_{n_k} \vert^2 - |\vN \vb_{n_k}|^2 \\
& - \vN\cdot\left( (\frac{\vert\vu_{n_k} \vert^2}2 + \frac{\vert\vb_{n_k} \vert^2}2)\vv_{n_k} \right) -\vN \cdot(p_{n_k} \vu_{n_k}) - \vN \cdot(q_{n_k} \vb_{n_k}) \\
&+ \vN \cdot ((\vu_{n_k} \cdot \vb_{n_k}) \vc_{n_k})  + \vu_{n_k} \cdot(\vN\cdot\mathbb{F}_{{n_k}}) +\vb_{n_k} \cdot(\vN\cdot\mathbb{G}_{n_k})- \mu_{n_k}.
 \end{split}     
 \end{equation*}
Then, by definition of $A_{n_k}$ we can write $ A_{n_k}= \vert\vN \vu_{n_k} \vert^2 + |\vN \vb_{n_k}|^2 + \mu_{n_k}$, and thus we have  $\displaystyle{ A_\infty= \lim_{n_k \to + \infty}  \vert\vN \vu_{n_k} \vert^2 + |\vN \vb_{n_k}|^2 + \mu_{n_k}}$.   \\

Now, let  $\Phi\in \mathcal{D}((0,T)\times \mathbb{R}^3)$ be a  non-negative function. As $\sqrt\Phi (\vN \vu_{n_k}+ \vN \vb_{n_k}) $ is weakly convergent to $\sqrt \Phi (\vN\vu_\infty+ \vN\vb_\infty)$ in $L^{2}_{t}L^{2}_{x}$, we have
\begin{equation*}
\begin{split}
\iint A_\infty \Phi\, dx\, ds = & \lim_{n_k\rightarrow +\infty}\iint A_{n_k} \Phi\, dx\, ds
\geq  \limsup_{n_k \to +\infty} \iint (\vert\vN \vu_{n_k} \vert^2 + |\vN \vb_{n_k}|^2) \Phi \, dx\, ds \\ 
\geq & \iint (\vert\vN \vu_{\infty} \vert^2 + |\vN \vb_{\infty}|^2) \Phi \, dx\, ds.
\end{split}    
\end{equation*} Thus, there   exists a non-negative locally finite measure $\mu_\infty$ on $(0,T)\times\mathbb{R}^3$ such that $A_\infty=(\vert\vN \vu_\infty\vert^2+\vert\vN \vb_\infty\vert^2) +\mu_\infty$, and then we obtain the desired local energy balance: 
 \begin{equation*}
 \begin{split}
& \partial_t(\frac {\vert\vu_\infty \vert^2 
 	+ |\vb_\infty|^2 }2)= \Delta(\frac {\vert\vu_\infty\vert^2 + |\vb_\infty|^2 }2)-\vert\vN \vu_\infty \vert^2 - |\vN \vb_\infty|^2 \\
& - \vN\cdot\left( (\frac{\vert\vu_\infty \vert^2}2 + \frac{\vert\vb_\infty \vert^2}2)\vv_\infty \right)
 -\vN \cdot(p_\infty \vu_\infty) - \vN \cdot(q_\infty \vb_\infty) \\
 &+ \vN \cdot ((\vu_\infty \cdot \vb_\infty) \vc_\infty)  + \vu_\infty \cdot(\vN\cdot\mathbb{F}_{\infty}) +\vb_\infty \cdot(\vN\cdot\mathbb{G}_{_\infty})- \mu_\infty.
 \end{split}     
 \end{equation*}
In order to finish this proof, it remains to prove the convergence to the initial data $(\vu_{0,\infty},\vb_{0,\infty} )$. Once we dispose of this local energy equality,  as in \eqref{localineq} we can write:

\begin{equation*} \begin{split}
\int & \frac{\vert \vu_n (t,x)\vert^2 +  \vert \vb_n(t,x)\vert^2}2 \phi_R \, dx    +  \int_0^t  \int ( \vert\vN \vu_n \vert^2+  \vert\vN \vb_n \vert^2)\, \,  \phi_R \, dx\, ds
 \\ \leq & \int \frac{\vert \vu_{0,n}(x)\vert^2 + \vert \vb_{0,n}(x)\vert^2 }2 \phi_R \, dx   + \int_0^t  \int \frac{ | \vu_n |^2 + | \vb_n |^2}{ 2 }  \,   \Delta \phi_R \, dx\, ds
\\
 &+ \sum_{i=1}^3 \int_{0}^{t}  \int [(\frac{\vert\vu_n \vert^2}2+ \frac{|\vb_n |^2}{2} ) v_{n, i} + p_n u_{n,i} ]   \partial_i \phi_R \, dx\, ds
\\
 &+ \sum_{i=1}^3 \int_{0}^{t} \int [(\vu_n \cdot \vb_n) c_{n, i}  + q_n b_{n,i} ]   \partial_i \phi_R \, dx\, ds 
\\&- \sum_{1 \leq i,j \leq 3} ( \int_{0}^{t} \int F_{n,i,j} u_{n,j}    \partial_i \phi_R \, dx\, ds + \int_{0}^{t} \int   F_{n,i,j}\partial_i u_{n,j}\   \phi_R \, dx\, ds )
\\&- \sum_{1 \leq i,j \leq 3} ( \int_{0}^{t} \int G_{n,i,j} b_{n,j}   \partial_i \phi_R \, dx\, ds + \int_{0}^{t}  \int  G_{n,i,j}\partial_i b_{n,j}\   \phi_R \, dx\, ds ) .
\end{split}\end{equation*}
Then we have: 
 \begin{equation*}\label{} \begin{split}
\limsup_{n_k\rightarrow +\infty} \int & \frac{\vert \vu_{n_k} (t,x)\vert^2 +  \vert \vb_{n_k} (t,x)\vert^2}2 \phi_R \, dx    +  \int_0^t  \int ( \vert\vN \vu_{n_k} \vert^2+  \vert\vN \vb_{n_k} \vert^2)\, \,  \phi_R \, dx\, ds
 \\ \leq & \int \frac{\vert \vu_{0}(x)\vert^2 + \vert \vb_{0}(x)\vert^2 }2 \phi_R \, dx   + \int_0^t  \int \frac{ | \vu_\infty |^2 + | \vb_\infty |^2}{ 2 }  \,   \Delta \phi_R \, dx\, ds
\\
 &+ \sum_{i=1}^3 \int_{0}^{t}  \int [(\frac{\vert\vu_\infty \vert^2}2+ \frac{|\vb_\infty |^2}{2} ) v_{\infty, i} + p_\infty u_{\infty,i} ]   \partial_i \phi_R \, dx\, ds
\\
 &+ \sum_{i=1}^3 \int_{0}^{t} \int [(\vu_\infty \cdot \vb_\infty) c_{\infty, i}  + q_\infty b_{\infty,i} ]   \partial_i \phi_R \, dx\, ds 
\\&- \sum_{1 \leq i,j \leq 3} ( \int_{0}^{t} \int F_{\infty,i,j} u_{\infty,j}    \partial_i \phi_R \, dx\, ds + \int_{0}^{t} \int   F_{\infty,i,j}\partial_i u_{\infty,j}\   \phi_R \, dx\, ds )
\\&- \sum_{1 \leq i,j \leq 3} ( \int_{0}^{t} \int G_{\infty,i,j} b_{\infty,j}   \partial_i \phi_R \, dx\, ds + \int_{0}^{t}  \int  G_{\infty,i,j}\partial_i b_{\infty,j}\   \phi_R \, dx\, ds ) .
\end{split}\end{equation*}

Recalling that  $ \vu_{n_k}= 
\vu_{0, n_k}+ \int_0^t \partial_t \vu_{n_k}\, ds $, we may observe that $\vu_{n_k}(t,.)$ converges to $\vu_\infty(t,.)$ in $\mathcal{D}'(\mathbb{R}^3)$, hence, it converges  weakly  in $L^2_{\rm loc}(\Rt)$ and we can write:
$$\ds{\int \frac{\vert \vu_\infty(t,x)\vert^2}2 \phi_R \, dx   \leq  
\limsup_{n_k\rightarrow +\infty}\int \frac{\vert \vu_{n_k}(t,x)\vert^2}2 \phi_R \, dx}.$$     
Moreover, this  weakly convergence gives
$$
\int_0^t \int \frac{\vert \vN\vu_\infty(s,x)\vert^2}2 \phi_R \, dx \, ds  \leq  
\limsup_{n_k\rightarrow +\infty} \int_0^t\int \frac{\vert \vN\vu_{n_k}(s,x)\vert^2}2 \phi_R \, dx \, ds,$$    
and we have the same estimates for  $\vb_\infty$. In this way we get 

 \begin{equation*} \begin{split}
 \int & \frac{\vert \vu_{\infty} (t,x)\vert^2 +  \vert \vb_{\infty} (t,x)\vert^2}2 \phi_R \, dx    +  \int_0^t  \int ( \vert\vN \vu_{\infty} \vert^2+  \vert\vN \vb_{\infty} \vert^2)\, \,  \phi_R \, dx\, ds
 \\ \leq & \int \frac{\vert \vu_{0}(x)\vert^2 + \vert \vb_{0}(x)\vert^2 }2 \phi_R \, dx   + \int_0^t  \int \frac{ | \vu_\infty |^2 + | \vb_\infty |^2}{ 2 }  \,   \Delta \phi_R\, dx\, ds
\\
 &+ \sum_{i=1}^3 \int_{0}^{t}  \int [(\frac{\vert\vu_\infty \vert^2}2+ \frac{|\vb_\infty |^2}{2} ) v_{\infty, i} + p_\infty u_{\infty,i} ]   \partial_i \phi_R \, dx\, ds
\end{split}\end{equation*}
 \begin{equation*} \begin{split}
 &+ \sum_{i=1}^3 \int_{0}^{t} \int [(\vu_\infty \cdot \vb_\infty) c_{\infty, i}  + q_\infty b_{\infty,i} ]   \partial_i \phi_R \, dx\, ds 
\\&- \sum_{1 \leq i,j \leq 3} ( \int_{0}^{t} \int F_{\infty,i,j} u_{\infty,j}    \partial_i \phi_R \, dx\, ds + \int_{0}^{t} \int   F_{\infty,i,j}\partial_i u_{\infty,j}\   \phi_R \, dx\, ds )
\\&- \sum_{1 \leq i,j \leq 3} ( \int_{0}^{t} \int G_{\infty,i,j} b_{\infty,j}   \partial_i \phi_R \, dx\, ds + \int_{0}^{t}  \int  G_{\infty,i,j}\partial_i b_{\infty,j}\   \phi_R \, dx\, ds ) .
\end{split}\end{equation*}

Finally, letting $t$ go to $0$,  we have:
 $$\limsup_{t\rightarrow 0}  \| (\vu_\infty, \vb_\infty )(t,. )\|_{L^2 (\phi_R (x) dx )}^2\leq  \| (\vu_{0,\infty}, \vb_{0,\infty} )  \|_{L^2 (\phi_R (x) dx )}^2 .$$   On the other hand, by   weakly convergence we also have  
  $$    \|(\vu_{0,\infty}, \vb_{0,\infty} ) \|_{L^2 (\phi_R (x) dx )}^2 \leq   \liminf_{t\rightarrow 0}  \| (\vu_\infty, \vb_\infty ) (t,. )\|_{L^2 (\phi_R (x) dx )}^2 .$$   
Thus we have the strong convergence to initial data in the Hilbert space $L^2 (\phi_R (x) dx )$.   
\section{Proof of Theorem \ref{B_2}}\label{Proof-Th}
\subsection{Local in time existence} 
\label{led3}
Following the ideas of \cite{PF_OJ}, for the function  $\phi_R(x)=\phi(\frac x R)$ given in \eqref{varphidef}, and the Leray's projector $\mathbb{P}$, we define $\vu_{0,R}=\mathbb{P}(\phi_R \vu_0)$, $\vb_{0,R}=\mathbb{P}(\phi_R \vb_0)$, $\mathbb{F}_R=\phi_R \mathbb{F}$, $\G_R=\phi_R \G$; and we consider the approximated problem ($ MHD_{R,\epsilon}$): \begin{equation*}\label{MHD}
 \left\{ \begin{array}{ll}\vspace{2mm} 
\partial_t \vu_{R,\epsilon} = \Delta \vu_{R,\epsilon} - ( (\vu_{R,\epsilon} *\theta_\epsilon ) \cdot \nabla) \vu_{R,\epsilon} + ( (\vb_{R,\epsilon} *\theta_\epsilon ) \cdot \nabla) \vb_{R,\epsilon} - \nabla p_{R,\epsilon} + \vN \cdot \F_R,  \\ \vspace{2mm}
\partial_{t} \vb_{R,\epsilon} = \Delta \vb_{R,\epsilon} - ( ( \vu_{R,\epsilon} *\theta_\epsilon ) \cdot \nabla) \vb_{R,\epsilon} + ( ( \vb_{R,\epsilon} *\theta_\epsilon ) \cdot \nabla) \vu_{R,\epsilon} -\vN q_{R,\epsilon} + 
\vN  \cdot \G_R 
, \\ \vspace{2mm}
\vN \cdot  \vu_{R,\epsilon}=0, \,  \vN \cdot \vb_{R,\epsilon} =0, \\ \vspace{2mm}
\vu_{R,\epsilon}(0,\cdot)=\vu_{0,R}, \, \vb_{R,\epsilon}(0,\cdot)=\vb_{0,R}. 
\end{array}
\right.     
\end{equation*}
By the Appendix in \cite{PF_OJ} (see the page 35) we know that $(\text{ MHD}_{R,\epsilon})$ has a unique solution $(\vu_{R,\epsilon}, \vb_{R,\epsilon} )$  in $L^\infty((0,+\infty), L^2)\cap L^2((0,+\infty),\dot H^1)$, and moreover, this solution belongs to $\mathcal{C}([0,+\infty), L^2)$ and it fulfills the hypothesis of the Theorem  \ref{estimatesB2}.  Applying this result (for the case $(\vv, \vc) = (\vu * \theta_\epsilon, \vb * \theta_\epsilon )$)  there exists a constant $C>0$ such that for every time $T_0$ small enough:
$$ C \left( 1 +\| ( \vu_{0, R}, \vb_{0, R} ) \|_{B_2}^2+  \| ( \mathbb{F_{R, \epsilon}}, \G_{R, \epsilon} ) \|_{B_2 L^2 (0,T_0)}^2 \right)^2 \, T_0\leq 1,$$ we have the controls:
$$ \sup_{0\leq t\leq T_0} \|\ ( \vu_{R, \epsilon}, \vb_{R, \epsilon} ) (t)\|_{B_2}^2 \leq
 C \left( 1 +\|( \vu_{0, R}, \vb_{0, R} )\|_{B_2}^2+  \|( \mathbb{F_{R, \epsilon}}, \G_{R, \epsilon} )\|_{B_2 L^2 (0,T_0)}^2 \right),$$ and 
$$  {  \|\vN ( \vu_{R, \epsilon}, \vb_{R, \epsilon} ) \|_{B_2L^2 (0,T_0)}^2 }\leq  
 C  \left( 1 +\|( \vu_{0, R}, \vb_{0, R} )\|_{B_2}^2+  \|( \mathbb{F_{R, \epsilon}}, \G_{R, \epsilon} )\|_{B_2 L^2 (0,T_0)}^2 \right).$$

Then, in the setting of Theorem \ref{stability}, we set $(\vu_{0,n},\vb_{0,n})=(\vu_{0,R_n}, \vb_{0,R_n}  )$, $\mathbb{F}_n=\mathbb{F}_{R_n}$, $\G_n=\G_{R_n}$ and $(\vu_n, \vb_n )=(\vu_{R_n,\epsilon_n}, \vb_{R_n,\epsilon_n} )$; and letting  $R_n\rightarrow +\infty$ and $\epsilon_n\rightarrow 0$  we find a local solution of the (MHD) equations which verifies  the desired properties stated in Theorem \ref{B_2}.
 \subsection{Global in time existence}\label{Sec:global}
 Let $\lambda>1$. For $n\in\mathbb{N}$ we consider the (MHD) equations  with initial value $$(\vu_{0,n}, \vb_{0,n})=(\lambda^n \vu_0(\lambda^n \cdot),\lambda^n \vb_0(\lambda^n \cdot)),$$ and the forcing tensors
 $$(\mathbb{F}_n, \G_n)= (\lambda^{2n}\mathbb{F}(\lambda^{2n} \cdot , \lambda^n \cdot ),\lambda^{2n} \G( \lambda^{2n} \cdot , \lambda^n \cdot ) ).$$ Then, by the local in time existence proved above, there exists a  solution $(\vv_n, \vc_n)$ on $(0,T_n)$, with 
 $$   C \left( 1+ \|(\vv_{0,n}, \vc_{0,n}) \|_{B_2}^2+ \|(\mathbb{F}_n, \G_n)\|_{B_2 L^2 (0,T_n) }^2\right)^2\, T_n =  1.$$
Remark also that  by the well-known  scaling properties of the (MHD) equations we have $$(\vv_n(t,x),\vc_n(t,x))=(\lambda^{n}\vu_n(\lambda^{2n}t,\lambda^n x), \lambda^{n}\vb_n(\lambda^{2n}t,\lambda^n x) 
),$$ where $(\vu_n, \vb_n)$ is a solution of the   (MHD) on $(0,\lambda^{2n}T_n)$ associated  with the  initial data $(\vu_0, \vb_0)$ and then forcing tensors $\mathbb{F}$ and $ \G$.\\

At this point, we need the following simple remark which will be  proved at the end of this section. 
\begin{Remark}\label{rf} 
If $\vu_0, \vb_0 \in B_{2,0}$ and $\mathbb{F}, \G \in B_{2,0} L^2 (0, +\infty )$, then for all $\lambda > 1$ we have: 

\begin{equation*}
    \lim_{n\rightarrow +\infty} \frac{\lambda^n}{ 1+ \| (\vv_{0,n}, \vc_{0,n} ) \|_{B_2}^2+ \|(\mathbb{F}_n, \G_n )\|_{B_2 L^2}^2 }=+\infty.
\end{equation*} 
\end{Remark}
Thus,  for fixed $\lambda>1$, we have $\lim_{n\rightarrow +\infty} \lambda^{2n} T_n=+\infty$. Then, for $T>0$, let $n_T$ such that for all  $n\geq n_T$,  $\lambda^{2n}T_n>T$, then $(\vu_n, \vb_n)$ is a solution of the (MHD) equations  on $(0,T)$. \\

We set now $(\vw_n(t,x), \vd_n(t,x) ) = (\lambda^{n_T} \vu_n(\lambda^{2n_T}t,  \lambda^{n_T}x) , \lambda^{n_T} \vb_n(\lambda^{2n_T}t,  \lambda^{n_T}x) )$, where we observe that  for  $n\geq n_T$ the couple $(\vw_n, \vd_n)$ is a solution of (MHD) equations on $(0,\lambda^{-2n_T} T)$ with initial value $(\vv_{0,n_T}, \vc_{0,n_T} )$ and forcing tensor $(\mathbb{F}_{n_T}, \G_{n_T})$. But, since we have $\lambda^{-2n_T}T\leq T_{n_T}$, then we obtain 
  $$   C \left( 1+  \| (\vv_{0,n_T}, \vc_{0,n_T} )\|_{B_2}^2 + \| ( \mathbb{F}_{n_T}, \G_{n_T} ) \|_{B_2 L^2(0,\lambda^{-2n_T}T)}^2  \right)^2\, \lambda^{-2n_T} T \leq  1,$$ and thus, by Theorem \ref{estimatesB2} we are able to write:
$$ \sup_{0\leq t\leq   \lambda^{-2n_T} T} \|\ (\vw_n, \vd_n )(t,.) \|_{L^2_{w_\gamma}}^2 \leq
 C (1+  \|(\vv_{0,n_T}, \vc_{0,n_T} ) \|_{B_2}^2 + \| ( \mathbb{F}_{n_T}, \G_{n_T} ) \|_{B_2 L^2(0,\lambda^{-2n_T}T)}^2 
),$$ and 
$$  {  \|\vN(\vw_n, \vd_n)\|_{ B_2 L^2 (0, \lambda^{-2n_T}T ) }^2  }\leq  
 C (1+  \|(\vv_{0,n_T}, \vc_{0,n_T} )\|_{B_2}^2 + \|( \mathbb{F}_{n_T}, \G_{n_T} )\|_{B_2 L^2(0,\lambda^{-2n_T}T)}^2  ).$$ From these estimates we get the following uniforms controls for $\vu_n$ and $\vb_n$: 
 $$ \| (\vw_n, \vd_n) (t) \|_{B_2}^2 \geq \lambda^{n_T} \| (\vu_n, \vb_n) (\lambda^{2n_T}t,.)\|_{B_2}^2 , $$
 and
 \begin{equation*} \begin{split}   {  \|\vN (\vw_n, \vd_n) \|_{ B_2 L^2 (0, \lambda^{-2n_T}T ) }^2  }    \geq&  \lambda^{n_T}  \| \vN (\vu_n, \vb_n) \|_{B_2L^2 (0,T)}^2.
 \end{split}\end{equation*} 
 In order to finish this proof, observe that  we have controlled uniformly $\vu_n, \vb_n$ and $\vN\vu_n, \vN\vb_n$ on $(0,T)$ for $n\geq n_T$. Then,  we may  apply Theorem \ref{stability} to obtain a solution on $(0,T)$. As $T>0$ is an arbitrary time,  we can use a diagonal argument to obtain a solution $\vu, \vb$ on $(0,+\infty)$. Finally, the control for the solution $(\vu, \vb, p, q)$ on $(0, + \infty)$ is given by Theorem \ref{estimatesB2}.
 \Endproof{} \\
    
\textbf{Proof of Remark \ref{rf}.} It is enough to detail the computations for the functions $\vu_{0,n}$ and $\mathbb{F}_n$ since the computations for $\vb_{0,n}$ and $\mathbb{G}_{n}$ follows the same lines. \\

It is straightforward to see that we have
  \begin{align*}
      \frac{\|\vv_{0,n}\|_{B_2}^2}{\lambda^n} &=  \sup_{R \geq 1} \frac{1}{\lambda^n R^2} \int_{|x|\leq R} |\lambda^n \vu_0 (\lambda^n x)|^2 \, dx =  \sup_{R \geq 1} \frac{1}{ (\lambda^n R)^2} \int_{|x|\leq \lambda^n R} | \vu_0 ( x)|^2 \, dx,
 \end{align*}
 and
 \begin{equation*}
     \lim_{P \rightarrow +\infty} \sup_{R \geq P} \frac{1}{ (\lambda^n R)^2} \int_{|x|\leq \lambda^n R} | \vu_0 ( x)|^2 \, dx = \lim_{ R \rightarrow +\infty} \frac{1}{  R^2} \int_{|x|\leq R} | \vu_0 ( x)|^2 \, dx = 0.
 \end{equation*}
 Moreover, remark that we have: 
  \begin{align*}
      \frac{\|\mathbb{F}_{n}\|_{B_2 L^2 (0,+\infty)}^2}{\lambda^n} &=   \sup_{R \geq 1} \frac{1}{\lambda^n R^2} \int_0^{+\infty} \int_{|x|\leq R} |\lambda^{2n} \mathbb{F} (\lambda^{2n} t, \lambda^n x)|^2 \, dx \, ds \\
     &=   \sup_{R \geq 1} \frac{1}{ (\lambda^n R)^2} \int_0^{+\infty} \int_{|x|\leq \lambda^n R} |\mathbb{F} (t,  x)|^2 \, dx,
 \end{align*}
 and 
  \begin{equation*}
     \lim_{P \to +\infty} \sup_{R \geq P} \frac{1}{R^2} \int_0^{+\infty} \int_{|x| \leq R} |\mathbb{F}(t,x)|^2 \, dx \, ds= \lim_{R \to +\infty} \frac{1}{R^2} \int_0^{+\infty} \int_{|x| \leq R} |\mathbb{F}(t,x)|^2 \, dx \, ds=0.
 \end{equation*}
\Endproof{}\\

\begin{appendices}
\section{The 2D case}\label{Sec:2D}
In this appendix we make a discussion on the 2D case which is more delicate to treat. In comparison to the 3D case, the Leray projector is not bounded on the space  $B_2(\mathbb{R}^2)$, and for this reason one must be careful in the handling of the pressure term in dimension $2$.  We  indicate  here how to treat the pressure and moreover, we give a sketch of the proof of the local and global existence of weak suitable solutions of the (MHD) system. It is worth remark that  for the (NS) equations,  A. Basson  obtained in his Ph. D. thesis \cite{Ba06} the local existence of weak solutions for  with initial data in $B_2(\mathbb{R}^2)$. Thus, our main contribution is the study of global weak solutions in the generalized setting of the  (MHD) system which contains the (NS) as a particular case.\\ 

Let us start recalling the Basson's idea to handle the pressure term in dimension $2$.  The main point consists in giving a useful decomposition for the pressure. For this we shall fix some notation. First, we denote $\mathcal {R}$  the vector field  of the Riesz transforms and we shall write $H_{i,j}$  the kernel of the operator $\mathcal{R}_i  \mathcal{R}_j $ and we denote $\mathbb{H}= (H_{i,j})$. On the other hand, let $\varphi \in \mathcal{D}(\mathbb{R}^2)$ be  a non negative function supported on $B(0,2)$ such that $\varphi = 1$ on $B(0,1)$.  For the function $\varphi$ given, and for $k\in \mathbb{N}$, we define the functions
$\psi_k (x) = \varphi (2^{-k-1} x  ) - \varphi (2^{-k} x  )$ and $\chi_k  =  \varphi(2^{-k-3} x  ) - \varphi (2^{-k+2} x)$. Hence we have,  $\psi_k(x) = 1$ for all $ 2^{k-1} \leq |x|\leq  2^{k+3}$, and moreover, we have
$$\text{sup}(\chi_k) \subset \{x \in \mathbb{R}^2 : \,  2^{k-2}  \leq |x|\leq2^{k+4} \},$$
and
$$\text{sup}(\psi_k) \subset \{x \in \mathbb{R}^2 : \,  2^{k} \leq |x|\leq  2^{k+2} \}.$$

Then, for a index-family $\mathcal{A}$ that we shall set conveniently later,   let $(\vu_\alpha)_{\alpha \in \mathcal{A}}$ and $(\vb_\alpha)_{\alpha \in \mathcal{A}}$ be two families of time dependent vector fields defined on $[0, T) \times \mathbb{R}^2$, and let $(\mathbb{F}_\alpha)_{\alpha \in A}$, $(\mathbb{G}_{\alpha})_ {\alpha \in \mathcal{A}}$  be  two families of tensors 
defined on $[0, T) \times \mathbb{R}^2$. For $\alpha \in \mathcal{A}$ we denote
\begin{equation}\label{defA}
A_{\alpha} =  \vu_\alpha \otimes \vu_\alpha - \vb_\alpha \otimes \vb_\alpha - \mathbb{F}_\alpha,
\end{equation}
and then  we define the terms $p_{\alpha,1}$, $\nabla p_{\alpha,2}$, $q_{\alpha,1}$ and $\nabla q_{\alpha,2}$ as follows:
\begin{equation}
\label{defp1}
    p_{\alpha,1} =  \varphi (x/8)  \mathcal{R} \otimes \mathcal{R} ( \varphi ( A_{\alpha} ) ) +  \sum_{k=1}^{+\infty} \chi_k  \mathcal{R} \otimes \mathcal{R} ( \psi_k (A_{\alpha}  ) ),
\end{equation}
\begin{equation}
\label{defp2}
\begin{split}
    {\vN} p_{\alpha,2} =&   \vN [(1- \varphi (x/8) )  \mathcal{R} \otimes \mathcal{R} ( \varphi (A_{\alpha})  )]  \\
    &+  \sum_{k=1}^{+\infty} \vN [(1- \chi_k )  \mathcal{R} \otimes \mathcal{R} ( \psi_k (A_{\alpha})  ) ] .
\end{split}
\end{equation}
and 
\begin{equation}\label{defq1}
    q_{\alpha,1} =  \varphi (x/8)  \mathcal{R} \otimes \mathcal{R} ( \varphi ( - \mathbb{G}_\alpha) ) +  \sum_{k=1}^{+\infty} \chi_k  \mathcal{R} \otimes \mathcal{R} ( \psi_k ( - \mathbb{G}_\alpha  ) ),
\end{equation}
\begin{equation}\label{defq2}
\begin{split}
     {\vN} q_{\alpha,2} =&   \vN [(1- \varphi (x/8) )  \mathcal{R} \otimes \mathcal{R} ( \varphi ( - \mathbb{G}_\alpha)  )]  \\
    &+  \sum_{k=1}^{+\infty} \vN [(1- \chi_k )  \mathcal{R} \otimes \mathcal{R} ( \psi_k ( - \mathbb{G}_\alpha)  ) ].
    \end{split}
\end{equation}
Now, let us explain the general idea  to study the local and global existence of solutions in the 2D case. First, we will   consider a (MHD)-type system of approximated  solutions $(\vu_n, \vb_n, p_n, q_n)$ where the key point is to split the terms $p_n$ and $q_n$ as  the expressions above.  Thereafter, using a local energy balance we will obtain a uniform bound on the approximated solutions. Finally, passing to the limit we will able to  get a solution $(\vu, \vb, p,q)$ of the (MHD) system. Following these ideas we have the next result. 

\begin{Theorem}[Local and global weak suitable solutions]\label{dim2} Let $0<T <+\infty$. Let  $\vu_{0}, \vb_{0} \in B_{2}(\mathbb{R}^2)$ be  divergence-free vector fields. Let    $\mathbb{F}$ and $\G$ be    tensors  belonging to $ B_2 L^2 (0,T)$.  Then, there exists a time $0<T_0 <  T$ such that  the system (MHD) has a solution $(\vu, \vb,p,q)$  which satisfies :
 \begin{itemize} 
 \item[$\bullet$]  $\vu, \vb$ belong to $L^\infty((0,T_0), B_2   )$ and $\vN\vu, \vN\vb $ belong to $B_2   L^2 (0,T_0)$.
\item[$\bullet$] The pressure $p$ and the term $q$ are related to $\vu, \vb$, $\mathbb{F}$ and $\G$ as follows. Let  $\tilde \varphi \in \mathcal{D}(\mathbb{R}^2)$  be a test function such that $\tilde\varphi(x )=1$ on a neighborhood of the origin. We define  $$\Phi_{i,j,\varphi}=(1-\tilde \varphi) \partial_i \partial_j G_2. $$ where $G_2 = \frac 1{2\pi}\ln(\frac 1{\vert x\vert})$ is the fundamental solution of the operator $-\Delta$ (we have  $-\Delta G_2=\delta_0$). Then $p$ and $q$ can be defined by : 

\begin{equation}\label{caracterizacion-p}
\begin{split}
 p_{\tilde \varphi}(t,x)=  & \sum_{i,j}(\tilde \varphi \partial_i \partial_j G_2) * (u_iu_j-b_ib_j-F_{i,j})(t,x)\cr &+   \sum_{i,j}   \int (\Phi_{i,j, \tilde\varphi}(x-y)-\Phi_{i,j,\tilde\varphi}(-y))  \left(u_i(t,y)u_j(t,y)\right. \cr & \left.-b_i(t,y)b_j(t,y)-F_{i,j}(t,y)\right)\, dy, 
 \end{split}
 \end{equation}
and
\begin{equation}\label{caracterizacion-q}
\begin{split}
 &q_{\tilde\varphi}(t,x)=   \sum_{i,j}(\tilde\varphi \partial_i \partial_j G_2) * (-G_{i,j})(t,x) \cr &+   \sum_{i,j}   \int (\Phi_{i,j,\tilde\varphi}(x-y)-\Phi_{i,j,\tilde\varphi}(-y))  (-G_{i,j}(t,y))\, dy.
 \end{split}
  \end{equation}
 
 \item[$\bullet$] The map $t\in [0,T)\mapsto (\vu(t,\cdot ), \vu(t,\cdot ) )$ is $*$-weakly continuous from $[0,T)$ to $B_2(\mathbb{R}^2) $, and for all compact set $K \subset \mathbb{R}^2$ we have:
 $$ \lim_{t\rightarrow 0}    \Vert  ( \vu(t,\cdot) - \vu_{0} , \vb(t,\cdot) - \vb_{0} ) \Vert_{L^{2}(K)} = 0.$$
 \item[$\bullet$] The solution $(\vu, \vb, p, q)$ is suitable : there exists a non-negative locally finite measure $\mu$ on $(0,T)\times\mathbb{R}^2$ such that:
\begin{equation*}
\begin{split}
   \partial_t(\frac {\vert\vu \vert^2 
    + |\vb|^2 }2)=&\Delta(\frac {\vert\vu \vert^2 + |\vb|^2 }2)-\vert\vN \vu \vert^2 - |\vN \vb|^2 - \vN\cdot\left( [\frac{\vert\vu \vert^2}2 + \frac{\vert\vb \vert^2}2 +p]\vu   \right) \\
  &  + \vN \cdot ([(\vu \cdot \vb)+ q] \vb ) + \vu \cdot(\vN\cdot\mathbb{F}) +\vb \cdot(\vN\cdot\mathbb{G})- \mu.
\end{split}     
\end{equation*}
 \end{itemize}
In particular we have the global control on the solution: for all $0\leq t \leq T_0$, 
 \begin{equation} \label{global-control-2D}
 \begin{split}  
 &\max \{ \|(\vu, \vb )(t) \|_{B_2}^2 , \|\vN(\vu, \vb) \|_{B_2 L^2 (0,T_0)}^2  \} \\
 &\leq  C\| (\vu _0, \vb _0 ) \|^2_{B_{2}} + C\|  (\mathbb{F},\G) \|^2_{B_2 L^2 (0,t)}+ C \int_0^t 1 + \| (\vu, \vb) (s) \|^4_{B_{2 }} ds. \end{split}
\end{equation}
\item[$\bullet$] Finally,  if the data verify: $$
    \lim_{R \to +\infty}  R^{-2} \int_{|x| \leq R} |\vu_0 (x)|^2 + |\vb_0 (x)|^2 \, dx =0,
$$
and  
$$
    \lim_{R \to +\infty}  R^{-2} \int_0^{+\infty} \int_{|x| \leq R} |\mathbb{F}(t,x)|^2 + |\mathbb{G}(t,x)|^2 \, dx \, dt=0,
$$ then $(\vu, \vb, p,q)$ is a global weak solution.\\
 \end{Theorem}

\textbf{Sketch of the proof}. This proof follows some of the main ideas of \cite{Ba06} and thus we will only detail  its main steps. The key is to obtain \emph{a priori} estimates for the following  approximated solutions. Let $\phi_R$ as in \eqref{varphidef}. Then, in \cite{Ba06}, Basson proves that we can take a sequence $R_n\rightarrow +\infty$ such that $\mathbb{P}(\phi_{R_n} \vu_0)$ converges *-weakly to $\vu_0$ and $\mathbb{P}(\phi_{R_n} \vb_0)$ converges *-weakly to $\vb_0$ in $B_2$. As in section \ref{led3},  we define $\vu_{0,n}=\mathbb{P}(\phi_{R_n} \vu_0)$, $\vb_{0,n}=\mathbb{P}(\phi_{R_n} \vb_0)$, $\mathbb{F}_n=\phi_{R_n} \mathbb{F}$, $\G_n=\phi_{R_n} \G$; and we consider the solutions of the  approximated problem ($ MHD_{n}$): \begin{equation*}\label{MHDaprox}
 \left\{ \begin{array}{ll}\vspace{2mm} 
\partial_t \vu_n = \Delta \vu_n - (\vu_n \cdot \nabla) \vu_n + (\vb_n \cdot \nabla) \vb_n - \nabla p_n + \vN \cdot \F_n,  \\ \vspace{2mm}
\partial_{t} \vb_n = \Delta \vb_n - (\vu_n \cdot \nabla) \vb_n + (  \vb_n  \cdot \nabla) \vu_n -\vN q_n + 
\vN  \cdot \G_n 
, \\ \vspace{2mm}
\vN \cdot  \vu_n=0, \,  \vN \cdot \vb_n =0, \\ \vspace{2mm}
\vu_n(0,\cdot)=\vu_{0,n}, \, \vb_n(0,\cdot)=\vb_{0,n}, 
\end{array}
\right.     
\end{equation*}
which belongs to $L^\infty((0,+\infty), L^2(\mathbb{R}^2))\cap L^2((0,+\infty),\dot H^1(\mathbb{R}^2))$, and moreover, belongs to $\mathcal{C}([0,+\infty), L^2(\mathbb{R}^2))$.
Now we write $\nabla p_n= \nabla p_{n,1}+\nabla p_{n,2}$, where the term $p_{n,1}$ is given by (\ref{defp1}) and the term $\nabla p_{n,2}$ is defined in (\ref{defp2}). Similarly, using the expressions (\ref{defq1}) and (\ref{defq2}) we write $\nabla q_n= \nabla q_{n,1}+\nabla q_{n,2}$. 
We need the following technical lemma. 
\begin{Lemma}
	\label{baspgcontrol}
	Let $\F=(F_{i,j})_{1\leq i,j \leq 2} \in L^1_{loc}$ be a tensor. Then we have
		\begin{align*}
	\label{press}
	\| \vN [ (1-\chi_k ) \mathcal{R}_i \mathcal{R}_j ( \psi_k F_{i,j} ) ] \|_{L^{\infty} (\mathbb{R}^2)} &\leq  C \int \frac{ \psi_k (y) |\F (y)|  }{(1+|y|)^3} \, dy \nonumber\\
	& \approx C 2^{-3k} \int  \psi_k (y) |\F (y)| \, dy.
	\end{align*}
	and 
	\begin{equation*}
	\| \vN [ ( 1-\varphi (x/8)) \mathcal{R}_i \mathcal{R}_j ( \varphi F_{i,j} ) ] \|_{L^{\infty} (\mathbb{R}^2)} \leq  C \int  \varphi (y) |\F (y)|  \, dy \nonumber
	\end{equation*}
\end{Lemma}
\pv
We will proceed as Proposition 1.2 in \cite{Basson06}.
In dimension $2$, the kernel $\mathbb{H}$ satisfies $|\mathbb{H}(x)| \leq \frac{C}{|x|^2}$ and $| \vN \mathbb{H}(x) | \leq \frac{C}{|x|^3}$. Remark that, for all $k \in \mathbb{N}$, by the localization properties of the functions $\psi_k$ and $\chi_k$,  for  $|x-y| < \frac{1+|y|}{16} $ we have   $(1-\chi_k )(x) \psi_k (y)=0$ and $\vN \chi_k(x) \psi_k (y)=0$. Then we can write 
\begin{equation*}
    \mathcal{R}_i \mathcal{R}_j ( \psi _k F_{i,j}) (x) = \int_{|x-y|>\frac{1+ |y|}{16}} H_{i,j} (x-y) \psi_k(y) F_{i,j} dy,
\end{equation*}
and moreover
\begin{align*}
|  \vN & [ (1-\chi_k ) \mathcal{R}_i \mathcal{R}_j ( \psi _k F_{i,j} ) ] (x)| \\
\leq &  | \vN (1-\chi_k )(x) \int_{|x-y|>\frac{1+ |y|}{16}} H_{i,j} (x-y) \psi_k(y) F_{i,j} dy \, | \\
& + | (1-\chi_k )(x) \int_{|x-y|>\frac{1+ |y|}{16}} \vN H_{i,j} (x-y) \psi_k(y) F_{i,j} dy \, | \\
\leq & C 2^{-k} \int \frac{ \psi_k (y) |\F (y)|  }{(1+|y|)^2} \, dy + C \int \frac{ \psi_k (y) |\F (y)|  }{(1+|y|)^3} \, dy \\
\approx &  C \int \frac{ \psi_k (y) |\F (y)|  }{(1+|y|)^3} \, dy\approx C 2^{-3k} \int  \psi_k (y) |\F (y)| \, dy.
\end{align*}
On the other hand,  by the localization properties of the function $\varphi$ we have  $(1-\varphi(x/8) ) \varphi (y)=0$ and $\vN \varphi(x/8) \varphi (y)=0$; and the second estimate stated in this lemma follows the same lines.\Endproof{} \\

With this lemma at hand, we get back to the expression (\ref{defA}) (where we set $\alpha=(n,i,j)\in \mathbb{N}\times \{ 1,2\}\times \{1,2\}$) to obtain 
\begin{align*} 
    & \| \vN [ (1-\chi_k ) \mathcal{R}_i \mathcal{R}_j  ( \psi_k A_{n,i,j}) ] \|_{L^{\infty} (\mathbb{R}^2)} \\
    &\leq C   \int \frac{ \psi_k  ( |\vu_n |^2+  |\vb_n |^2 ) }{(1+|y|)^3} \, dy  + C  \int \frac{ \psi_k  | \mathbb{F}_n |  }{(1+|y|)^3} \, dy  \\
    & \approx 2^{-3 k} \int \psi_k ( |\vu_n |^2+  |\vb_n |^2 ) dy + 2^{-3k} \int \psi_k  | \mathbb{F}_n | dy.
\end{align*}
Then, after summation over $k$  and using H\"older inequality in the term with the forcing tensor,  for $2<\gamma_0<4$ we obtain 
\begin{equation}
\begin{split}
        \|  \vN p_{n,2} \|_{L^{\infty}} &<  C \int \frac{|\vu_n|^2 + |\vb_n|^2 }{(1+|x|)^3} \, dy + C \left( \int \frac{|\mathbb{F}_n|^2}{(1+|x|)^{\gamma_0} } \right)^{1/2},
\end{split}
\label{p2infty}
\end{equation}
hence $ \|  \vN p_{n,2} \|_{L^{\infty}} $ is uniformly bounded. The  same statement holds for $q_{n,2}$.\\

We study now the terms $p_{n,1}$ and $q_{n,1}$. First, let $R \geq 1$ fix, and we set  $k_0 \in \mathbb{N}$ such that  $2^{k_0 -1} \leq 2R \leq 2^{k_0}$. Thereafter,  for the expression $A_n$ given in (\ref{defA}), by the localization properties of the function $\chi_k$  we have: 

\begin{equation}
\label{desig}
\begin{split}
   \int_{|x| \leq 2R} |p_{n,1}|^{3/2}\, dx   \leq & C \int_{|x| \leq 2R} |\varphi (x/8)  \Ri \otimes \Ri ( \varphi (A_n ) )|^{3/2} \, dx \\
  &  + C \sum_{k=1}^{k_0 +1} \int_{|x| \leq 2R}  | \chi_k  \Ri \otimes \Ri ( \psi_k (A_n ) ) |^{3/2}\, dx  \\
    \leq &  C \int_{|x| \leq 2R} | \varphi ( \vu_n \otimes \vu_n + \vb_n \otimes \vb_n )  |^{3/2} dx \\
    & + C \sum_{k=1}^{k_0 +1 } \int_{|x| \leq 2R}  |  \psi_k (\vu_n \otimes \vu_n + \vb_n \otimes \vb_n )  |^{3/2} \, dx\\
        & + C \int | \varphi \F_n |^{3/2} + C \sum_{k=1}^{k_0 +1 } \int  |  \psi_k \F_n  |^{3/2} \\
    \leq &  C \int_{|x| \leq 2^5 R} |  \vu_n   |^{3} + |  \vb_n   |^{3} dx + C \int_{|x| \leq 2^5 R} |  \F_n   |^{3/2} dx.
\end{split}
\end{equation}
Similarly, for the term $q_{n,1}$ we have the estimate:
\begin{equation*}
    \int_{|x| \leq 2R} |q_{n,1}|^{3/2} \,dx 
    \leq    C \int_{|x| \leq 2^5 R} |  \G_n   |^{3/2}\,dx.
\end{equation*}
Once we have all these estimates, we are able to establish the  global energy controls (\ref{global-control-2D}) for the approximated solutions $(\vu_n, \vb_n, p_n, q_n)$, where, for the sake of simplicity, we will get rid of the index $n$ and we shall write $(\vu, \vb, p,q)$. \\

In this part of the proof we will proceed as in dimension $3$, so we  only detail the main computations.  Let $0< t_0 < t_1 <T$, we define $\alpha_{\eta,t_0,t_1}(t)$ as in \eqref{alpha} so that $\alpha_{\eta,t_0,t_1}$ converges almost everywhere to $\mathds 1 _{[t_0, t_1]}$ when $\eta \to 0$  and $\partial_t \alpha_{\eta,t_0,t_1}$ is the difference between two identity approximations, the first one in $t_0$ and the second one in $t_1$. Moreover, we always define $\phi_R$ as in \eqref{varphidef}. Thus, we use the local energy balance to write the following estimate:

 \begin{equation*} \begin{split}
 \iint   &\partial_t (\frac{| \vu |^2 +  | \vb |^2}2)  \alpha_{\eta,t_0,t_1} \phi_R \, dx\, ds      +    \iint  |\vN\vu|^2 + |\vN\vb|^2 \, \,  \alpha_{\eta,t_0,t_1} \phi_R \,dx\,ds
 \\ &\leq  \iint \frac{ | \vu |^2 + | \vb |^2  }{ 2 }  \,   \alpha_{\eta,t_0,t_1} \Delta \phi_R \, dx \, ds 
\\&  + \sum_{i=1}^2  \iint   (\frac{|\vu|^2 + |\vb|^2  }2 )   u_i      \alpha_{\eta,t_0,t_1} \partial_i \phi_R\,  dx \, ds + \iint \vu \cdot \vN p \phi_R\,  dx \, ds 
\\&  + \sum_{i=1}^2  \iint  ( \vu \cdot \vb   )b_i      \alpha_{\eta,t_0,t_1}\partial_i \phi_R \, dx \,ds + \iint \vb \cdot \vN q \phi_R\,  dx\, ds 
\\&- \sum_{1 \leq  i,j \leq 2}  \iint    F_{i,j} u_j    \alpha_{\eta,t_0,t_1} \partial_i \phi_R \,dx\,ds  - \sum_{1 \leq  i,j \leq 2}  \iint    F_{i,j}\partial_iu_j  \alpha_{\eta,t_0,t_1} \phi_R  \,dx\,ds
\\&- \sum_{1 \leq  i,j \leq 2}  \iint    G_{i,j} b_j    \alpha_{\eta,t_0,t_1} \partial_i \phi_R \,dx\,ds  - \sum_{1 \leq  i,j \leq 2}  \iint    G_{i,j}\partial_ib_j\   \alpha_{\eta,t_0,t_1} \phi_R \,dx\,ds. 
\end{split}\end{equation*}
In each term in  this estimate, first we  divide by $R^2$. Then we use the estimates $\Vert \nabla \phi_R \Vert_{L^{\infty}} \leq c/ R$ and $\Vert \Delta \phi_R \Vert_{L^{\infty}} \leq c / R^2$ (where we recall that $R \geq 1$). Finally, for the terms involving the tensor $\F$ and $\G$, we apply the Cauchy-Schwarz inequalities and the Young inequalities with two constants $C_0>0$ and $C=C(C_0)>0$, where the first constant $C_0$ is arbitrarily small. Thus we get: 

 \begin{equation}\label{localineq}
 \begin{split} 
      &\frac{1}{R^2}  \iint  \partial_t (\frac{| \vu |^2 + |\vb|^2 }2   ) \alpha_{\eta,t_0,t_1}\phi_R\,dx\,ds     +  \frac{1}{R^2}  \iint  (|\vN\vu|^2 + |\vN \vb|^2)\, \,  \alpha_{\eta,t_0,t_1} \phi_R \,dx\,ds\\ 
     \leq & \frac{C}{R^2}  \iint_{|x|\leq 2R} (| \vu |^2  + | \vb |^2) \alpha_{\eta,t_0,t_1}\,dx\,ds  + \frac{1}{R^3} \iint_{|x|\leq 2R}  (\frac{|\vu|^3 + |\vb|^2 |\vu|}2) \alpha_{\eta,t_0,t_1} \,dx\,ds \\
     &  + \underbrace{  \frac{2}{R^2} \iint   \phi_R  (\vu \cdot \vN p + \vb \cdot  \vN q)  \alpha_{\eta,t_0,t_1} \,dx\,ds}_{(a)}\\
     & + \frac{C}{R^2} \iint_{|x|\leq 2R} |\F |^2 \alpha_{\eta,t_0,t_1}\,dx\,ds + \frac{C}{R^2} \iint_{|x|\leq 2R} |\G |^2 \alpha_{\eta,t_0,t_1}\,dx\,ds\\
     &+ \frac{C_0}{R^2}  \iint_{|x|\leq 2R}  (|\vN\vu|^2 + |\vN\vb|^2)\, \alpha_{\eta,t_0,t_1} \,dx\,ds.
\end{split}\end{equation}
Now, we need to estimate the term $(a)$. We will only treat   the term  $\vu \cdot \nabla p$ since the other term $\vb \cdot \nabla q$ follows the same estimates.  To control this term we use the decomposition $\nabla p= \nabla p_1+ \nabla p_2$ where $p_1$ and $\nabla p_2$ are always given by (\ref{defp1}) and (\ref{defp2}) respectively. First, as we have $div(\vu)=0$ then we write
\begin{equation*}
    \int   \phi_R  \vu \cdot \vN p \, dx = - \int p_{1} \vu    \cdot \vN \phi_R \, dx  + \int   \phi_R  \vu \cdot  \vN p_2 \, dx.
\end{equation*}
On the other hand, we will need the following estimate:
by the interpolation inequalities we have 
\begin{align}
\label{u3b22}
     \int_{|x| \leq 2^5 R} |  \vu  |^{3} dx& \leq  \| \varphi_{2^5 R} \vu_n \|_{L^3}^3 \nonumber \leq  C \| \varphi_{2^5 R} \vu \|_{L^2}^2 \| \vN (\varphi_{2^5 R} \vu) \|_{L^2} \nonumber  \\
    & \leq  C \| \varphi_{2^5 R} \vu \|_{L^2}^3  +C' \| \varphi_{2^5 R} \vu \|_{L^2}^2\|  \varphi_{2^5 R} \vN \vu \|_{L^2}. 
\end{align}
Thus, getting back to the previous identity, using the estimates  \eqref{p2infty} and \eqref{desig}, and moreover, using the estimate  \eqref{u3b22}, for $2<\gamma_0<4$  and $2 < \gamma_1 < 10/3$,  we obtain
\begin{equation*}
\begin{split}
        \frac{1}{R^2} & \left| \int   \phi_R  \vu \cdot \vN p \, dx \right|  \leq \frac{C}{R^3} \int_{\vert x \vert \leq 2R} (|\vu|^3 + |p_1|^{3/2}) dx + C\frac{ \|  \vN p_2 \|_{\infty}}{R} \left(\int \phi_R |\vu|^2\,dx\right)^{1/2} \\
         \leq & \frac{C}{R^3} \int_{|x| \leq 2^5 R} (|  \vu   |^{3}+ |  \vb   |^{3})dx  + C \int \frac{|  \F   |^{3/2}}{ (1+|x|)^3} dx \\
         &+ \frac{C}{R} (\|\vu\|^2_{B_2} + \|\vb\|^2_{B_2} )  \| \phi_R \vu \|_{L^2}  + \frac{C}{R^2} \| \phi_R \vu \|^2_{L^2} + C  \int \frac{|\mathbb{F}|^2}{(1+|x|)^{\gamma_0} } dx  \\
         \leq &  \frac{C}{R^3}  \| \varphi_{2^5 R} \vu \|_{L^2}^3  +\frac{C'}{R^3}  \| \varphi_{2^5 R} \vu \|_{L^2}^2\|  \varphi_{2^5 R} \vN \vu \|_{L^2}  + C \left( \int \frac{|  \F  |^{2}}{ (1+|x|)^{\gamma_1}}\, dx  \right)^{3/4}  \\
         & +  \frac{C}{R^3}  \| \varphi_{2^5 R} \vb \|_{L^2}^3   +\frac{C'}{R^3}  \| \varphi_{2^5 R} \vb \|_{L^2}^2\|  \varphi_{2^5 R} \vN \vb \|_{L^2}  \\
         &  + \frac{C}{R} (\|\vu\|^2_{B_2} + \|\vb\|^2_{B_2} ) \left(\int_{|x|\leq 2R}  | \vu |^2  dx\right)^{1/2}+ \frac{C}{R^2} \| \phi_R \vu \|^2_{L^2} dx  + C \int \frac{|\mathbb{F}|^2}{(1+|x|)^{\gamma_0} } dx.
\end{split}
\end{equation*}
Thus, using this control on the term $(a)$, we get back to estimate (\ref{localineq}) and we obtain:

     \begin{equation*} \begin{split}  
     \frac{1}{R^2} & \iint  \partial_t (\frac{| \vu |^2 + | \vb |^2 }2 )  \alpha_{\eta,t_0,t_1} \phi_R  \,dx\,ds      +  \frac{1}{R^2}  \iint  (|\vN\vu|^2 + |\vN\vb|^2) \,  \alpha_{\eta,t_0,t_1} \phi_R \,dx\,ds   \\
     \leq &  \frac{C}{R^2}  \iint_{|x|\leq 2R} (| \vu |^2 + | \vb |^2) \alpha_{\eta,t_0,t_1} \,dx\,ds \\
     & + \frac{C}{R^3} (\int \| \varphi_{2^5 R} \vu \|_{L^2}^3 + \| \varphi_{2^5 R} \vb \|_{L^2}^3) \alpha_{\eta,t_0,t_1} \,ds + C \int (\|\vu\|_{B_2}^3 +  \|\vb\|_{B_2}^3 ) \alpha_{\eta,t_0,t_1} \,ds  \\
     &    + C  \int \frac{1}{R^4} \|   \varphi_{2^5 R} \vu \|_{L^2}^4 \alpha_{\eta,t_0,t_1} \, ds +  \frac{C
     _0}{R^2} \int  \|  \varphi_{2^5 R} \vN \vu \|^2_{L^2} \alpha_{\eta,t_0,t_1} \, ds   \\
     &    + C  \int \frac{1}{R^4} \|   \varphi_{2^5 R} \vb \|_{L^2}^4 \alpha_{\eta,t_0,t_1} \, ds  + \frac{C_0}{R^2} \int  \|  \varphi_{2^5 R} \vN \vb \|^2_{L^2} \alpha_{\eta,t_0,t_1}\, ds \\
     & + C \int \left( 1 + \int \frac{| \F |^2 + | \G |^2 }{ (1+|x|)^{\gamma_1}} \alpha_{\eta,t_0,t_1} dx\right) ds + C \int   \int \frac{|\F|^2 + |\G|^2 }{(1+|x|)^{\gamma_0} } \alpha_{\eta,t_0,t_1} \,dx\,ds \\
     & + \frac{C}{R^2} \iint_{|x|\leq 2R}(|\F |^2 + |\G |^2) \alpha_{\eta,t_0,t_1} \,ds\,ds + \frac{C_0}{R^2}  \iint_{|x|\leq 2R}  (|\vN\vu|^2 + |\vN\vb|^2 )\,  \alpha_{\eta,t_0,t_1} \,dx\,ds. 
     \end{split}
     \end{equation*}
Now, we let $\eta $ goes to $0$, and moreover, since for $t_0,t_1$ two Lebesgue points of the function  $ A_{R}(s)=\ds{\int | \vu(s,x)|^2 + | \vb(s,x)|^2   \phi_R(x)  \, dx}$, we have 
 $$-\iint \frac{|\vu|^2+|\vb|^2}2 \partial_t\alpha_{\eta,t_0,t_1}   \phi_R  \, dx\, ds=-\frac 1 2 \int \partial_t \alpha_{\eta,t_0,t_1} 
A_{R}(s) \, ds, $$
and
 $$ \lim_{\eta\rightarrow 0}   -\iint \frac{|\vu|^2+|\vb|^2}2 \partial_t\alpha_{\eta,t_0,t_1}  \phi_R\, dx\, ds=\frac 1 2 (  A_{R}(t_1)- A_{R}(t_0)),$$ 
 then we get 
     \begin{equation*} \begin{split}  
     \frac{1}{R^2}  & \int   ( \frac{| \vu(t_1) |^2 + | \vb(t_1) |^2 }2 - \frac{|\vu(t_0)|^2 + |\vb(t_0)|^2}{2})  \phi_R\, dx\\
     &  +  \frac{1}{R^2}  \int_{t_0}^{t_1} \int   (|\vN\vu|^2 +|\vN \vb|^2  \,   \phi_R) \,dx\,ds  \\
     \leq &   C \int_{t_0}^{t_1} 1 +    \|\vu\|_{B_2}^4 + \|\vb\|_{B_2}^4 \,ds   \\
     & + C \|  \mathbb{F} \|^2_{B_2 L^2 (t_0,t_1)} + C \|  \mathbb{G} \|^2_{B_2 L^2 (t_0,t_1)} \\
     &  + \frac{C_0}{R^2}  \int_{t_0}^{t_1} \int_{|x|\leq 2^5 R}  |\vN\vu|^2 + |\vN\vb|^2 \, dx\,ds. \end{split}\end{equation*}

Thereafter, by the continuity at $0$ of the map $t \in [0,T) \mapsto (\vu, \vb)(t) \in L^2_{loc}(\mathbb{R}^2)$, we can let $t_0$ go to zero. Moreover,  by the *-weak  continuity of this map  we can let $t_1$ go to $t \in (0,T)$. Thus, for all $t \in (0,T)$ we find 

    \begin{equation*} 
    \begin{split}   
    \frac{1}{R^2}  \int &  ( \frac{| \vu(t) |^2+| \vb(t) |^2}2)  \phi_R \, dx     +  \frac{1}{R^2}  \int_0^t \int  ( |\vN\vu|^2 + |\vN\vb|^2)\,   \phi_R\, dx   \\
    \leq &  C (\| \vu_0 \|_{B_2}^2 +  \| \vb_0 \|_{B_2}^2)  + C \|  \mathbb{F} \|^2_{B_2 L^2 (0,t)} + C \|  \mathbb{G} \|^2_{B_2 L^2 (0,t)} \\
    & +   C \int_{t_0}^{t_1} 1 +    \|\vu\|_{B_2}^4 + \|\vb\|_{B_2}^4 \,ds   \\
    &  + \frac{C_0}{R^2}  \int_0^t \int_{|x|\leq 2^5 R}  (|\vN\vu|^2 +  |\vN\vb|^2 )\,dx\,ds . \end{split}\end{equation*}
In this estimate, we use first the Young inequalities and moreover, we take the supreme on $R\geq 1$ to obtain  the global energy control (\ref{global-control-2D}) for the approximated solutions. Once we have this global energy control, 
following the ideas of A. Basson in \cite{Ba06}, and using a slightly modification of Lemma 3.1 in \cite{PF_OJ},  we will able to obtain a subsequence $(\vu_{n_k}, \vb_{n_k}, p_{n_k}, q_{n_k}) $ which converges in the sense of distributions to a local solution $(\vu, \vb, p, q)$ of the (MHD) system on $[0,T_0]$, where  $$T_0 \approx \frac{1}{1+ \|(\vu_0,\vb_0 ) \|_{B_2}^2 + \|(\F,\G ) \|_{B_2 L^2 (0, + \infty)}  ^2}.$$

The pressure terms are given by the expressions $\vN p = \vN p_1 + \vN p_2$, with  $\ds{p_1 = \lim_{k \to \infty } p_{n_k ,1} }$ and $\ds{\vN p_2 = \lim_{k \to \infty } \vN p_{n_k ,2}}$; and $\vN q = \vN q_1 + \vN q_2$, with  $\ds{q_1 = \lim_{k \to \infty } q_{n_k ,1} }$ and $\ds{\vN q_2 = \lim_{k \to \infty } \vN q_{n_k ,2}}$. Moreover, $p_1$ and $\nabla p_2$ satisfy \eqref{defp1} and \eqref{defp2} (with $\vu_\alpha = \vu$) respectively. Similarly, the terms $q_1$ and $\nabla q_2$ satisfy (\ref{defq1}) and (\ref{defq2}). Thus, by Proposition \ref{prp} (with $d=2$) we get that  $p$ can be written as  (\ref{caracterizacion-p}) and $q$ can be written as  (\ref{caracterizacion-q}). Finally, rescaling the local solution, we can obtain a global solution in the same way as in dimension 3. See Section \ref{Sec:global} for the details. \Endproof{}\\
\end{appendices}

\textbf{Acknowledgement.} We are very grateful to Pierre Gilles Lemarié-Rieusset for his encouragement in the study of the MHD equations through these interesting new methods developed for the Navier-Stokes equations. We are also grateful for his useful comments to get over some technical difficulties. On the other hand, we are grateful to the referee for the valuable comments to improve this paper.

\end{document}